\documentclass[a4paper,reqno]{amsart}

\usepackage{tgschola}

\usepackage{latexsym}
\usepackage[english]{babel}
\usepackage{fancyhdr}
\usepackage[mathscr]{eucal}
\usepackage{amsmath}
\usepackage{mathrsfs}
\usepackage{mathtools}
\usepackage{amsthm}
\usepackage{amsfonts}
\usepackage{amssymb}
\usepackage{amscd}
\usepackage{bbm}
\usepackage{graphicx}
\usepackage{graphics}
\usepackage{latexsym}
\usepackage{color}
\usepackage{pifont}
\usepackage{tikz}
\usepackage[normalem]{ulem}
\usepackage{changepage}

\usepackage{multirow}

\usepackage{geometry}

\newcommand{\cH}{\mathcal{H}}

\theoremstyle{plain}

\theoremstyle{definition}

\newtheorem*{remark*}{Remark}

\numberwithin{equation}{section}

\begin{document}

\title[Complexity and degeneracy in QRL]
{Complexity scaling and optimal policy degeneracy in quantum reinforcement
learning via analytically solvable unitary-control-then-measure models}

\author[A.~Cintio]{Andrea Cintio}
\address[A.~Cintio]{Institute for Chemical and Physical processes, National Research Council (CNR) \\ via G.~Moruzzi 1 \\ 56124 Pisa.}
\email{andrea.cintio@pi.ipcf.cnr.it}

\author[A.~Michelangeli]{Alessandro Michelangeli}
\address[A.~Michelangeli]{Mathematics and Science Department, American University in Bulgaria, AUBG, 1 Georgi Izmirliev Sq., 2700 Blagoevgrad, \\
and Alexander von Humboldt Foundation, Bonn\\ and  Trieste Institute for Theoretical Quantum Technologies, TQT Trieste}
\email{amichelangeli@aubg.edu}

\author[D.~Tsutskov]{Dmitrii V.~Tsutskov}
\address[D.~Tsutskov]{American University in Bulgaria, AUBG, 1 Georgi Izmirliev Sq., 2700 Blagoevgrad}
\email{DNT230@aubg.edu}

\date{\today}



\thanks{\emph{Acknowledgements.} For this work, A.M.~gratefully acknowledges support from the Italian National Institute for Higher Mathematics INdAM, the David Flanagan funds at the AUBG, American University Bulgaria, and the Alexander von Humboldt Foundation, Germany.}

\begin{abstract}
We propose and analyse a class of analytically solvable models of quantum
reinforcement learning (QRL), formulated as finite-horizon Markov decision processes
in finite-dimensional Hilbert spaces. The models are built around a
`unitary-control-then-measure' protocol, in which a learning agent applies unitary
transformations to a quantum state and interleaves each control step with a
projective measurement onto a prescribed reference basis. Exact closed-form
expressions for trajectory probabilities, rewards, and the expected return are
derived for four concrete realisations: a closed-chain and an anti-periodic qubit
implementation, a qutrit model with ladder coupling, and a four-level two-qubit
system. Two structural features of these QRL protocols are then analysed.
First, we identify and quantify the reduction in the computational complexity
of the expected return, from the nominally exponential $O(e^N)$ scaling in the
trajectory length~$N$ to an explicit power-law $O(N^{\mathcal{I}})$, driven by
two rigorously established mechanisms, a trajectory equivalence and a sparsity
of the transition graph, besides a third, conjectured one: a spectral
concentration of the return, at the optimal policy, onto the polynomially
populated trajectory classes.
Second, we characterise the degeneracy of optimal policies. The low-dimensional
models exhibit unique optima whose asymptotic behaviour with~$N$ is governed by the
quantum Zeno effect, while the four-level system displays both plateau-type
quasi-degeneracy at large horizons and genuine discrete degeneracy at critical energy
parameters -- phenomena with no counterpart in the measurement-free quantum
optimal control landscape.
\end{abstract}

\maketitle

 \section{Introduction}\label{sec:intro}

This work addresses two main goals, of relevance to both physicists and computer
scientists working at the intersection of quantum mechanics and machine learning.

On one hand, we propose and analyse some representative realisations of a class of transparent, physically rigorous, analytically solvable `sandbox' models of quantum reinforcement learning. They are `unitary control + measure' protocols that retain the core features of a learning agent's exploration, exploitation, and reward maximisation.

Within our simplified and controlled framework we derive exact, closed-form expressions
for the expected return, trajectory probabilities, and optimal policies, replacing
black-box numerical evaluation with transparent analytical structure grounded in basic
linear algebra and finite-dimensional quantum mechanics.

In addition to that, we discuss in a rigorous mathematical framework the emergence of two relevant features of quantum reinforcement learning, the order of complexity and the degeneracy of optimisations, which are both critical for the identification of optimal strategies.

 We can infer sufficiently general conclusions and formulate problems that are applicable to more robust and sophisticated settings.

\emph{Reinforcement learning} (RL) \cite{Kaelbling1996reinforcement,SuttonBarto1998,Dong-Chen-Li-Tarn-2008,arulkumaran2017deep,Bertsekas2019reinforcement,Szepesvari2022algorithms,Bertsekas2024course} is a framework for an agent to learn how to make sequential decisions by interacting with an environment, and mathematically it consists of certain types of optimisation problems to find an optimal policy that maximises the expected cumulative reward.

From early foundations the field really coalesced in the 1980's through the connections of ideas from dynamic programming, control theory, and machine learning, and has since seen tremendous growth, especially in the 2010's with the rise of deep neural networks, machine learning, and artificial intelligence.

The classical formalisation of sequential decision making in which the actions affect immediate rewards as well as subsequent situations, and hence transition probabilities and rewards are initially unknown and must be learned through interaction with the environment, are the \emph{Markov decision processes} \cite{puterman2014markov}, \cite[Chapter 3]{SuttonBarto1998}, \cite[Chapter 1]{Szepesvari2022algorithms}. Within this mathematical idealisation, precise theoretical statements can be made: in their essence, Markov decision processes consist of optimisation over policy-dependent probability measures on a suitable space of sequential decisions.

Reinforcement learning has emerged as a powerful paradigm in modern physics, particularly for the control and optimisation of quantum systems.

 In fact, in a variety of applications (including quantum chemistry simulations, quantum error correction, and even financial portfolio optimisation and drug discovery) protocols of RL effectively consist of an interaction of the agent with a quantum system. This is customarily referred to as \emph{quantum reinforcement learning} \cite{dong2008quantum,meyer2022survey,Schenk2024hybrid,Wu2025quantum}, QRL, especially from the perspective of combining quantum computing principles with RL to potentially achieve computational advantages.

 Naturally QRL encompasses a variety of approaches such as quantum-enhanced classical RL (the use of quantum algorithms to speed up components like policy search and linear algebra operations), quantum environments where agents learn to control quantum systems and the state space consists of quantum states, and quantum agents that leverage superposition and entanglement for decision-making and parallel exploration of action spaces.

The present discussion sits at the overlap of two different disciplines -- and cultures, in a sense: machine learning and quantum physics. Already providing a thorough discussion of such an intersection is a fascinating subject in itself, particularly given the unprecedented surge of activity in this area and the current lack of a systematic theoretical framework.

 In our analysis we have three main points of focus. The first is the very mathematical modelling and the analytic solvability. The analytical approach reveals structural features that brute-force numerical methods would probably miss. In particular, working out explicit combinatorial summations and implementing flow conservation principles we replace `black box' numerics with transparent maths.

 The second is on how the numerical computation of the expected return scales with $N$, in order to pinpoint which specific parts of the QRL protocols drive the customary reduction from a nominal exponential complexity in $N$ to a much faster power-law dependence. Two such mechanisms are rigorously identified, a \emph{trajectory equivalence} and a \emph{sparsity of the transition graph}, both operating for every admissible policy; the evidence from the models further points to a third mechanism, of conjectural character, a \emph{spectral concentration at the optimum} whereby only a polynomially small minority of the trajectory equivalence classes effectively supports the optimisation.

 The third is on the issue of possible discrete or also continuous multiplicity of
optimal strategies. This is a genuinely under-explored problem in the QRL literature.
 Numerically, degeneracy may be monitored by means of various advanced, specialised techniques, carefully sampling and analysing the expected return function. However, this can get tricky if distinct optimal policies are in some sense far apart or also are arranged in a continuous range. In this regard, we want to explore simple, yet insightful, instances where (non-)degeneracy arises, tracing it back to the various internal components of the model.

 After presenting the general model in Section \ref{sec:generalQRLmodel}, we will be able to rephrase the above questions more precisely. In Sections \ref{sec:closequbit} through \ref{sec:4levels} we then specialise the general model for several concrete `toy-model' realisations: each of these is analysed analytically and then optimised numerically. The evidence from these sample cases is finally re-examined in the concluding Sections: in Section \ref{sec:complexitycomparison} we isolate the distinct mechanisms responsible for the drop of the computational complexity of the expected return, and in Section \ref{sec:deg} we organise the evidence on the degeneracy of the optimal policies -- in both cases elaborating specific answers and formulating general problems deferred to subsequent investigation.

\section{The Quantum Reinforcement Learning model}\label{sec:generalQRLmodel}

 Let us introduce in this Section a class of simplified QRL models, formulated as finite-horizon Markov decision processes in a finite-dimensional quantum framework.

 Specific realisations of such models will be then analysed in the sequel.

 In the language of quantum reinforcement learning, the models consist of an agent executing the following protocol and finding the strategy that maximises the final reward gained.

 The agent must explore multiple times sequences of $N$ consecutive controlled manipulations of a quantum state according to a general strategy, each state manipulation being followed by the observation (collapse) of the resulting state onto a prescribed reference state that will be the starting point for the next manipulation.

 The reference states are selected out of a given orthonormal basis
\begin{equation}\label{eq:statesS}
 \mathcal{S}\,=\,\{\Psi^{(1)},\dots,\Psi^{(d)}\}
\end{equation}
of vectors of an underlying $d$-dimensional Hilbert space $\cH$ over $\mathbb{C}$. Each state manipulation is a unitary transformation (a rotation) that creates a superposition of basis vectors. The observation is a collapse that produces a pure state ready for the next rotation.

The agent explores all possible paths of `manipulations + observations' from the $0$-th to the $N$-th state; for all paths, initial and final states are prescribed or left free. Thus, a path is of the type
\begin{equation}\label{eq:path}
 \Psi_{\textrm{initial}}\,\equiv\,|\psi_0\rangle\to|\psi_1\rangle\to|\psi_2\rangle\to\cdots\to|\psi_{N-1}\rangle\to|\psi_N\rangle\,\equiv\,\Psi_{\textrm{final}}
\end{equation}
for some choice $\psi_1,\dots,\psi_{N-1}\in\mathcal{S}$.

For each trial, the player chooses a policy $\pi$ (a `deterministic' policy) that consists of declaring $d$ unitary transformations
\begin{equation}\label{eq:choosepol}
 \pi(\Psi^{(1)}),\dots , \pi(\Psi^{(d)})
\end{equation}
in the Hilbert space $\cH$, each one labelled by a state from $\mathcal{S}$. In experimental quantum systems, these unitaries would correspond to applied control fields (like laser pulses, microwave fields, or magnetic fields) that induce unitary evolution. They can be seen as preparation operations before each measurement, analogous to choosing a measurement basis in quantum information protocols.

At each trial, along a path like \eqref{eq:path} the agent performs the unitary transformation $\pi(\psi_i)$ of a state $|\psi_i\rangle$ and the resulting state $\pi(\psi_i)|\psi_i\rangle$ is collapsed onto the next measurement state $|\psi_{i+1}\rangle$. This identifies for each path a trajectory $\tau$ of alternating states and actions of the type
\begin{equation}\label{eq:tautrak}
\begin{split}
  \tau:\; \Psi_{\textrm{initial}}\equiv\,|\psi_0\rangle&\xrightarrow{\;\pi(\psi_0)\;}\pi(\psi_0)|\psi_0\rangle\xrightarrow{\,\textrm{(collapse)}\,}|\psi_1\rangle \\
  &\xrightarrow{\;\pi(\psi_1)\;}\pi(\psi_1)|\psi_1\rangle\xrightarrow{\,\textrm{(collapse)}\,}|\psi_2\rangle \\
  &\cdots \xrightarrow{\;\pi(\psi_{N-1})\;}\pi(\psi_{N-1})|\psi_{N-1}\rangle\xrightarrow{\,\textrm{(collapse)}\,}|\psi_N\rangle\,\equiv\,\Psi_{\textrm{final}}\,.
\end{split}
\end{equation}

The likelihood of each jump $|\psi_i\rangle\xrightarrow{\;\pi(\psi_i)\;}\pi(\psi_i)|\psi_i\rangle\xrightarrow{\,\textrm{(collapse)}\,}|\psi_{i+1}\rangle$ is given by the quantum transition probability $|\langle\psi_{i+1}|\pi(\psi_i)\psi_i\rangle|^2$. The product of all such probabilities gives the likelihood of the whole trajectory $\tau$:
\begin{equation}\label{eq:probtraj}
 P_\pi(\tau)\,=\,\prod_{i=0}^{N-1} \big|\langle \psi_{i+1}|\pi(\psi_i)\psi_i\rangle\big|^2\,.
\end{equation}
In fact, \eqref{eq:probtraj} naturally defines an actual probability on the set $\Gamma_{\pi,N}$ of length-$N$ trajectories with fixed initial and final state.

 Along each trajectory $\tau$ an energy assessment is made: the energy difference between the post-collapse state and the pre-collapse state determines the rewards. More precisely: an energy observable (quantum Hamiltonian) $H$ is assigned on $\cH$, which for concreteness is diagonal  on the basis $\mathcal{S}$, i.e., the $\Psi^{(j)}$'s are eigenstates of $H$ and in such basis $H$ is represented as
 \begin{equation}
  H\,=\,\sum_{k=1}^d E_k|\Psi^{(k)}\rangle\langle\Psi^{(k)}|\,=\,
  \begin{pmatrix}
    E_d    &        &        &    \\
           & \ddots &        &     \\
           &        & E_2    &     \\
           &        &        & E_1
  \end{pmatrix}
 \end{equation}
 with energy levels conventionally ordered as $E_1< E_2 < \cdots< E_d$. The reward when transitioning from state $\psi_i$ to state $\psi_{i+1}$ is
 \[
  \langle \psi_{i+1} , H \psi_{i+1}\rangle - \langle \pi(\psi_i)\psi_i, H \pi(\psi_i)\psi_i\rangle\,,
 \]
 namely the difference in the energy expectations before and after the collapse. All such individual rewards add up to the total reward for the considered trajectory $\tau$,
 \begin{equation}\label{eq:rewardRnodisc}
  \begin{split}
   R_\pi(\tau)\,&=\,\sum_{i=0}^{N-1}\Big(\langle \psi_{i+1} , H \psi_{i+1}\rangle - \langle \pi(\psi_i)\psi_i, H \pi(\psi_i)\psi_i\rangle\Big)\,.
  \end{split}
 \end{equation}

 By exploring all possible trajectories with the same chosen policy $\pi$ the player computes the expected return $J(\pi)$ for that policy, defined as the expectation of the reward $R_\pi(\tau)$ with respect to the probability distribution of the admissible trajectories:
 \begin{equation}\label{eq:Jpidi}
  J(\pi)\,=\sum_{\tau\in\Gamma_{\pi,N}} R_\pi(\tau) P_\pi(\tau)
 \end{equation}

 The goal finally is: \emph{to identify the policy $\pi^*$ that maximises the expected return}. Thus,
 \begin{equation}\label{eq:optpolpistar}
  \pi^*\,=\,\mathrm{argmax}_\pi\, J(\pi)\,.
 \end{equation}

 This is a typical \emph{finite-horizon RL} problem with stationary deterministic policy, where
 \begin{itemize}
  \item the collection $\mathcal{S}$ from \eqref{eq:statesS} is the set of \emph{states},
  \item $U(d)$ (the group of unitary operators on $\cH$) is the set of \emph{actions},
  \item the maps $\pi:\mathcal{S}\to U(d)$ as in \eqref{eq:choosepol} are the \emph{policies},
  \item $|\langle\psi_{i+1}|\pi(\psi_i)\psi_i\rangle|^2$ is the \emph{state transition probability function},
  \item $\langle \psi_{i+1} , H \psi_{i+1}\rangle - \langle \pi(\psi_i)\psi_i, H \pi(\psi_i)\psi_i\rangle$ is the \emph{reward function}.
 \end{itemize}
 As trajectories have all finite length $N$ (finite horizon), the multiplicative \emph{discount factor} $\gamma^i$, $\gamma\in(0,1]$, that customarily appears in the $i$-th summand of \eqref{eq:rewardRnodisc} to ensure convergence for infinite trajectories, is omitted here ($\gamma\equiv 1$).

\section{Connections to quantum control, quantum walks, and reinforcement learning}
\label{sec:context}

From a quantum-mechanical perspective, the protocol described constitutes a
specific instance of a \emph{quantum random walk}, or \emph{measured quantum
walk} \cite{AharonovDavZa1993,Venegas2012,Attal2012}. Structurally, this is
a `unitary-then-measurement' experiment characterised by an alternating cycle
of
\begin{itemize}
    \item initialisation, where the system begins in a state
      $|\Psi_{\textrm{initial}}\rangle$,
    \item control phase, where the agent applies unitary operators,
      representing the `control knobs' available to the experimenter to
      direct the evolution,
    \item measurement phase, where the agent performs projective
      measurements in the computational basis, introducing the probabilistic
      collapse inherent to quantum mechanics.
\end{itemize}

The unitary operators in the control phase are the central objects of the
\emph{quantum optimal control landscape} programme (see, e.g., \cite{Hsieh-Rabitz-2007,Hsieh-Wu-Rabitz-Lidar-2010,Sugny-Kontz-2008,
Volkov-Myachkova-Pechen-2025}). The class of models considered here operate in the
complementary regime to measurement-free unitary control: the systematic
interleaving of projective measurements converts what would otherwise be a
fidelity optimisation over a group of unitaries into a stochastic
expected-return maximisation over distributions on trajectory space. As we
shall analyse in Section~\ref{sec:deg}, it is precisely
this interplay that can generate genuine degeneracies of the optimal policy
with no counterpart in the measurement-free control landscape.

In the framework of artificial intelligence, this architecture aligns with
\emph{projective simulation} \cite{Briegel-DelasCuevas-2012,Paparo2014}: the
reference state collection $\mathcal{S}$ corresponds to the clip network, the
coupling structures analysed in Sections~\ref{sec:closequbit}-\ref{sec:4levels}
correspond to specific decision graphs, and the `glow' or `damping' parameters
correspond to the forgetting mechanisms used to handle non-stationary
environments.

More broadly, the `unitary control + measure' protocol is one among several
structurally distinct routes to applying reinforcement learning to
quantum-mechanical problems. A representative approach at the opposite end of
the modelling spectrum exploits the Feynman-Kac representation of the
imaginary-time Schr\"odinger equation \cite{Barr-Gispen-Lamacraft-2020},
yielding a variational route to ground-state energies via neural-network
function approximators. The two approaches differ at every structural level --
continuous versus discrete state space, imaginary versus real time, stochastic
differential equations versus projective collapse, neural representations
versus exact analytical formulae -- and are accordingly complementary rather
than competing. The contrast serves to highlight what is structurally
distinctive about the present protocol: the projective measurement step is not
an incidental feature but the element that simultaneously defines the
probability law on trajectories \eqref{eq:probtraj}, drives the reduction in
computational complexity discussed in
Section~\ref{sec:complexitycomparison}, and, in conjunction with the policy
structure, is responsible for the landscape degeneracies identified therein.

\section{Closed-chain qubit implementation}\label{sec:closequbit}

 From this Section until Section \ref{sec:4levels} the preceding general model shall be specialised to certain simplified settings that retain a level of instructiveness and can be worked out in detail. The simplest non-trivial case is the setting where the states are qubits.

 \subsection{The model}\label{sec:2levelclosed}~

 The agent explores chains of qubit states across $\mathbb{C}^2$:
 \begin{equation}
  d=\mathrm{dim}\cH=2\,,\qquad \cH=\mathbb{C}^2\,.
 \end{equation}

 \medskip

 \textbf{States.} For concreteness, let the set of states be
 \begin{equation}\label{eq:states2}
  \mathcal{S}\,=\,\{\Psi^{(1)},\Psi^{(2)}\}
 \end{equation}
 with
 \begin{equation}
  \Psi^{(1)}\equiv\Psi^-\equiv |-\rangle =:\begin{pmatrix} 0 \\1 \end{pmatrix}\,,\qquad \Psi^{(2)}\equiv\Psi^+\equiv |+\rangle =:\begin{pmatrix} 1 \\ 0 \end{pmatrix}\,,
 \end{equation}
 and the Hamiltonian be
 \begin{equation}\label{eq:Hamilt22}
  H\,=\,E_+|+\rangle\langle +|+ E_-|-\rangle\langle -|\,=\,
  \begin{pmatrix}
   E_+ & 0 \\
   0 & E_-
  \end{pmatrix}\,.
 \end{equation}
 In fact, \eqref{eq:rewardRnodisc}, hence also \eqref{eq:Jpidi}, is \emph{linear} in $H$, and there are only two energy levels $E_+> E_-$. Thus, upon a non-restrictive overall shift in the energy and energy-rescaling, there is no loss of generality in setting $E_+=1$, $E_-=0$.

 \medskip

 \textbf{Actions.} The set of actions is $U(2)$, the unitary $2\times 2$ matrices. Recall that an element $U\in U(2)$ is characterised by four real parameters $\alpha,\beta,\gamma,\theta\in[0,2\pi)$ as
 \begin{equation}\label{eq:unitaries2}
  U \,=\, e^{i\alpha} \begin{pmatrix}
e^{i\beta}\cos\theta & e^{i\gamma}\sin\theta \\
-e^{-i\gamma}\sin\theta & e^{-i\beta}\cos\theta
\end{pmatrix}\,.
 \end{equation}
 Obviously, the overall phase $e^{i\alpha}$ is always cancelled both in the evaluation of the transition probability $|\langle \psi',U\psi\rangle|^2$ and of the energy expectation $\langle U\psi, H U\psi\rangle$, so one can directly set $\alpha\equiv 0$.

\medskip

 \textbf{Policies.} Choosing a policy $\pi$ amounts to declaring two unitaries
 \begin{equation}\label{eq:policies2}
   \begin{split}
     \pi_+\,\equiv\,\pi(\Psi^+)\,, & \qquad\textrm{and therefore parameters $\beta^+,\gamma^+,\theta^+$}\,, \\
     \pi_-\,\equiv\,\pi(\Psi^-)\,, & \qquad\textrm{and therefore parameters $\beta^-,\gamma^-,\theta^-$}\,.
   \end{split}
 \end{equation}

\medskip

 \textbf{End-points.} For the remaining part of this Section, let us consider the closed-chain scenario
 \begin{equation}\label{eq:initialstate2pp}
  |\Psi_{\textrm{initial}}\rangle\,=\,|\Psi_{\textrm{final}}\rangle\,=\,|+\rangle\,.
 \end{equation}

 \medskip

 \textbf{Trajectories.}
 The set $\Gamma_{N,\pi}$ consists of $2^{N-1}$ trajectories, as many as the distinct ways to allocate (with ordering) an amount $N-1$ of $+$'s and $-$'s into
\begin{center}
 \begin{tabular}{|*{11}{c|}}
\hline
$+$ & $\cdots$ & $\cdots$ & $\cdots$ & $\cdots$ & $\cdots$ & $\cdots$ & $\cdots$ & $\cdots$ & $\cdots$ & $+$ \\
\hline
\end{tabular}
\end{center}
In diagrams like the one above, the convention shall be that the trajectory's steps go \emph{from left to right}. Along a trajectory only four types of jumps occur, denoted as
\[
 [++]\,,\qquad [-+]\,,\qquad[+-]\,,\qquad[--]\,
\]
respectively. Thus: $[-+]$ is the (left-to-right) jump $|-\rangle\to |+\rangle$, etc.
Now, for a length-$N$ trajectory with fixed $+$ at the edges, with $\textsf{p}$ states $|+\rangle$ ($\textsf{p}\in\{2,\dots,N+1\}$) and hence $N+1-\textsf{p}$ states $|-\rangle$, the numbers
\begin{equation}\label{eq:trajparamabcd}
  \begin{split}
   \textsf{a}\,&:=\,\#[++]\,, \\
   \textsf{b}\,&:=\,\#[-+]\,, \\
   \textsf{c}\,&:=\,\#[+-]\,, \\
   \textsf{d}\,&:=\,\#[--]\, \\
 \end{split}
\end{equation}
are constrained by the conditions
\begin{equation}\label{eq:2condparam1}
 \begin{split}
   \textsf{a}+\textsf{b}+\textsf{c}+\textsf{d}\,&=\,N\,, \\
   \textsf{a}+\textsf{b}\,&=\,\textsf{p}-1\,, \\
   \textsf{c}+\textsf{d}\,&=\,N+1-\textsf{p}\,, \\
   \textsf{b}\,&=\,\textsf{c}\,.
 \end{split}
\end{equation}
Example:
\begin{center}
 \begin{tabular}{|*{11}{c|}}
\hline
$+$ & $+$ & $-$ & $-$ & $+$ & $-$ & $+$ & $-$ & $-$ & $+$ & $+$ \\
\hline
\end{tabular}
\end{center}
\begin{equation}
 N=10\,,\quad \textsf{p}=6\,,\quad \textsf{a}=2\,,\quad \textsf{b}=3 \,,\quad \textsf{c}=3\,,\quad \textsf{d}=2\,.
\end{equation}
Therefore, length-$N$ trajectories with $\textsf{p}$ states of type $|+\rangle$ have
\begin{equation}\label{eq:abdtocp}
  \begin{split}
  \textsf{a}\,&=\,\textsf{p}-1-\textsf{c}\,, \\
  \textsf{b}\,&=\,\textsf{c}\,, \\
  \textsf{c}\,&\in\big\{0\,,\,\min\{\textsf{p}-1,N+1-\textsf{p}\}\big\}\,, \\
  \textsf{d}\,&=\,N+1-\textsf{p}-\textsf{c}\,.
 \end{split}
\end{equation}
The number of valid trajectories with fixed endpoints, containing $\textsf{p}$ positive states and $2\mathsf{c}$ sign changes, is determined by the number of ways to arrange the states into alternating blocks. Since the chain starts and ends with $|+\rangle$, there must be $\mathsf{c}+1$ blocks of $|+\rangle$ states and $\mathsf{c}$ blocks of $|-\rangle$ states interspersed between them. Using the stars and bars method for compositions, one distributes $\textsf{p}$ items into $\mathsf{c}+1$ non-empty bins and $N+1-\textsf{p}$ items into $\mathsf{c}$ non-empty bins. This yields the exact multiplicity:
\begin{equation}\label{eq:Ncount2pc}
 \mathcal{N}(N,\textsf{p},\mathsf{c}) \,=\, \binom{\textsf{p}-1}{\mathsf{c}} \cdot \binom{N-\textsf{p}}{\mathsf{c}-1}\,.
\end{equation}

\medskip

 \textbf{Probability of a trajectory.} Next, here is how \eqref{eq:probtraj} and \eqref{eq:rewardRnodisc} are specialised for this scenario. Concerning the probability of a trajectory $\tau$ with $\textsf{p}$ positive states and $2\mathsf{c}$ sign changes, this is now the quantity
\begin{equation}
 \begin{split}
   P_\pi(\textsf{p},\textsf{c})\,&=\, \prod_{i=0}^{N-1} \big|\langle \psi_{i+1}|\pi(\psi_i)\psi_i\rangle\big|^2\\
  & =\, \big|\langle +|\pi_+|+\rangle\big|^{2\textsf{a}}\cdot \big|\langle +|\pi_-|-\rangle\big|^{2\textsf{b}}\cdot \big|\langle -|\pi_+|+\rangle\big|^{2\textsf{c}}\cdot \big|\langle -|\pi_-|-\rangle\big|^{2\textsf{d}}
 \end{split}
\end{equation}
with $\textsf{a},\textsf{b},\textsf{d}$ given by \eqref{eq:abdtocp}. From \eqref{eq:unitaries2}-\eqref{eq:policies2},
\begin{equation}\label{eq:2probabilitiespm}
 \begin{split}
  |\langle +|\pi_+|+\rangle|^2\,&=\,|\cos\theta_+|^2\,,\qquad |\langle +|\pi_-|-\rangle|^2\,=\,|\sin\theta_-|^2\,, \\
  |\langle -|\pi_+|+\rangle|^2\,&=\,|\sin\theta_+|^2\,,\qquad\, |\langle -|\pi_-|-\rangle|^2\,=\,|\cos\theta_-|^2\,,
 \end{split}
\end{equation}
 whence
 \begin{equation}\label{eq:ProbTraj2}
  P_\pi(\textsf{p},\textsf{c})\,=\,|\cos\theta_+|^{2(\textsf{p}-1-\textsf{c})}\cdot|\sin\theta_-|^{2\textsf{c}}\cdot|\sin\theta_+|^{2\textsf{c}}\cdot|\cos\theta_-|^{2(N+1-\textsf{p}-\textsf{c})}\,.
 \end{equation}

\medskip

 \textbf{Reward along a trajectory.} According to \eqref{eq:rewardRnodisc}, and recalling that in \eqref{eq:Hamilt22} one has $E_+=1$, $E_-=0$, the total reward along a trajectory $\tau$ with $\textsf{p}$ positive states is now the quantity
 \begin{equation}\label{eq:RewardTraj2}
  \begin{split}
   R_\pi(\textsf{p},\textsf{c})\,&=\,\sum_{i=0}^{N-1}\Big(\langle \psi_{i+1} , H \psi_{i+1}\rangle - \langle \psi_i, \pi(\psi_i)^*H \pi(\psi_i)\psi_i\rangle\Big) \\
   &=\,(\textsf{p}-1)-(\textsf{p}-1)\langle +|\pi_+^*H\pi_+|+\rangle-(N+1-\textsf{p})\langle -|\pi_-^*H\pi_-|-\rangle \\
   &=\,(\textsf{p}-1)-(\textsf{p}-1)|\cos\theta_+|^2-(N+1-\textsf{p})|\sin\theta_-|^2\,.
  \end{split}
 \end{equation}
 Recall that the above expression is dressed by an overall multiplicative factor $\Delta E$ in those applications where one intends to monitor the dependence of the reward on the energy difference $\Delta E$ between the two energy levels (here above set to $\Delta E\equiv 1$).

\medskip

 \textbf{Expected return.} In view of \eqref{eq:abdtocp}, \eqref{eq:Ncount2pc}, \eqref{eq:ProbTraj2}, and \eqref{eq:RewardTraj2}, now \eqref{eq:Jpidi} yields
  \begin{equation}\label{eq:Jpidi2}
   \begin{split}
    &J(\pi)\,=\sum_{\tau\in\Gamma_{\pi,N}} R_\pi(\tau) P_\pi(\tau) \\
    &=\,\sum_{\textsf{p}=2}^{N} \;\sum_{\mathsf{c}=1}^{\min\{\textsf{p}-1,N+1-\textsf{p}\}} \mathcal{N}(N,\textsf{p},\mathsf{c}) \cdot R_{\pi}(\textsf{p},\mathsf{c}) \cdot P_{\pi}(\textsf{p},\mathsf{c}) \;\; + \;\; R_\pi(N+1,0)P_\pi(N+1,0) \\
    &=\,\sum_{\textsf{p}=2}^{N} \;\sum_{\mathsf{c}=1}^{\min\{\textsf{p}-1,N+1-\textsf{p}\}}\binom{\textsf{p}-1}{\mathsf{c}} \cdot \binom{N-\textsf{p}}{\mathsf{c}-1}\times \\
    &\qquad\qquad\qquad\times\Big(|\cos\theta_+|^{2(\textsf{p}-1-\textsf{c})}\cdot|\sin\theta_-|^{2\textsf{c}}\cdot|\sin\theta_+|^{2\textsf{c}}\cdot|\cos\theta_-|^{2(N+1-\textsf{p}-\textsf{c})}\Big) \\
    &\qquad\qquad\qquad\times\big((\textsf{p}-1)(1-\cos^2\theta_+) -(N+1-\textsf{p})\sin^2\theta_-\big) \\
    &\qquad + N(1-\cos^2\theta_+)|\cos\theta_+|^{2N}\,.
   \end{split}
  \end{equation}
 Observe that the term separated out corresponds to the unique trajectory with $\textsf{p}=N+1$ and $\mathsf{c}=0$ (all $|+\rangle$ states), as the combinatorial factor for $\mathsf{c} \geqslant 1$ covers all trajectories with at least one excursion to $|-\rangle$.

 This shows that the expected return for the policy $\pi$ only depends on the \emph{two} parameters $\sin\theta_\pm$ of the policy, which can be non-restrictively taken in $[0,\pi]$. Re-labelling
 \begin{equation}
  x_\pm\,:=\,\sin^2\theta_\pm\,\in[0,1]
 \end{equation}
 one re-writes
   \begin{equation}\label{eq:Jpidi2bis}
    \begin{split}
     J(\pi)\,&=\,\sum_{\textsf{p}=2}^{N} \;\sum_{\mathsf{c}=1}^{\min\{\textsf{p}-1,N+1-\textsf{p}\}}\binom{\textsf{p}-1}{\mathsf{c}} \cdot \binom{N-\textsf{p}}{\mathsf{c}-1}\big((\textsf{p}-1)x_+ -(N+1-\textsf{p})x_-\big) \\
     &\qquad\qquad\qquad\times\Big(x_+^{\textsf{c}}(1-x_+)^{(\textsf{p}-1-\textsf{c})}\,x_-^{\textsf{c}}(1-x_-)^{(N+1-\textsf{p}-\textsf{c})} \Big) \\
     &\qquad + N x_+ (1-x_+)^N \\
     &=:\,J(x_+,x_-)\,,
    \end{split}
   \end{equation}
   as a function over the square $(x_+,x_-)\in[0,1]\times[0,1]$.

 A crucial observation for later discussion is that while the first line of \eqref{eq:Jpidi2} formally scales as $O(2^N)$, \eqref{eq:Jpidi2bis} actually scales as $O(N^2)$.

 The latter conclusion is valid in general, but it does not prevent that for certain policies the double sum over the two labelling parameters $(\textsf{p},\textsf{c})$ in \eqref{eq:Jpidi2bis} collapses to a single sum over just one of them, thereby dropping further the complexity to $O(N)$.

 An example of this occurrence is the set of policies $\widetilde{\pi}$ defined by
 \begin{equation}\label{eq:pi1-constraint}
 x_+ + x_- \,=\, 1\,,\qquad x\,:=\,x_+\,\in\,[0,1]\,,
\end{equation}
i.e., $\sin^2\theta_+ + \sin^2\theta_- = 1$ (hence $\theta_- = \pm\big(\frac{\pi}{2}-\theta_+\big)$, modulo the $2\pi$-periodicity of the rotation angles).
Under the constraint \eqref{eq:pi1-constraint}, the probability factor in \eqref{eq:Jpidi2bis} simplifies as
\begin{equation}\label{eq:Ppi1}
\begin{split}
 x_+^{\textsf{c}}(1-x_+)^{\textsf{p}-1-\textsf{c}}\,x_-^{\textsf{c}}(1-x_-)^{N+1-\textsf{p}-\textsf{c}}\,&=\,x^{\textsf{c}}(1-x)^{\textsf{p}-1-\textsf{c}}\,(1-x)^{\textsf{c}}\,x^{N+1-\textsf{p}-\textsf{c}} \\
 &=\,x^{N+1-\textsf{p}}\,(1-x)^{\textsf{p}-1}\,,
\end{split}
\end{equation}
and the reward factor in \eqref{eq:Jpidi2bis} simplifies as
\begin{equation}\label{eq:Rpi1}
 (\textsf{p}-1)\,x_+ - (N+1-\textsf{p})\,x_- \,=\, (\textsf{p}-1)\,x - (N+1-\textsf{p})(1-x) \,=\,(\textsf{p}-1) - N(1-x)\,.
\end{equation}
 Both \eqref{eq:Ppi1} and \eqref{eq:Rpi1} depend on $\textsf{p}$ only, the transition-count parameter $\textsf{c}$ having dropped out entirely.
Accordingly, \eqref{eq:Jpidi2bis} takes the form
\begin{equation}\label{eq:Jpi1}
 J(\widetilde{\pi})\,=\,\sum_{\textsf{p}=2}^{N} x^{N+1-\textsf{p}}\,(1-x)^{\textsf{p}-1}\,\Big[(\textsf{p}-1) - N(1-x)\Big]\,\mathcal{S}_{\textsf{p}} \;+\; N\,x\,(1-x)^N
\end{equation}
with
\begin{equation}
 \mathcal{S}_{\textsf{p}}\,:=\sum_{\textsf{c}=1}^{\min\{\textsf{p}-1,\,N+1-\textsf{p}\}}\!\binom{\textsf{p}-1}{\textsf{c}}\binom{N-\textsf{p}}{\textsf{c}-1}\,=  \sum_{\textsf{c}=1}^{\min\{\textsf{p}-1,\,N+1-\textsf{p}\}}\mathcal{N}(N,\textsf{p},\mathsf{c})\,.
\end{equation}
 \eqref{eq:Jpi1} collapses to a single summation over $\textsf{p}$. The last summand therein, $N x(1-x)^N$, is the boundary contribution $Nx_+(1-x_+)^N$ of \eqref{eq:Jpidi2bis}.
 Actually, $\mathcal{S}_{\textsf{p}}$ is the cardinality of the trajectory sub-class with $\textsf{p}$ positive states, obtained by summing the multiplicity \eqref{eq:Ncount2pc} over all admissible $\textsf{c}$. Thus, the computational complexity of evaluating $J(\widetilde{\pi})$ drops from the $O(N^2)$ of the general policy case \eqref{eq:Jpidi2bis} to $O(N)$.

 \subsection{Policy optimisation and discussion}~

\begin{figure}[t!]
    \centering
   \includegraphics[width=6.8cm]{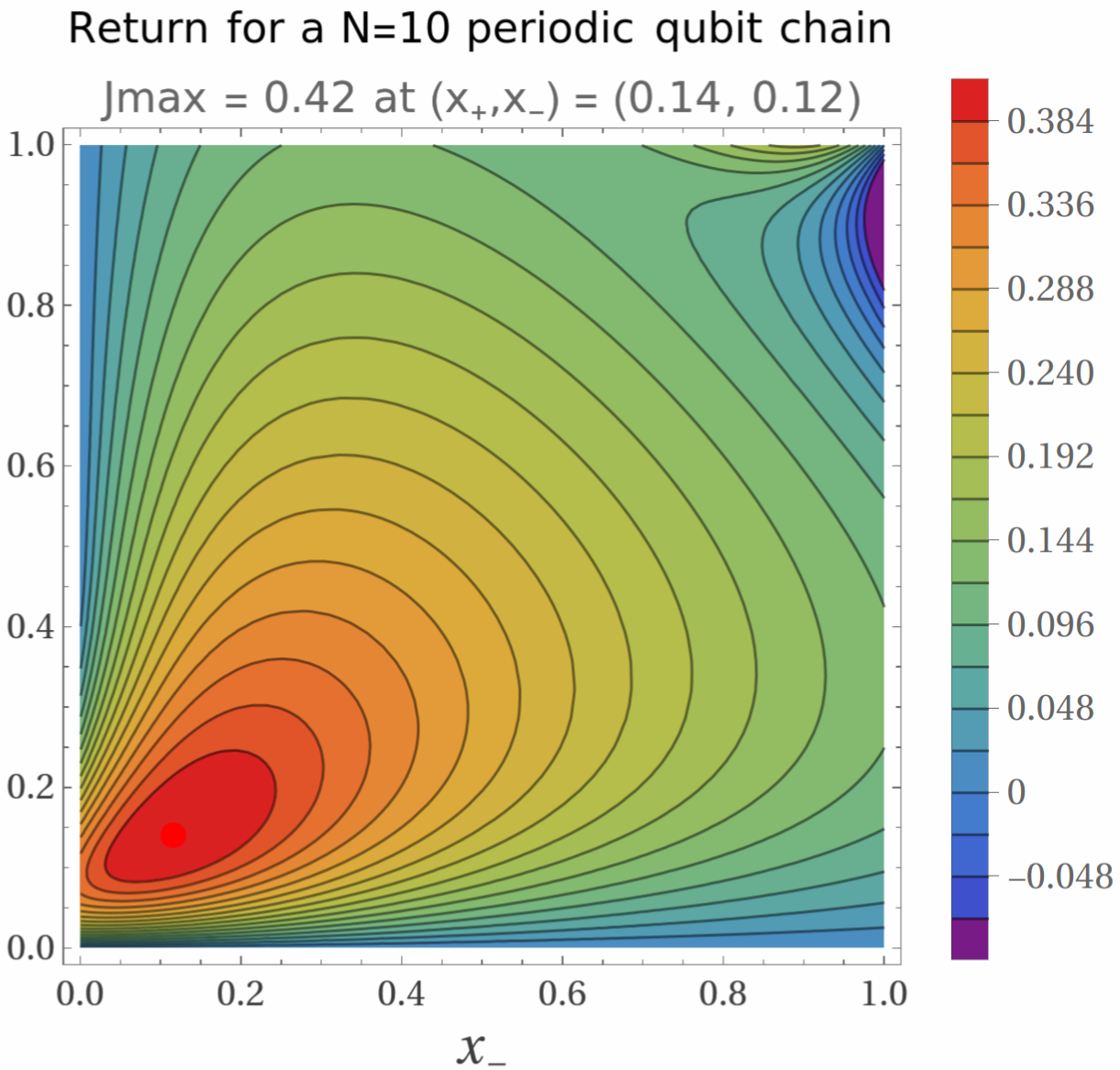}\quad
        \includegraphics[width=7.05cm]{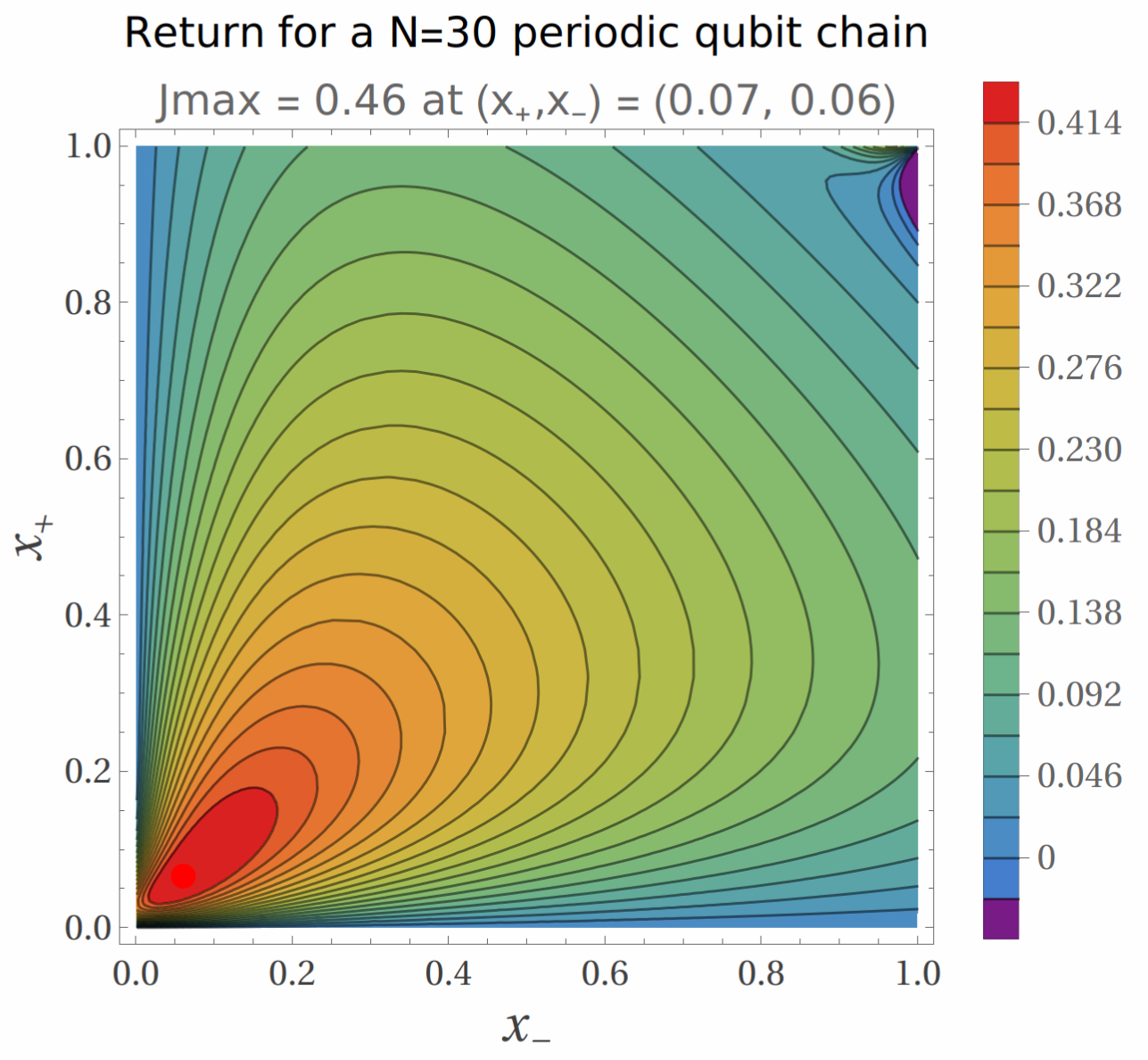}
    \caption{Numerical optimisation of the Expected Return \eqref{eq:Jpidi2bis}}
    \label{fig:Jpi-2levelPER}
\end{figure}

 The numerical optimisation of \eqref{eq:Jpidi2bis} over $(x_+,x_-)\in[0,1]\times[0,1]$ shows that the unique absolute maximum is reached at $(x_+,x_-)=(x_+^{\textrm{max}},x_-^{\textrm{max}})$ with $(x_+^{\textrm{max}},x_-^{\textrm{max}})$ approaching $0$ as $N$ grows (Figure \ref{fig:Jpi-2levelPER}). Thus, the optimising policies are those for which (modulo $2\pi$-periodicity) both $\theta_+$ and $\theta_-$ attain values $\pm\vartheta$ or $\pi\pm\vartheta$ for positive and suitable small $\vartheta$.

Mathematically, such a behaviour is explained by observing the structure of \eqref{eq:Jpidi2bis} in the regime of small parameters.
\begin{itemize}
    \item The scaling is driven by the boundary term: the isolated term $N x_+ (1-x_+)^N$ dominates the behaviour of the function near the origin. Elementary calculus shows that this term attains its maximum at $x_+ = \frac{1}{N+1} \approx \frac{1}{N}$. This imposes an overall scaling of $x_+ \sim O(1/N)$, forcing the peak towards zero as $N$ grows.

    \item The coupling requires non-zero $x_-$: the double sum, which represents the contribution from trajectories with multiple sign changes (governed by the power $\mathsf{c} \geqslant 1$), contains the factor $x_-^{\mathsf{c}}$. If $x_- = 0$, this entire `bulk' contribution vanishes. To gain the additional positive return from these terms (specifically from the $(\textsf{p}-1)x_+$ parts), the optimisation pushes $x_-$ to be strictly positive.

    \item The penalty limits $x_-$: the reward kernel inside the sum contains the negative term $-(N+1-\textsf{p})x_-$. Since this coefficient scales as $O(N)$, $x_-$ must scale as $O(1/N)$ to keep the penalty finite (order $O(1)$). Consequently, $x_-$ tracks the order of magnitude of $x_+$, balancing the activation of the sum's positive contribution against the large length-dependent penalty.
\end{itemize}

From a quantum mechanical perspective, the condition $x_\pm \ll 1$ implies that the optimal policies $\pi_\pm$ (see \eqref{eq:unitaries2}--\eqref{eq:policies2} above) involve rotation angles $\theta_\pm$ very close to zero. This corresponds to a `soft' driving regime that minimises the probability of sign changes (transitions between $|+\rangle$ and $|-\rangle$), thereby maintaining high state fidelity. By ensuring minimal disturbance during intermediate steps, the system preserves a well-defined energy state throughout the trajectory; this strategy leverages the constraint of fixed equal initial and final states ($|+\rangle$) to maximise the probability amplitude. The reward term, which scales with $N$, further incentivises this behaviour as the horizon length increases.

The optimisation results strongly evoke the \emph{quantum Zeno effect}, the phenomenon where frequent measurements of a quantum system inhibit its evolution away from the initial state \cite{FacchiPascazio_Zeno2008}. By adopting optimal policies, the system is effectively `frozen' close to its initial state through the synergy of minimal unitary rotation and frequent projective measurement. From the viewpoint of the quantum Zeno effect, the optimal policy exploits Zeno subspaces to maintain the system within high-reward energy states.

From a quantum control perspective, when managing the studied model over extended horizons with intermediate measurements, the optimal strategy approaches a regime of `minimal intervention'. This is somewhat counter-intuitive; in quantum control problems involving measurement, there is typically a fundamental trade-off between control aggressiveness and measurement-induced decoherence. In essence, this QRL model demonstrates that the optimal policy converges to a Zeno-like strategy, where the system is gently steered while measurements effectively `pin' it to the desired energy manifolds.

\section{Anti-periodic qubit implementation}\label{sec:qubitantiper}
 The counterpart version of the previous implementation with reversed endpoints consists of the same states, actions, and policies \eqref{eq:states2}-\eqref{eq:policies2}, whereas \eqref{eq:initialstate2pp} is replaced by
  \begin{equation}\label{eq:initialstate2pm}
  \Psi_{\textrm{initial}}\,=\,|+\rangle\,,\qquad \Psi_{\textrm{final}}\,=\,|-\rangle\,.
 \end{equation}
 The set $\Gamma_{N,\pi}$ still consists of $2^{N-1}$ trajectories whose parameters \eqref{eq:trajparamabcd} now obey
\begin{equation}\label{eq:2condparam2}
 \begin{split}
   \textsf{a}+\textsf{b}+\textsf{c}+\textsf{d}\,&=\,N\,, \\
   \textsf{a}+\textsf{c}\,&=\,\textsf{p}\,, \\
   \textsf{b}+\textsf{d}\,&=\,N-\textsf{p}\,, \\
   \textsf{b}\,&=\,\textsf{c}-1\,.
 \end{split}
\end{equation}
The key difference from the closed-chain scenario is the constraint $\textsf{b}=\textsf{c}-1$, as one must have one more $|+\rangle\to|-\rangle$ transition than $|-\rangle\to|+\rangle$ transitions to end with $|-\rangle$.
Additionally, we group transitions by their starting state (sources): all $\textsf{p}$ instances of $|+\rangle$ act as sources (hence $\textsf{a}+\textsf{c}=\textsf{p}$), whereas only $N-\textsf{p}$ instances of $|-\rangle$ act as sources (hence $\textsf{b}+\textsf{d}=N-\textsf{p}$).
Therefore,
\begin{equation}
  \begin{split}
  \textsf{a}\,&=\,\textsf{p}-\textsf{c}\,, \\
  \textsf{b}\,&=\,\textsf{c}-1\,, \\
  \textsf{c}\,&\in\big\{1\,,\,\min\{\textsf{p},N+1-\textsf{p}\}\big\}\,, \\
  \textsf{d}\,&=\,N-\textsf{p}-\textsf{c}+1\,.
 \end{split}
\end{equation}
 Along the same reasoning, the number of trajectories with $\textsf{p}$ positive states and $2\textsf{c}-1$ sign changes is now
\begin{equation}
 \mathcal{N}_{N,\textsf{p},\mathsf{c}} \,=\, \binom{\textsf{p}-1}{\mathsf{c}-1} \cdot \binom{N-\textsf{p}}{\mathsf{c}-1}\,,
\end{equation}
 which counts the ways to arrange the $\textsf{p}$ positive states into $\mathsf{c}$ distinct blocks and the $N+1-\textsf{p}$ negative states into $\mathsf{c}$ distinct blocks.
 This yields the following updated expression for the probability of a trajectory $\tau$ with $\textsf{p}$ positive states,
 \begin{equation}
  P_\pi(\textsf{p},\textsf{c})\,=\,|\cos\theta_+|^{2(\textsf{p}-\textsf{c})}\cdot|\sin\theta_-|^{2(\textsf{c}-1)}\cdot|\sin\theta_+|^{2\textsf{c}}\cdot|\cos\theta_-|^{2(N-\textsf{p}-\textsf{c}+1)}\,.
 \end{equation}
 The total reward along that trajectory must be updated to reflect the correct number of source states. The energy expectation of the final state is 0 (since $\Psi_N=|-\rangle$), and the sum of initial energies minus rotation costs involves $\textsf{p}$ operations of $\pi_+$ and $N-\textsf{p}$ operations of $\pi_-$. Thus:
 \begin{equation}
  \begin{split}
   R_\pi(\textsf{p},\textsf{c})\,&=\,\sum_{i=1}^{N} \langle \Psi_i, H \Psi_i\rangle - \sum_{i=0}^{N-1} \langle \psi_i, \pi(\psi_i)^*H \pi(\psi_i)\psi_i\rangle \\
   &=\,(\textsf{p}-1) - \Big( \textsf{p}\langle +|\pi_+^*H\pi_+|+\rangle + (N-\textsf{p})\langle -|\pi_-^*H\pi_-|-\rangle \Big) \\
   &=\,(\textsf{p}-1)-\textsf{p}|\cos\theta_+|^2-(N-\textsf{p})|\sin\theta_-|^2\,.
  \end{split}
 \end{equation}
 Therefore, the expected return for this model is
   \begin{equation}\label{eq:Jpidi2pm}
    \begin{split}
     J(\pi)\,&=\,\sum_{\textsf{p}=1}^{N} \sum_{\mathsf{c}=1}^{\min\{\textsf{p},N+1-\textsf{p}\}}\binom{\textsf{p}-1}{\mathsf{c}-1} \cdot \binom{N-\textsf{p}}{\mathsf{c}-1}\big((\textsf{p}-1) - \textsf{p}(1-x_+) -(N-\textsf{p})x_-\big) \\
     &\qquad\qquad\qquad\times\Big(x_+^{\textsf{c}}(1-x_+)^{(\textsf{p}-\textsf{c})}\,x_-^{(\textsf{c}-1)}(1-x_-)^{(N-\textsf{p}-\textsf{c}+1)} \Big) \\
     &=\,\sum_{\textsf{p}=1}^{N} \sum_{\mathsf{c}=1}^{\min\{\textsf{p},N+1-\textsf{p}\}}\binom{\textsf{p}-1}{\mathsf{c}-1} \cdot \binom{N-\textsf{p}}{\mathsf{c}-1}\big(\textsf{p} x_+ - 1 - (N-\textsf{p})x_- \big) \\
     &\qquad\qquad\qquad\times\Big(x_+^{\textsf{c}}(1-x_+)^{(\textsf{p}-\textsf{c})}\,x_-^{(\textsf{c}-1)}(1-x_-)^{(N-\textsf{p}-\textsf{c}+1)} \Big) \\
     & =:\,J(x_+,x_-)\,,
    \end{split}
   \end{equation}
 where
 \begin{equation}
 x_\pm\,:=\,\sin^2\theta_\pm\,\in[0,1]\,.
 \end{equation}

\begin{figure}[t!]
    \centering
   \includegraphics[width=6.8cm]{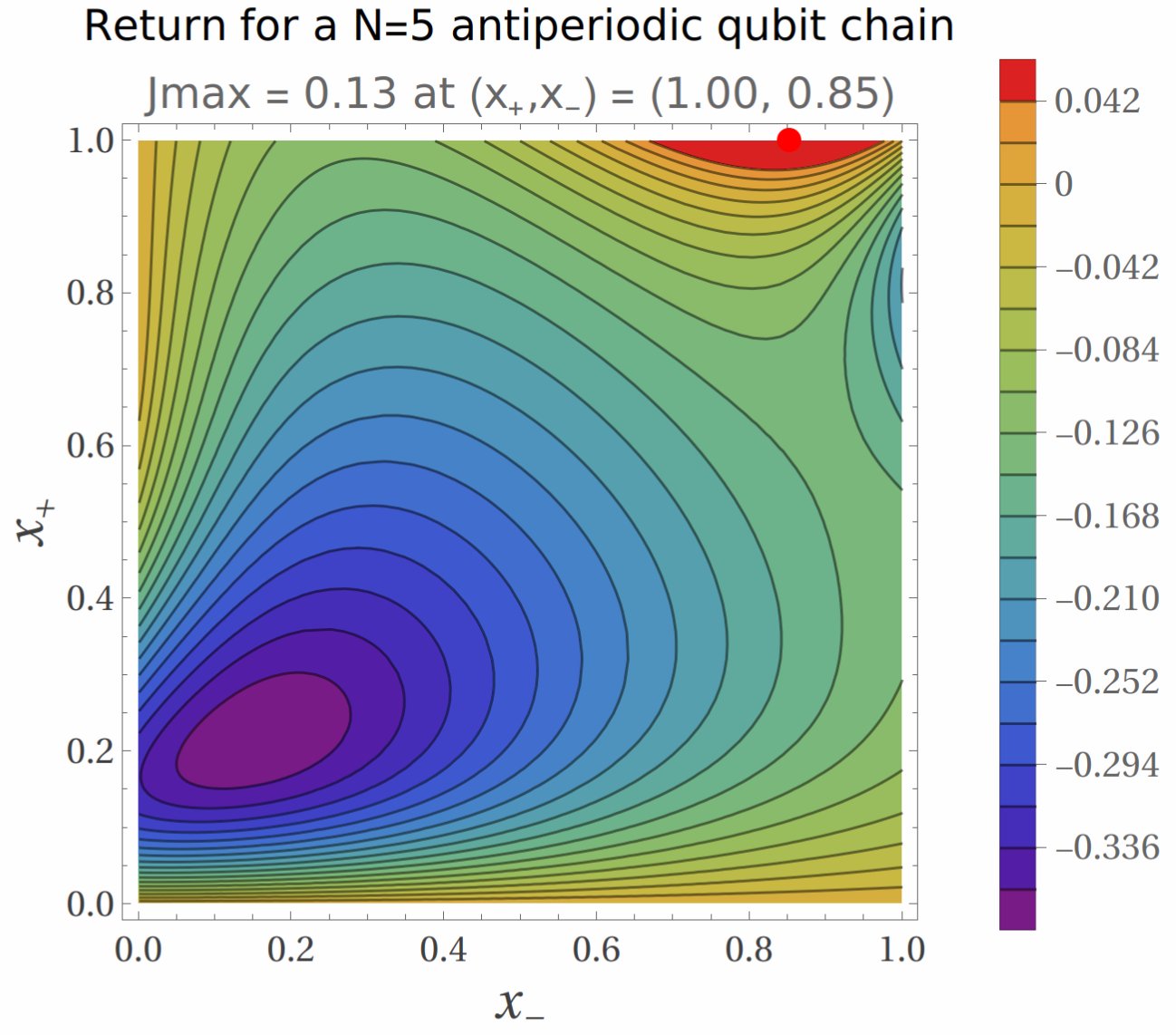}\quad
        \includegraphics[width=7.05cm]{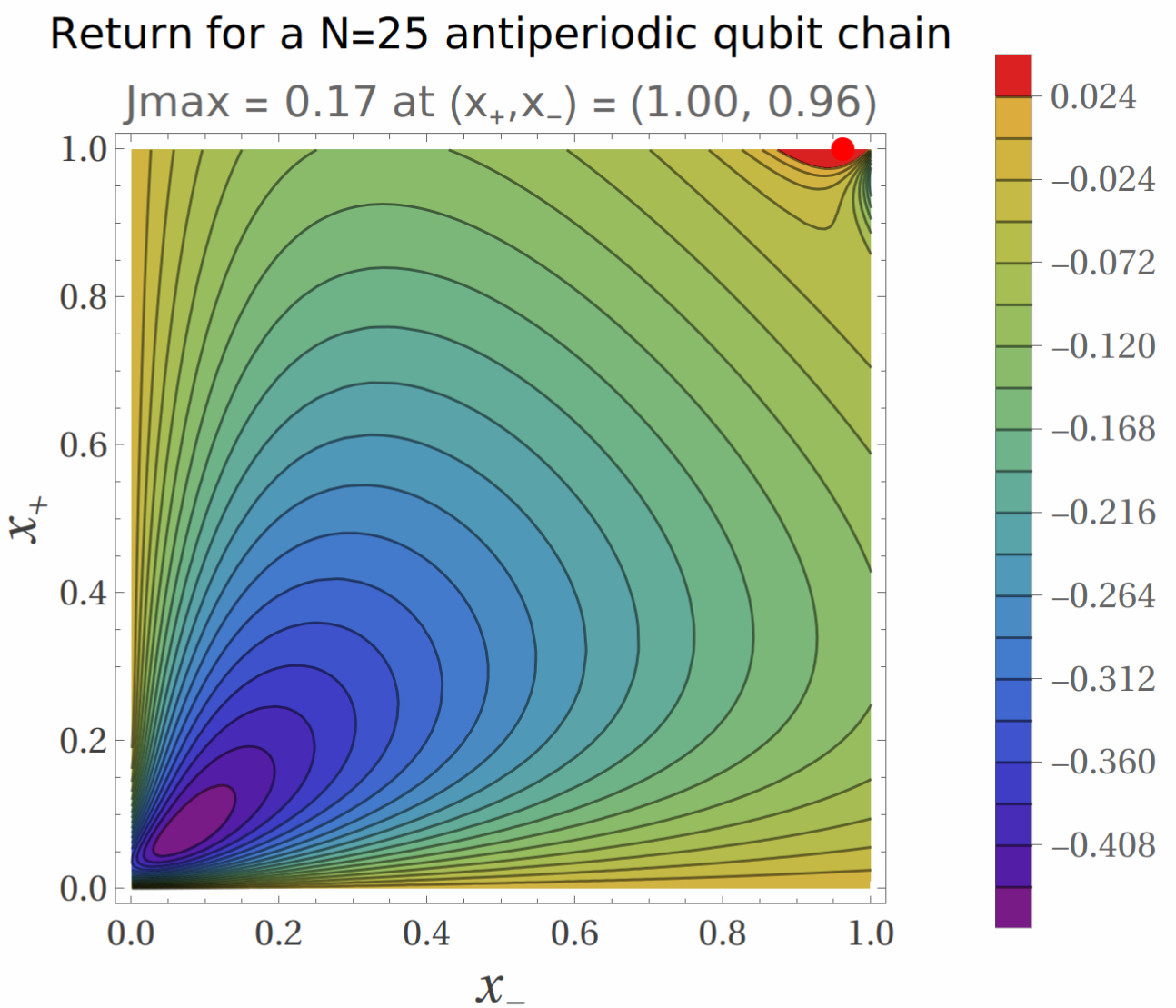}
    \caption{Numerical optimisation of the Expected Return \eqref{eq:Jpidi2pm}}
    \label{fig:Jpi-2levelANTIPER}
\end{figure}

  The numerical optimisation of \eqref{eq:Jpidi2pm} over $(x_+,x_-)\in[0,1]\times[0,1]$ shows that the unique absolute maximum is reached at $(x_+,x_-)=(x_+^{\textrm{max}},x_-^{\textrm{max}})$ with $x_+^{\textrm{max}}=1$ and $x_-^{\textrm{max}}$ approaching $1$ as $N$ grows (Figure \ref{fig:Jpi-2levelANTIPER}). Thus, the optimising policies are those for which (modulo $2\pi$-periodicity) $\theta_+=\pm\frac{\pi}{2}$ and $\theta_-=\pm\frac{\pi}{2}\pm\vartheta$ for positive and suitable small $\vartheta$.

Mathematically, this distinct behaviour is explained by observing how the antiperiodic boundary conditions alter the structure of the sum in \eqref{eq:Jpidi2pm}.
\begin{itemize}
    \item The condition $x_+^{\textrm{max}} = 1$ is not merely an asymptotic trend but an exact feature for $N \geqslant 2$. In the term $(1-x_+)^{(\textsf{p}-\textsf{c})}$, setting $x_+ = 1$ forces the contribution to vanish unless $\textsf{p} = \textsf{c}$. This effectively collapses the inner summation, selecting only those trajectories where the number of sign changes ($\mathsf{c}$) matches the position index ($\textsf{p}$). This simplifies the optimisation landscape by filtering out `partial' transitions, effectively locking the system into a specific coherent subspace.

    \item Unlike the periodic case, where translational invariance imposes $x_+ \approx x_-$, the antiperiodic constraint introduces a frustration that breaks this symmetry. While $x_+$ is pinned to the boundary of the domain, $x_-$ acts as a relaxation parameter. The optimisation decouples the two controls: $x_+$ is maximised to exploit the positive reward terms $\textsf{p}x_+$, while $x_-$ adjusts to mitigate the penalty terms $-(N-\textsf{p})x_-$ and satisfy the probability constraints.

    \item In stark contrast to the periodic case where $x_\pm \sim O(N^{-1})$, here the optimal parameters saturate towards unity. The `bulk' contribution does not vanish as $x \to 1$; rather, it becomes dominant. As $N$ grows, $x_-$ must approach $1$ to maintain the overlap with the trajectories selected by $x_+=1$, ensuring that the probability weight remains concentrated on the high-return paths defined by the boundary condition.
\end{itemize}

From a quantum mechanical perspective, the saturation $x_\pm \to 1$ implies that the optimal policies involve rotation angles $\theta_\pm$ close to $\pi/2$. This corresponds to a `hard' driving regime, effectively implementing a sequence of $\pi$-pulses (spin flips) or `bang-bang' controls. This is the precise antithesis of the strategy observed in the periodic chain. Here the optimal policy switches to a regime of `maximal intervention', where the system is vigorously driven (flipped) at every step to align with the antiperiodic constraint while maximising the accumulated reward.


 \section{Qutrits extension}\label{sec:Qtrits}

The extension from qubits to qutrits represents a natural progression in the present class of QRL models. It maintains analytical tractability while significantly enriching the physical and mathematical structure. Three-level quantum systems are ubiquitous in atomic physics, quantum optics, and solid-state systems, making this extension both theoretically motivated and experimentally relevant.

\subsection{The model}\label{sec:modelQtrits}

 One now has:
 \begin{equation}
  d=\mathrm{dim}\cH=3\,,\qquad \cH=\mathbb{C}^3\,,
 \end{equation}
 \begin{equation}\label{eq:states3}
  \begin{split}
   & \qquad\qquad \mathcal{S}\,=\,\big\{|0\rangle,|1\rangle,|2\rangle\big\}\,, \\
   & \textrm{with }\quad |0\rangle:=\begin{pmatrix} 1 \\ 0 \\ 0 \end{pmatrix}\,,\quad
   |1\rangle:=\begin{pmatrix} 0 \\ 1 \\ 0 \end{pmatrix}\,,\quad
   |2\rangle:=\begin{pmatrix} 0 \\ 0 \\ 1 \end{pmatrix}\,,
  \end{split}
 \end{equation}
 \begin{equation}\label{eq:Hamilt3}
  \begin{split}
   & H\,=\,E_0|0\rangle\langle 0|+E_1|1\rangle\langle 1|+E_2|2\rangle\langle 2| \\
   & \qquad \textrm{with }\quad E_0\,<\,E_1\,<\,E_2\,,
  \end{split}
 \end{equation}
 and, without loss of generality,
 \begin{equation}\label{eq:Hamilt3bis}
  E_0=0\,,\quad E_1=\varepsilon\in(0,1)\,,\quad E_2=1\,.
 \end{equation}

 \medskip

 \textbf{Policies.} As for the unitaries implementing the policies of the QRL  model, rather than considering the full $U(3)$ group, it is meaningful to restrict to physically motivated two-level subspace controls, representing realistic experimental constraints where control fields couple only specific level pairs.

 For concreteness, select a ladder structure with
 \begin{itemize}
  \item  a $|0\rangle\leftrightarrow |1\rangle$ coupling when acting on $|0\rangle$,
  \item  a $|1\rangle\leftrightarrow |2\rangle$ coupling when acting on $|1\rangle$,
  \item  a $|2\rangle\leftrightarrow |0\rangle$ coupling when acting on $|2\rangle$,
 \end{itemize}
 namely
\begin{equation}\label{eq:unitaries3-1}
\begin{split}
  |0\rangle\,\mapsto\,\pi_0(\theta_0)|0\rangle\,, \quad & \pi_0(\theta)\,:=\,
  \begin{pmatrix}
   \cos\theta & -\sin\theta & 0 \\
   \sin\theta & \cos\theta & 0 \\
   0 & 0 & 1
  \end{pmatrix}, \\
 |1\rangle\,\mapsto\,\pi_1(\theta_1)|1\rangle\,, \quad & \pi_1(\theta)\,:=\,
  \begin{pmatrix}
   1 & 0 & 0 \\
   0 & \cos\theta & -\sin\theta \\
   0 & \sin\theta & \cos\theta
  \end{pmatrix}, \\
 |2\rangle\,\mapsto\,\pi_2(\theta_2)|2\rangle\,, \quad & \pi_2(\theta)\,:=\,
  \begin{pmatrix}
   \cos\theta & 0 & -\sin\theta \\
   0 & 1 & 0 \\
   \sin\theta & 0 & \cos\theta
  \end{pmatrix}.
\end{split}
\end{equation}
Thus, a policy $\pi$ is now a triplet
\begin{equation}
 \pi\,\equiv\,\big(\pi_0(\theta_0),\pi_1(\theta_1),\pi_2(\theta_2)\big)\,\equiv\,(\theta_0,\theta_1,\theta_2)\,.
\end{equation}

 Each class of unitaries represents a control field (laser, microwave, magnetic field) with specific selection rules determining allowed transitions. For instance, in alkali atoms in magnetic fields $|0\rangle$, $|1\rangle$, $|2\rangle$ correspond to different Zeeman sub-levels, the unitaries \eqref{eq:unitaries3-1} represent, respectively, $\sigma^+$, $\sigma^-$, and $\pi$ polarized light, and the energy structure is controlled by external magnetic field strength \cite{Cohen-Tannoudji_QM}. For trapped ions, $|0\rangle$, $|1\rangle$, $|2\rangle$ represent an electronic ground state and two metastable excited states, with laser-driven transitions with precise frequency control \cite{Leibfried-trappedions2003}. For superconducting transmons, $|0\rangle$, $|1\rangle$, $|2\rangle$ are charge states in a Cooper pair box and microwave pulses implement unitary rotations \cite{ClarkWilhelm2008,Wendin_2017}. In quantum dots the above tripartite system models a spin-1 system with three $m_s$ states, with external electric and magnetic field control, and tunable energy splittings via gate voltages \cite{Reimann-Manninen_Qdots2002}.

\medskip
\textbf{Trajectories.} A trajectory $\tau$ with $N+1$ states has $\textsf{n}_j$ vectors of type $|j\rangle$, with $\textsf{n}_0+\textsf{n}_1+\textsf{n}_2=N+1$. Consider for concreteness RL to control an overall transition from $|0\rangle$ to $|2\rangle$. Lowest-to-highest-energy-level trajectories of length $N$ (i.e., initial state $|0\rangle$, final state $|2\rangle$) constitute the set $\Gamma_{N,\pi}$, where
\begin{equation}\label{eq:njconstraints3}
 \begin{split}
   \textsf{n}_0\,&\in\,\{1,\dots,N-1\}\,, \\
   \textsf{n}_2\,&\in\,\{1,\dots,N-\textsf{n}_0\}\,, \\
   \textsf{n}_1\,&=\,N+1-\textsf{n}_0-\textsf{n}_2\,.
 \end{split}
\end{equation}
(Observe that the constraint $\textsf{n}_0+\textsf{n}_2\leqslant N$ ensures that $\textsf{n}_1\geqslant 1$, guaranteeing at least one visit to the intermediate level.)

Indicating, analogously to \eqref{eq:trajparamabcd},
\begin{equation}
 \mathsf{a}_{ij}\,:=\,\#[ij] \qquad \big(\textrm{number of transitions $|i\rangle\xrightarrow{\pi_i}\pi_i(\theta_i)|i\rangle\xrightarrow{\textrm{(collapse)}}|j\rangle$}\big)\,
\end{equation}
one has
\begin{equation}\label{eq:constr3-1}
\mathsf{a}_{00} + \mathsf{a}_{01} + \mathsf{a}_{02} + \mathsf{a}_{10} + \mathsf{a}_{11} + \mathsf{a}_{12} + \mathsf{a}_{20} + \mathsf{a}_{21} + \mathsf{a}_{22} = N
\end{equation}
with outgoing constraints (accounting for fixed start at $|0\rangle$ and end at $|2\rangle$)
\begin{equation}\label{eq:constr3-2}
\begin{split}
\text{(from $|0\rangle$)} \quad &\mathsf{a}_{00} + \mathsf{a}_{01} + \mathsf{a}_{02} = \textsf{n}_0\,, \\
\text{(from $|1\rangle$)} \quad &\mathsf{a}_{10} + \mathsf{a}_{11} + \mathsf{a}_{12} = \textsf{n}_1\,, \\
\text{(from $|2\rangle$)} \quad &\mathsf{a}_{20} + \mathsf{a}_{21} + \mathsf{a}_{22} = \textsf{n}_2 - 1\,,
\end{split}
\end{equation}
and incoming constraints
\begin{equation}\label{eq:constr3-3}
\begin{split}
\text{(to $|0\rangle$)} \quad &\mathsf{a}_{00} + \mathsf{a}_{10} + \mathsf{a}_{20} = \textsf{n}_0 - 1\,, \\
\text{(to $|1\rangle$)} \quad &\mathsf{a}_{01} + \mathsf{a}_{11} + \mathsf{a}_{21} = \textsf{n}_1\,, \\
\text{(to $|2\rangle$)} \quad &\mathsf{a}_{02} + \mathsf{a}_{12} + \mathsf{a}_{22} = \textsf{n}_2\,.
\end{split}
\end{equation}

It is worth remarking that, in view of \eqref{eq:constr3-2}-\eqref{eq:constr3-3}, subtracting the incoming constraint from the outgoing constraint for each state yields:
\begin{equation}\label{eq:flowbalance3}
\begin{split}
\text{State } |0\rangle: \quad &(\mathsf{a}_{00} + \mathsf{a}_{01} + \mathsf{a}_{02}) - (\mathsf{a}_{00} + \mathsf{a}_{10} + \mathsf{a}_{20}) = \textsf{n}_0 - (\textsf{n}_0 - 1) = +1\,, \\
\text{State } |1\rangle: \quad &(\mathsf{a}_{10} + \mathsf{a}_{11} + \mathsf{a}_{12}) - (\mathsf{a}_{01} + \mathsf{a}_{11} + \mathsf{a}_{21}) = \textsf{n}_1 - \textsf{n}_1 = 0\,, \\
\text{State } |2\rangle: \quad &(\mathsf{a}_{20} + \mathsf{a}_{21} + \mathsf{a}_{22}) - (\mathsf{a}_{02} + \mathsf{a}_{12} + \mathsf{a}_{22}) = (\textsf{n}_2 - 1) - \textsf{n}_2 = -1\,.
\end{split}
\end{equation}
Simplifying, the self-loop terms $\mathsf{a}_{ii}$ cancel, yielding:
\begin{equation}\label{eq:netflow3}
\begin{split}
|0\rangle: \quad &\mathsf{a}_{01} + \mathsf{a}_{02} - \mathsf{a}_{10} - \mathsf{a}_{20} = +1\,, \\
|1\rangle: \quad &\mathsf{a}_{10} + \mathsf{a}_{12} - \mathsf{a}_{01} - \mathsf{a}_{21} = 0\,, \\
|2\rangle: \quad &\mathsf{a}_{20} + \mathsf{a}_{21} - \mathsf{a}_{02} - \mathsf{a}_{12} = -1\,.
\end{split}
\end{equation}
These equations express a `flow conservation principle': each state has a net flow (outgoing minus incoming transitions, excluding self-loops). State $|0\rangle$ has net outflow $+1$ (the trajectory starts there), state $|2\rangle$ has net inflow $-1$ (the trajectory ends there), and state $|1\rangle$ has zero net flow (pure transit state). This is analogous to Kirchhoff's current law in circuit theory: the algebraic sum of net flows must equal zero ($+1 + 0 - 1 = 0$), ensuring conservation of probability flow through the quantum system.

Notably, the self-loop terms $\mathsf{a}_{00}$, $\mathsf{a}_{11}$, $\mathsf{a}_{22}$ do not appear in \eqref{eq:netflow3} because they represent internal circulation within each state, contributing equally to both outgoing and incoming counts and thus canceling in the net flow calculation.

\medskip

\textbf{Allowed transitions.} Formula \eqref{eq:probtraj} requires the evaluation of factors of type
\begin{equation}
 |\langle j | \pi_i(\theta_i) | i \rangle |^2 \,=:\,P_{ij}
\end{equation}
for each transition $|i\rangle\xrightarrow{\pi_i}\pi_i(\theta_i)|i\rangle\xrightarrow{\textrm{(collapse)}}|j\rangle$. These $P_{ij}$'s are computed directly by means of \eqref{eq:states3} and \eqref{eq:unitaries3-1}; they are assembled into the transition probability matrix
\begin{equation}\label{eq:3transP}
\mathbf{P} \,:=\,(P_{ij})\,=\, \begin{pmatrix}
\cos^2\theta_0 & \sin^2\theta_0 & 0 \\
0 & \cos^2\theta_1 & \sin^2\theta_1 \\
\sin^2\theta_2 & 0 & \cos^2\theta_2
\end{pmatrix}.
\end{equation}
The zero-entries in $\mathbf{P}$ correspond to transitions between orthogonal states and are an obvious consequence of the special rotations \eqref{eq:unitaries3-1}. They act as super-selection rules for the model. The fact that $P_{10}=P_{21}=P_{02}=0$ results in the additional constraint
\begin{equation}\label{eq:3additionalConstr}
 \mathsf{a}_{02}\,=\,\mathsf{a}_{10}\,=\,\mathsf{a}_{21}\,=\,0\,.
\end{equation}

\medskip

\textbf{Independent degrees of freedom.} Eventually, one degree of freedom only remains in the transition counts $\mathsf{a}_{ij}$. This is straightforward to determine by direct identification of the independent constraints. Actually, a more systematic constraint count elucidates the natural graph-like structure  of the above relationships. Indeed, the three outgoing constraints
\eqref{eq:constr3-2} are mutually independent and the three net-flow equations
\eqref{eq:netflow3} are precisely Kirchhoff's node law for the directed graph
$\mathcal{G}_3$ with node set $\{|0\rangle,|1\rangle,|2\rangle\}$ and edge set
$\{\mathsf{a}_{ij}\,|\,i\neq j\}$ (Figure \ref{fig:graph-qutrits}). The coefficient matrix of \eqref{eq:netflow3} is the
incidence matrix of $\mathcal{G}_3$. Since $\mathcal{G}_3$ is connected, its incidence
matrix has rank $3-1=2$ \cite[Theorem~2.3]{Bapat2014}, so exactly two of the three
equations \eqref{eq:netflow3} are independent. They are moreover independent of
\eqref{eq:constr3-2}, since each outgoing constraint involves the self-loop term
$\mathsf{a}_{ii}$, which does not appear in any of the net-flow equations
\eqref{eq:netflow3}. Finally, the ladder topology of the policy unitaries
\eqref{eq:unitaries3-1} enforces three forbidden transitions \eqref{eq:3additionalConstr} (zero probability). In total, one has $3+2+3=8$ independent conditions on the $9$ transition counts
$\mathsf{a}_{ij}$, leaving exactly $9-8=1$ degree of freedom.

\begin{figure}[t!]
\centering
\begin{tikzpicture}[>=stealth, semithick]
  \node[draw, circle, inner sep=5pt] (n0) at ( 0.0,  2.2) {$|0\rangle$};
  \node[draw, circle, inner sep=5pt] (n1) at (-1.9,  0.0) {$|1\rangle$};
  \node[draw, circle, inner sep=5pt] (n2) at ( 1.9,  0.0) {$|2\rangle$};
  \draw[->] (n0) to[bend right=20] (n1);
  \draw[->] (n1) to[bend right=20] (n0);
  \draw[->] (n1) to[bend right=20] (n2);
  \draw[->] (n2) to[bend right=20] (n1);
  \draw[->] (n2) to[bend right=20] (n0);
  \draw[->] (n0) to[bend right=20] (n2);
\end{tikzpicture}
\qquad\qquad
\begin{tikzpicture}[>=stealth, semithick]
  \node[draw, circle, inner sep=5pt] (n0) at ( 0.0,  2.2) {$|0\rangle$};
  \node[draw, circle, inner sep=5pt] (n1) at (-1.9,  0.0) {$|1\rangle$};
  \node[draw, circle, inner sep=5pt] (n2) at ( 1.9,  0.0) {$|2\rangle$};
  \draw[->] (n0) to[bend right=20] node[left,  midway] {$1+\mathsf{c}\;$} (n1);
  \draw[->] (n1) to[bend right=20] node[below, midway] {$1+\mathsf{c}$} (n2);
  \draw[->] (n2) to[bend right=20] node[right, midway] {$\;\mathsf{c}$}   (n0);
\end{tikzpicture}
\caption{Left: the directed graph $\mathcal{G}_3$ on node set
$\{|0\rangle,|1\rangle,|2\rangle\}$ with all six off-diagonal edges, whose incidence
matrix has rank $3-1=2$ \cite[Theorem~2.3]{Bapat2014}. Right: the
reduced transition graph $\widetilde{\mathcal{G}}_3$, obtained from $\mathcal{G}_3$
by retaining only the three edges $|i\rangle\to|j\rangle$ with nonzero probability
after imposing the forbidden-transition constraints \eqref{eq:3additionalConstr};
edge labels indicate the transition multiplicities $\mathsf{a}_{ij}$ from
\eqref{eq:3ParamSIMPL}, for $\mathsf{c}\geqslant 1$. The right graph illustrates the
$1+\mathsf{c}$ cyclic runs $|0\rangle\to|1\rangle\to|2\rangle$ punctuated by
$\mathsf{c}$ feedback returns $|2\rangle\to|0\rangle$; for $\mathsf{c}=0$ the
feedback arc is absent and $\widetilde{\mathcal{G}}_3$ reduces to the directed path
$|0\rangle\to|1\rangle\to|2\rangle$.}
\label{fig:graph-qutrits}
\end{figure}

One may choose, say, the feedback count $\mathsf{a}_{20}$ as the independent parameter:
\begin{equation}\label{eq:c20def}
\mathsf{c} \,:=\, \mathsf{a}_{20}\,,
\end{equation}
and the remaining variables are determined by
\begin{equation}\label{eq:3ParamSIMPL}
\begin{aligned}
\mathsf{a}_{00} \,&=\, \textsf{n}_0 - 1 - \mathsf{c} \,, & \mathsf{a}_{01} \,&=\, 1 + \mathsf{c} \,, \\
\mathsf{a}_{02} \,&=\, 0 \,, & \mathsf{a}_{10} \,&=\, 0 \,, \\
\mathsf{a}_{11} \,&=\, \textsf{n}_1 - 1 - \mathsf{c} \,, & \mathsf{a}_{12} \,&=\, 1 + \mathsf{c} \,, \\
\mathsf{a}_{20} \,&=\, \mathsf{c} \,, & \mathsf{a}_{21} \,&=\, 0 \,, \\
\mathsf{a}_{22} \,&=\, \textsf{n}_2 - 1 - \mathsf{c} \,. & &
\end{aligned}
\end{equation}

For a given $N$, the admissible ranges for $\textsf{n}_0, \textsf{n}_1, \textsf{n}_2$, $\mathsf{c}$ are given by
\begin{equation}\label{eq:3constraintsn}
\begin{split}
\textsf{n}_0\,&\in\, \{1,\dots,N-1\}\,, \\
\textsf{n}_2\,&\in\, \{1,\dots,N-\textsf{n}_0\}\,, \\
\textsf{n}_1\,&=\, N+1-\textsf{n}_0-\textsf{n}_2\,, \\
\mathsf{c} \,&\in\, \{0, \dots, \min\{\textsf{n}_0, \textsf{n}_1, \textsf{n}_2\} - 1\}\,.
\end{split}
\end{equation}
The ranges for $\textsf{n}_0, \textsf{n}_1, \textsf{n}_2$ were determined in \eqref{eq:njconstraints3}. The admissible range for $\mathsf{c}$ is determined by requiring all $\mathsf{a}_{ij} \geqslant 0$.

\medskip

\textbf{Trajectory multiplicity.} The forbidden transitions \eqref{eq:3additionalConstr} impose a rigid sequential
structure: a valid trajectory must consist of $1+\mathsf{c}$ distinct `runs' of visits
to state $|0\rangle$, followed by state $|1\rangle$, followed by state $|2\rangle$,
cycling $1+\mathsf{c}$ times (with the last run of $|2\rangle$ terminating the
trajectory). The allowed non-self-loop transitions, with their multiplicities from
\eqref{eq:3ParamSIMPL}, are depicted in Figure~\ref{fig:graph-qutrits}. The number of
such valid trajectories is exactly the number of ways to distribute the $\textsf{n}_j$
total visits into these $1+\mathsf{c}$ non-empty groups for each state. Using the
combinatorial composition formula (`stars and bars' counting), one obtains the
closed-form expression
\begin{equation}\label{eq:3trajmult}
\mathcal{N}(N,\textsf{n}_0,\textsf{n}_2,\mathsf{c}) \,=\,
\binom{\textsf{n}_0-1}{\mathsf{c}}\binom{\textsf{n}_1-1}{\mathsf{c}}\binom{\textsf{n}_2-1}{\mathsf{c}}\,.
\end{equation}
(Recall: $\textsf{n}_1= N+1-\textsf{n}_0-\textsf{n}_2$.)
This formula naturally evaluates to zero if $\mathsf{c}>\textsf{n}_j-1$, consistent with
the constraints.

\medskip

\textbf{Trajectory probability.} For a lowest-to-highest-energy-level trajectory $\tau$ with given $N,\textsf{n}_0,\textsf{n}_1,\textsf{n}_2,\mathsf{c}$ and given policy $\pi\equiv(\theta_0,\theta_1,\theta_2)$, the probability \eqref{eq:probtraj} reads
\begin{equation}\label{eq:ProbTraj3}
\begin{split}
P_\pi(\tau) \,&=\, (\cos^2\theta_0)^{\textsf{a}_{00}}(\sin^2\theta_0)^{\textsf{a}_{01}}(\cos^2\theta_1)^{\textsf{a}_{11}}(\sin^2\theta_1)^{\textsf{a}_{12}}(\sin^2\theta_2)^{\textsf{a}_{20}}(\cos^2\theta_2)^{\textsf{a}_{22}} \\
&=\, (\cos^2\theta_0)^{\textsf{n}_0 - 1 - \mathsf{c}}(\sin^2\theta_0)^{1 + \mathsf{c}}(\cos^2\theta_1)^{\textsf{n}_1 - 1 - \mathsf{c}}(\sin^2\theta_1)^{1 + \mathsf{c}} \times \\
& \qquad \times (\sin^2\theta_2)^{\mathsf{c}}(\cos^2\theta_2)^{\textsf{n}_2 - 1 - \mathsf{c}} \\
&=:\,P(N,\textsf{n}_0,\textsf{n}_1,\textsf{n}_2,\mathsf{c},\theta_0,\theta_1,\theta_2)\,,
\end{split}
\end{equation}
where we used \eqref{eq:3ParamSIMPL}.

\medskip

 \textbf{Reward along a trajectory.} The reward structure consists now, for each transition $|i\rangle\xrightarrow{\pi_i}\pi_i(\theta_i)|i\rangle\xrightarrow{\textrm{(collapse)}}|j\rangle$, of energy differences
\begin{equation}
  E_j - \langle\pi_i(\theta_i)i|H|\pi_i(\theta_i)i\rangle\,=:\, R_{ij}
\end{equation}
between post-collapse and pre-collapse states.
Such $R_{ij}$'s are computed directly by means of \eqref{eq:states3}-\eqref{eq:Hamilt3bis} and \eqref{eq:unitaries3-1}; they are assembled into the transition reward matrix
\begin{equation}\label{eq:3rewmatr}
 \mathbf{R} \,:=\,(R_{ij})\,=\,\begin{pmatrix}
-\varepsilon\sin^2\theta_0 & \varepsilon\cos^2\theta_0 & 1 - \varepsilon\sin^2\theta_0 \\
-\varepsilon\cos^2\theta_1 - \sin^2\theta_1 & -(1-\varepsilon)\sin^2\theta_1 & (1-\varepsilon)\cos^2\theta_1 \\
-\cos^2\theta_2 & \varepsilon - \cos^2\theta_2 & \sin^2\theta_2
\end{pmatrix}.
\end{equation}
For a lowest-to-highest-energy-level trajectory $\tau$, plugging the entries of \eqref{eq:3rewmatr} and the transition counts from \eqref{eq:3ParamSIMPL} into the reward sum yields $R_\pi(\tau)$.
Remarkably, the linear dependence on the feedback parameter $\mathsf{c}$ cancels out exactly (the energy cost of the feedback loop balances with the change in self-loops). The final result depends solely on the occupation numbers:
\begin{equation}\label{eq:reward3}
\begin{split}
R_\pi(\tau) \,&=\, \sum_{i,j} \mathsf{a}_{ij} R_{ij} \\
&=\, (\textsf{n}_0 - 1 - \mathsf{c})(-\varepsilon\sin^2\theta_0) + (1 + \mathsf{c})(\varepsilon\cos^2\theta_0) \\
&\quad + (\textsf{n}_1 - 1 - \mathsf{c})(-(1-\varepsilon)\sin^2\theta_1) + (1 + \mathsf{c})((1-\varepsilon)\cos^2\theta_1) \\
&\quad + \mathsf{c}(-\cos^2\theta_2) + (\textsf{n}_2 - 1 - \mathsf{c})(\sin^2\theta_2) \\
&=\, \textsf{n}_0 (-\varepsilon\sin^2\theta_0) + \textsf{n}_1 (-(1-\varepsilon)\sin^2\theta_1) + \textsf{n}_2 (\sin^2\theta_2) + \cos^2\theta_2 \\
&=:\,R_\varepsilon(\textsf{n}_0,\textsf{n}_1,\textsf{n}_2,\theta_0,\theta_1,\theta_2)\,.
\end{split}
\end{equation}

\medskip

\textbf{Expected return formula.} The expected return combines trajectory probabilities, combinatorial weights, and rewards, according to \eqref{eq:Jpidi}, and takes here the form
\begin{equation}\label{eq:Jpidi3}
\begin{split}
  J(\pi) \:&= \sum_{ \substack{ \textsf{n}_0,\textsf{n}_1,\textsf{n}_2,\mathsf{c} \\ \textrm{with } \eqref{eq:3constraintsn} }}   \mathcal{N}(N,\textsf{n}_0,\textsf{n}_1,\textsf{n}_2,\mathsf{c})\,R_\varepsilon(\textsf{n}_0,\textsf{n}_1,\textsf{n}_2,\theta_0,\theta_1,\theta_2) \, P(N,\textsf{n}_0,\textsf{n}_1,\textsf{n}_2,\mathsf{c},\theta_0,\theta_1,\theta_2) \\
  &=: J_\varepsilon(\theta_0,\theta_1,\theta_2)\,.
\end{split}
\end{equation}

\medskip

\textbf{Remarks on the model.}
\begin{enumerate}
 \item[(1)] Meaningful variants include actions that couple the ground state $|0\rangle$ with two excited states $|1\rangle$ and $|2\rangle$ (such as in `$\Lambda$-configurations'), or two ground states coupled to a common excited state (`V-configuration'). Here the standard ladder structure is adjusted so as to couple the highest energy state back to the ground state (the state $|2\rangle$ is rotated into $\cos\theta_2|2\rangle-\sin\theta_2|0\rangle$), a cyclic structure that disturbs the natural `climb-and-stay' strategy of reaching $|2\rangle$ from $|0\rangle$ as soon as possible and maintaining this state until the end of the trajectory, and is therefore expected to produce a richer optimisation scenario.

 \item[(2)] In either case, it is noticeable once again, as with \eqref{eq:Jpidi2bis} or \eqref{eq:Jpidi2pm}, that \eqref{eq:Jpidi3} displays a substantial reduction of computational complexity with $N$, since the sum scales as $O(N^3)$ (summing over $\textsf{n}_0,\textsf{n}_2,\mathsf{c}$ with constraints) as compared to the general $O(3^N)$-scaling of the expression \eqref{eq:Jpidi}. Moreover, the closed-form multiplicity formula \eqref{eq:3trajmult} eliminates the need for combinatorial enumeration, making the computation fully analytical.

 \item[(3)] Formulas \eqref{eq:ProbTraj3}, \eqref{eq:reward3}, and \eqref{eq:Jpidi3} show that the expected return $J(\pi)$ depends smoothly on the three variables $\theta_0,\theta_1,\theta_2$ (and the additional energy parameter $\varepsilon$), which can be taken in $[0,\frac{\pi}{2}]$ without loss of generality. This makes the policy optimisation of $J(\pi)$ quite tractable numerically.
\end{enumerate}

\medskip

\textbf{On the cardinality of a $(N,\textsf{n}_0,\textsf{n}_2,\mathsf{c})$-subclass of trajectories.} A sub-class of trajectories with the same counts $\textsf{n}_0,\textsf{n}_2,\mathsf{c}$ has cardinality \eqref{eq:3trajmult}. There are $O(N^3)$ such sub-classes because each labelling parameter runs on a set of $O(N)$-integers, as indicated in the constraints \eqref{eq:3constraintsn}. This indicates that there must be sub-classes with exponential-in-$N$ cardinality, in order to match the overall amount of $3^{N-1}$ distinct trajectories (having 3 possible distinct states in each of the $N-1$ free nodes of a trajectory). Technically, $3^{N-1}$ counts both the allowed and the forbidden (probability-zero) trajectories, but this does not alter the conclusion. In fact, there are subclasses whose cardinality scales in $N$ exponentially or with power-law. For example, it is not hard to find the following explicit estimates.
\begin{enumerate}
 \item Trajectory sub-classes with
\begin{equation}\label{eq:expsubclass}
\begin{split}
 \textsf{n}_0\,&=\,\alpha\,N\,,\qquad \textsf{n}_2\,=\,\beta\,N\,,\qquad \mathsf{c}\,=\,\gamma\,N\,, \\
 \alpha,\beta,\gamma &>0\,,\qquad \alpha+\beta<1\,, \qquad \gamma<\min\{\alpha,\beta,1-\alpha-\beta\}
\end{split}
\end{equation}
 (which is compatible with \eqref{eq:3constraintsn}, up to $O(1)$ integer correction), have cardinality
 \begin{equation}\label{eq:cardexp}
  \mathcal{N}(N,\textsf{n}_0,\textsf{n}_2,\mathsf{c}) \stackrel{(N\to\infty)}{=}
\frac{C(\alpha,\beta,\gamma)}{N^{3/2}}\;e^{\,N r(\alpha,\beta,\gamma)}\,\big(1+O(N^{-1})\big)
 \end{equation}
 for certain constants $C(\alpha,\beta,\gamma),r(\alpha,\beta,\gamma)>0$.
 \item Trajectory subclasses in which the count $\mathsf{c}$ is kept $N$-independent are only polynomially-in-$N$ populated. If
 \begin{equation}\label{eq:polysubclassA}
 \textsf{n}_0\,,\;\textsf{n}_2\,,\;\mathsf{c}\, = \,N\text{-independent integers}
\end{equation}
(which is compatible with the ranges in \eqref{eq:3constraintsn}), then
 \begin{equation}\label{eq:cardA}
 \mathcal{N}(N,\textsf{n}_0,\textsf{n}_2,\mathsf{c})\,\stackrel{(N\to\infty)}{=}\,A(\textsf{n}_0,\textsf{n}_2,\mathsf{c})\, N^{\mathsf{c}}\,+\,O(N^{\mathsf{c}-1})
\end{equation}
 for certain $A(\textsf{n}_0,\textsf{n}_2,\mathsf{c})>0$.
 \item Trajectory sub-classes with
\begin{equation}\label{eq:polysubclassB}
 \begin{split}
   \textsf{n}_0\,=\,\alpha\,N\,,\qquad \textsf{n}_2\,&=\,\beta\,N\,,\qquad \mathsf{c}\,=\, N\text{-independent integer}\,, \\
    \alpha\,,\beta\,&>\,0\,,\qquad \alpha+\beta\,<\,1\,,
 \end{split}
\end{equation}
 (which is compatible with \eqref{eq:3constraintsn}) have cardinality
 \begin{equation}\label{eq:cardB}
 \mathcal{N}(N,\textsf{n}_0,\textsf{n}_2,\mathsf{c})\,\stackrel{(N\to\infty)}{=}\,B(\alpha,\beta,\mathsf{c})\,N^{3\mathsf{c}}\,+\,O(N^{3\mathsf{c}-1})
\end{equation}
for certain $B(\alpha,\beta,\mathsf{c})>0$.
 \item The boundary case $\mathsf{c}=0$ reduces all three binomials in \eqref{eq:3trajmult} to unity, whence
\begin{equation}\label{eq:cardC0}
 \mathcal{N}(N,\textsf{n}_0,\textsf{n}_2,0)\,\equiv\,1
\end{equation}
exactly for every admissible $(\textsf{n}_0,\textsf{n}_2)$: such a sub-class contains the single ``climb-and-stay'' trajectory $|0\rangle^{\textsf{n}_0}|1\rangle^{\textsf{n}_1}|2\rangle^{\textsf{n}_2}$, consisting of one contiguous block per energy level and no feedback excursion.
\end{enumerate}

\subsection{Optimisation with policies of common rotation}\label{sec:optim3commonrot}

This is the scenario with constraint $\theta_0=\theta_1=\theta_2=:\theta$: the rotation angle is the same for all trajectory's states. This reduces the optimisation to just $\theta\in[0,\frac{\pi}{2}]$. Under this constraint, \eqref{eq:ProbTraj3} and \eqref{eq:reward3} simplify, respectively, to
\begin{equation}\label{eq:Ppi3commontheta}
\begin{split}
  P(N,\textsf{n}_0,\textsf{n}_2,\mathsf{c},\theta)\,:=&\,(\cos^2\theta)^{\textsf{a}_{00}+\textsf{a}_{11}+\textsf{a}_{22}}(\sin^2\theta)^{\textsf{a}_{01}+\textsf{a}_{12}+\textsf{a}_{20}} \\
  =&\,(\cos^2\theta)^{N-3\mathsf{c}-2}(\sin^2\theta)^{3\mathsf{c}+2}
\end{split}
\end{equation}
and
\begin{equation}\label{eq:Rpi3commontheta}
\begin{split}
R_\varepsilon(\textsf{n}_0,\textsf{n}_2,\theta) \, &= \, \textsf{n}_0 (-\varepsilon\sin^2\theta) + \textsf{n}_1 (-(1-\varepsilon)\sin^2\theta) + \textsf{n}_2 (\sin^2\theta) + \cos^2\theta \\
&= \, \big(N\varepsilon - N - 2\textsf{n}_0\varepsilon + \textsf{n}_0 - \textsf{n}_2\varepsilon + 2\textsf{n}_2 + \varepsilon - 1\big)\sin^2\theta + \cos^2\theta
\end{split}
\end{equation}
(recall again that $\textsf{n}_1= N+1-\textsf{n}_0-\textsf{n}_2$).

The expected return \eqref{eq:Jpidi3} takes the form $J_\varepsilon(x)\equiv J(\pi)$ in the variable
\begin{equation}
 x\,:=\,\sin^2\theta\,\in\,[0,1]
\end{equation}
given by
\begin{equation}\label{eq:Jpi3commontheta}
\begin{split}
J_\varepsilon(x)\,&=\sum_{\textsf{n}_0=1}^{N-1}\sum_{\textsf{n}_2=1}^{N-\textsf{n}_0}\sum_{\mathsf{c}=0}^{\min\{\textsf{n}_0,N+1-\textsf{n}_0-\textsf{n}_2,\textsf{n}_2\}-1} \binom{\textsf{n}_0-1}{\mathsf{c}}\binom{N-\textsf{n}_0-\textsf{n}_2}{\mathsf{c}}\binom{\textsf{n}_2-1}{\mathsf{c}} \\
&\qquad\times \bigg[ \Big( \varepsilon(N+1 - 2\textsf{n}_0-\textsf{n}_2) - (N+2-\textsf{n}_0 - 2\textsf{n}_2) \Big)x + 1 \bigg] \: (1-x)^{N-3\mathsf{c}-2}x^{3\mathsf{c}+2}\,,
\end{split}
\end{equation}
as follows by plugging \eqref{eq:3trajmult}, \eqref{eq:Ppi3commontheta}, and \eqref{eq:Rpi3commontheta} into \eqref{eq:Jpidi3}.

The latter formula has a symmetry that is explicitly seen by restoring $\textsf{n}_1$ from $\textsf{n}_0+\textsf{n}_1+\textsf{n}_2=N+1$. One then re-writes
\begin{equation}\label{eq:Jpi3simplified-2}
\begin{split}
J_\varepsilon(x)\,&=\sum_{\textsf{n}_0=1}^{N-1}\sum_{\textsf{n}_2=1}^{N-\textsf{n}_0}\sum_{\mathsf{c}=0}^{\min\{\textsf{n}_0,\textsf{n}_1,\textsf{n}_2\}-1} \binom{\textsf{n}_0-1}{\mathsf{c}}\binom{\textsf{n}_1-1}{\mathsf{c}}\binom{\textsf{n}_2-1}{\mathsf{c}} \\
&\qquad\times \bigg[ \Big( \varepsilon(\textsf{n}_1 - \textsf{n}_0) + (\textsf{n}_0 + 2\textsf{n}_2 - N - 2) \Big)x + 1 \bigg] \: (1-x)^{N-3\mathsf{c}-2}x^{3\mathsf{c}+2}\,, \\
&\textrm{with }\;\textsf{n}_0+\textsf{n}_1+\textsf{n}_2=N+1\,.
\end{split}
\end{equation}
Now, swapping $\textsf{n}_0\leftrightarrow \textsf{n}_1$ leaves the product of the binomials unchanged, and inverts the sign of the $\varepsilon(\textsf{n}_1 - \textsf{n}_0)$ term, and there are precisely as many summands with $ \textsf{n}_1 > \textsf{n}_0$ as with $ \textsf{n}_1 < \textsf{n}_0$. The net result is a cancellation of all the $\varepsilon$-dependent terms and a final formula that is independent of $\varepsilon$:
\begin{equation}\label{eq:Jpi3simplified-3}
\begin{split}
J_\varepsilon(x)\,&=\sum_{\textsf{n}_0=1}^{N-1}\sum_{\textsf{n}_2=1}^{N-\textsf{n}_0}\sum_{\mathsf{c}=0}^{\min\{\textsf{n}_0,N+1-\textsf{n}_0-\textsf{n}_2,\textsf{n}_2\}-1} \binom{\textsf{n}_0-1}{\mathsf{c}}\binom{N-\textsf{n}_0-\textsf{n}_2}{\mathsf{c}}\binom{\textsf{n}_2-1}{\mathsf{c}} \\
&\qquad\times \big[ (\textsf{n}_0 + 2\textsf{n}_2 - N - 2)x + 1 \big] \: (1-x)^{N-3\mathsf{c}-2}x^{3\mathsf{c}+2}\,\equiv\,J(x)\,.
\end{split}
\end{equation}

  \begin{figure}[t!]
    \centering
   \includegraphics[width=7cm]{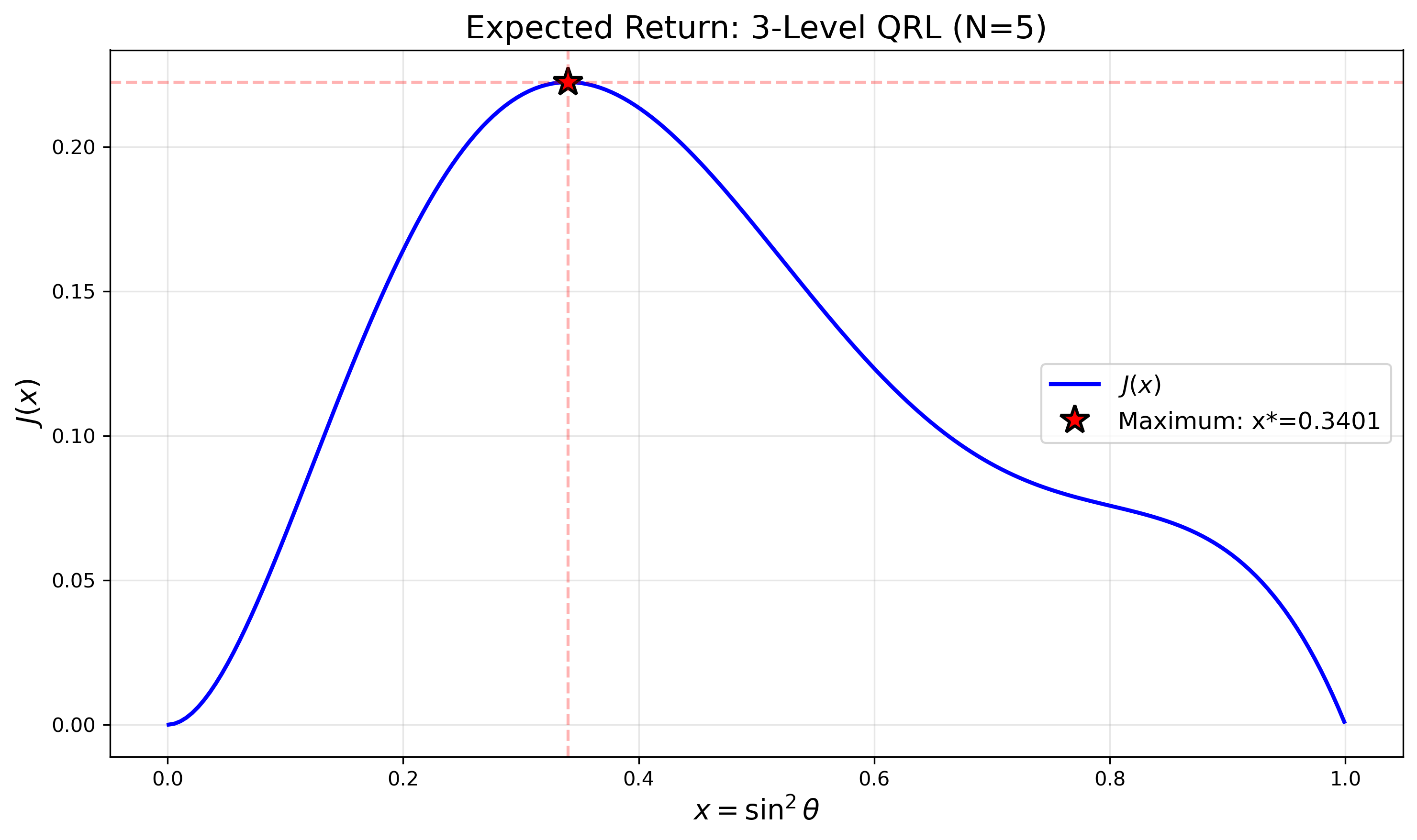}\quad
        \includegraphics[width=7cm]{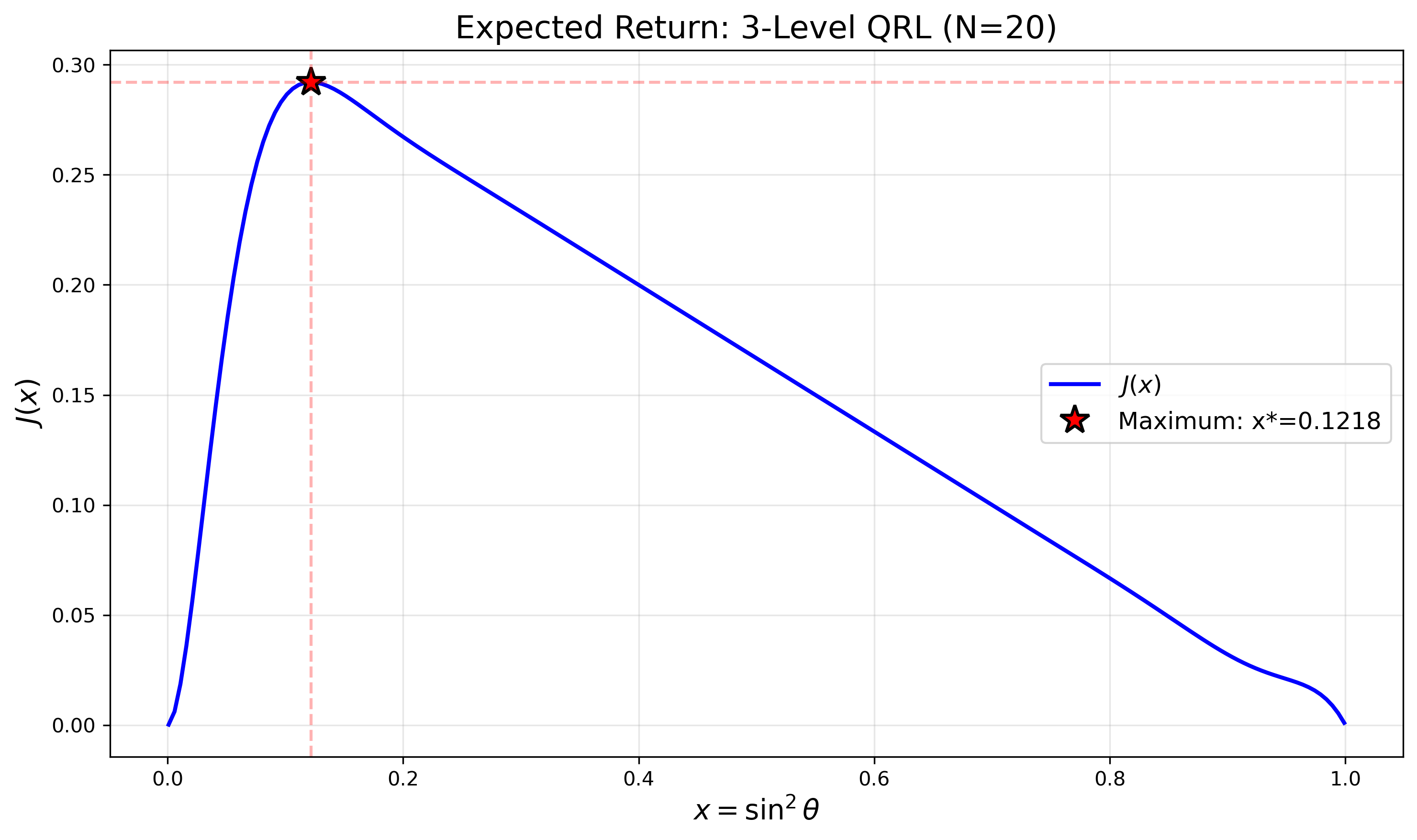}\quad
        \includegraphics[width=7cm]{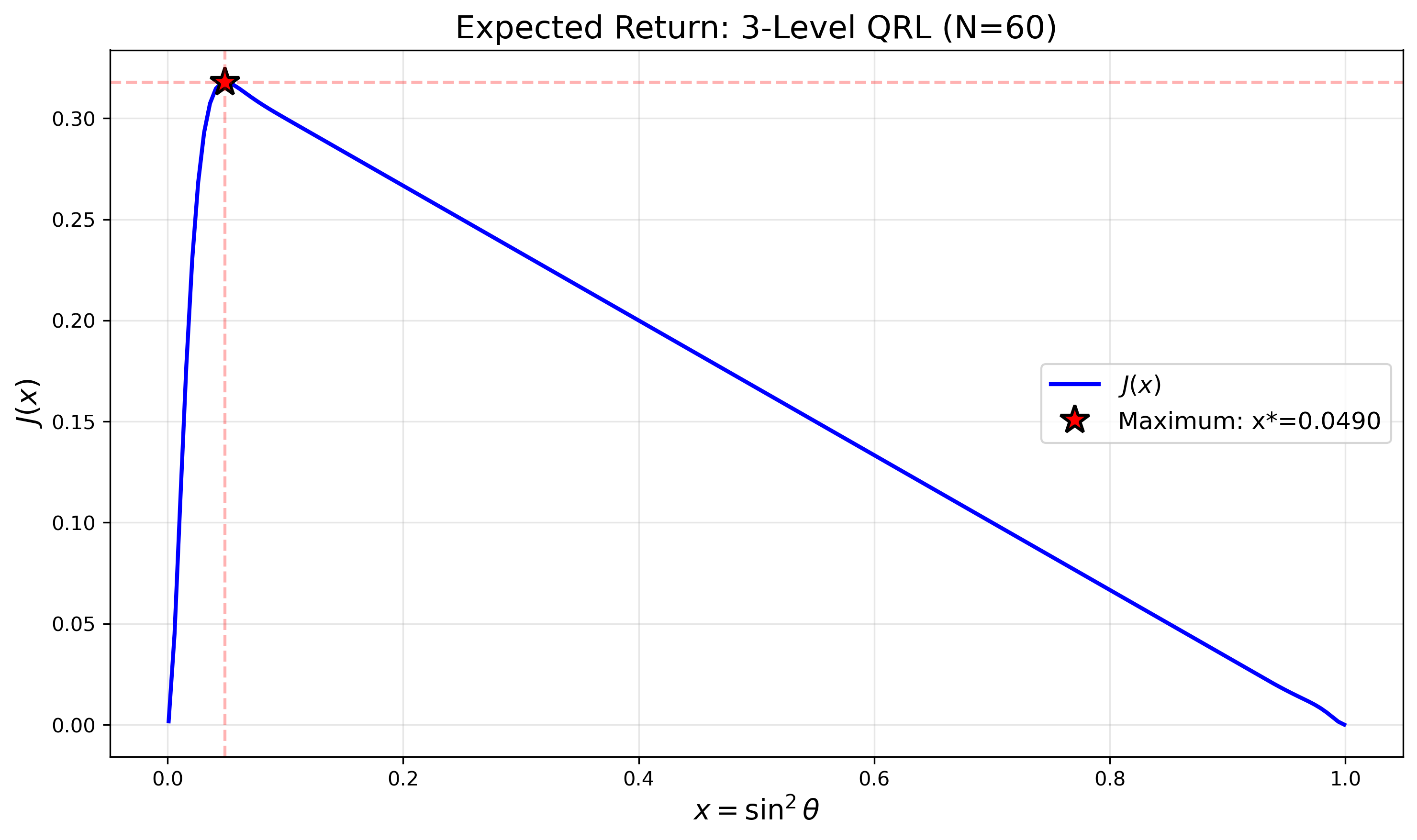}
    \caption{Numerical optimisation of the Expected Return \eqref{eq:Jpi3simplified-3}}
    \label{fig:Jpi-3levelCOMMON}
\end{figure}

 Figure \ref{fig:Jpi-3levelCOMMON} displays different cases of numerical optimisation of \eqref{eq:Jpi3simplified-3}.\begin{itemize}
    \item The return vanishes at the boundaries. At $x=0$ ($\theta=0$), the evolution is the identity, freezing the system in the initial state $|0\rangle$. At $x=1$ ($\theta=\pm\pi/2$), the dynamics become a deterministic cyclic permutation $|0\rangle \to |1\rangle \to |2\rangle \to |0\rangle \dots$, preventing stable accumulation of probability in the target state (inability of the system to remain in the target state $|2\rangle$ once it arrives there).

    \item As $N$ increases, the optimal control $x^*$ shifts towards zero. Short trajectories require strong pulses to populate the target state quickly, whereas long trajectories favour a `weak coupling' regime. In this limit, small rotations allow a steady, diffusive-like probability transfer while suppressing the rapid oscillatory behaviour that would depopulate the target state.
\end{itemize}

\subsection{Optimisation with policies of multiple independent rotations}\label{sec:optim3indeprot}

The general scenario consists of the optimisation of the expected return \eqref{eq:Jpidi3} as a function $J_\varepsilon(x,y,z)$ of the three independent variables
\begin{equation}
 x\,:=\,\sin^2\theta_0\,,\qquad y\,:=\,\sin^2\theta_1\,,\qquad z\,:=\,\sin^2\theta_2
\end{equation}
ranging in the cube $[0,1]\times[0,1]\times[0,1]$.

By substituting the multiplicity \eqref{eq:3trajmult}, the trajectory probability \eqref{eq:ProbTraj3}, and the reward \eqref{eq:reward3} into the general expectation formula \eqref{eq:Jpidi3}, one obtains
\begin{equation}\label{eq:Jpi3threethetas}
 \begin{split}
 J_\varepsilon(x,y,z)\,&=\sum_{\textsf{n}_0=1}^{N-1}\sum_{\textsf{n}_2=1}^{N-\textsf{n}_0}\sum_{\mathsf{c}=0}^{\min\{\textsf{n}_0,N+1-\textsf{n}_0-\textsf{n}_2,\textsf{n}_2\}-1} \binom{\textsf{n}_0-1}{\mathsf{c}}\binom{N-\textsf{n}_0-\textsf{n}_2}{\mathsf{c}}\binom{\textsf{n}_2-1}{\mathsf{c}} \\
 &\qquad\times \big[1 -\varepsilon \textsf{n}_0 x - (1-\varepsilon)(N+1-\textsf{n}_0-\textsf{n}_2) y + (\textsf{n}_2 - 1) z \big] \\
 &\qquad\times (1-x)^{\textsf{n}_0 - 1 - \mathsf{c}} x^{1 + \mathsf{c}} \,(1-y)^{N-\textsf{n}_0-\textsf{n}_2 - \mathsf{c}} y^{1 + \mathsf{c}} \,(1-z)^{\textsf{n}_2 - 1 - \mathsf{c}} z^{\mathsf{c}} \,,
 \end{split}
\end{equation}
or, equivalently,
\begin{equation}\label{eq:Jpi3threethetas2}
 \begin{split}
 J_\varepsilon(x,y,z)\,&=\sum_{\textsf{n}_0=1}^{N-1}\sum_{\textsf{n}_2=1}^{N-\textsf{n}_0}\sum_{\mathsf{c}=0}^{\min\{\textsf{n}_0,\textsf{n}_1,\textsf{n}_2\}-1} \binom{\textsf{n}_0-1}{\mathsf{c}}\binom{\textsf{n}_1-1}{\mathsf{c}}\binom{\textsf{n}_2-1}{\mathsf{c}} \\
 &\qquad\times \big[1 -\varepsilon \textsf{n}_0 x - (1-\varepsilon)\textsf{n}_1 y + (\textsf{n}_2 - 1) z \big] \\
 &\qquad\times (1-x)^{\textsf{n}_0 - 1 - \mathsf{c}} x^{1 + \mathsf{c}} \,(1-y)^{\textsf{n}_1 - 1 - \mathsf{c}} y^{1 + \mathsf{c}} \,(1-z)^{\textsf{n}_2 - 1 - \mathsf{c}} z^{\mathsf{c}} \\
 & \textrm{with }\; \textsf{n}_1 = N+1-\textsf{n}_0-\textsf{n}_2\,.
 \end{split}
\end{equation}

This expression highlights the distinct roles of the three controls: $x$ and $y$ drive the population upwards along the ladder (appearing with power $1+\mathsf{c}$), while $z$ drives the feedback loop (appearing with power $\mathsf{c}$). The reward structure similarly separates the energy costs: penalties proportional to occupancy $\textsf{n}_0$ and $\textsf{n}_1$ weighted by the partial steps $\varepsilon$ and $1-\varepsilon$, versus a gain proportional to the occupancy of the target state $\textsf{n}_2$.

Observe that
\begin{equation}
    J_\varepsilon(x, y, z) \,=\, J_{1-\varepsilon}(y, x, z) \,.
\end{equation}
 Indeed, the term $(1-\varepsilon) \textsf{n}_0 y + \varepsilon \textsf{n}_1 x$ in \eqref{eq:Jpi3threethetas2} is invariant under simultaneous replacement $\varepsilon \mapsto 1-\varepsilon$ and exchange of controls $x \leftrightarrow y$, and the summation is symmetric over all valid trajectories $(\textsf{n}_0, \textsf{n}_1)$. As a consequence, if $(x^*, y^*, z^*)$ is the maximiser for $\varepsilon$, then $(y^*, x^*, z^*)$ is the maximiser for $1-\varepsilon$.

\begin{figure}[t!]
 \centering
 \includegraphics[height=4.8cm]{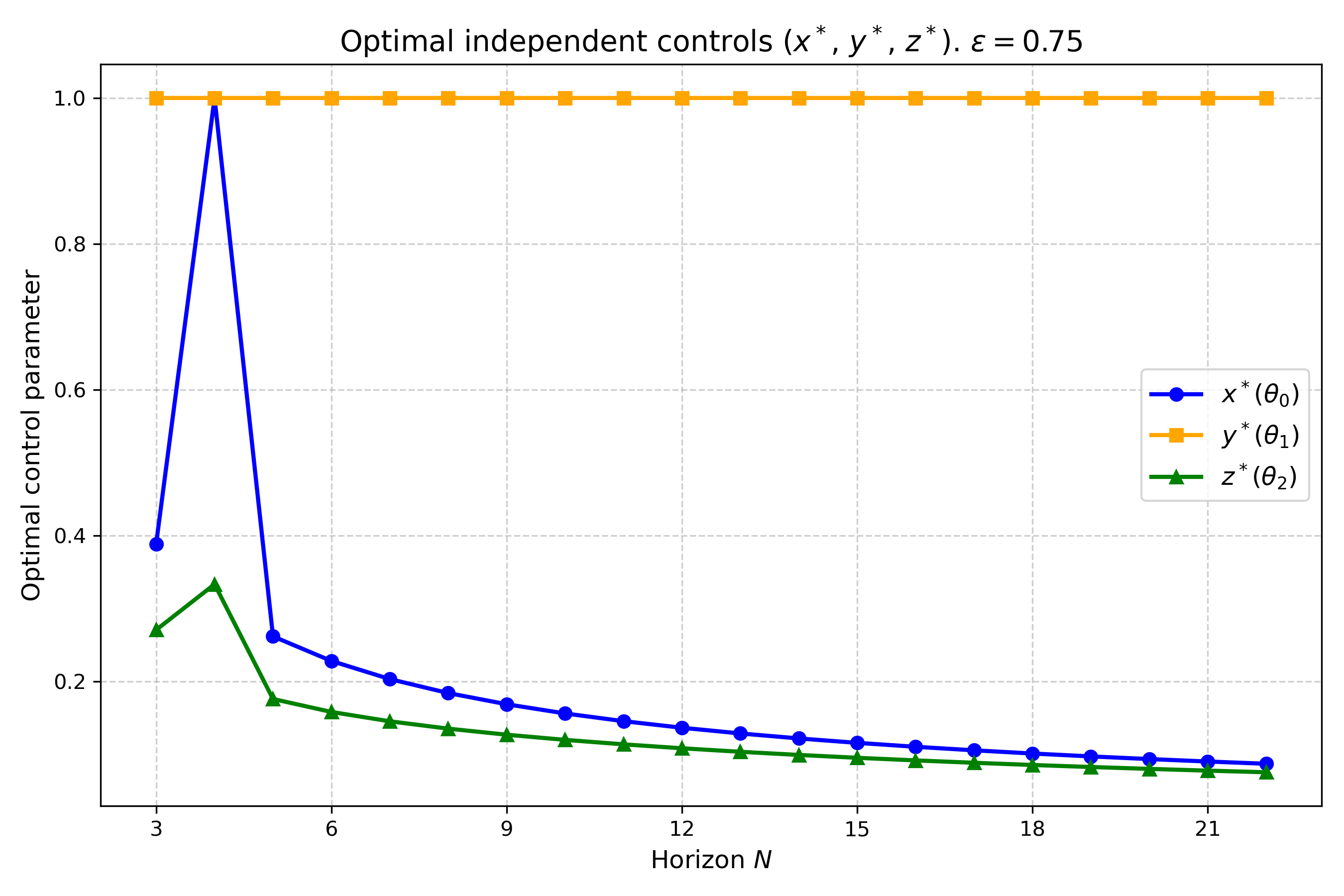}\;\;
 \includegraphics[height=4.8cm]{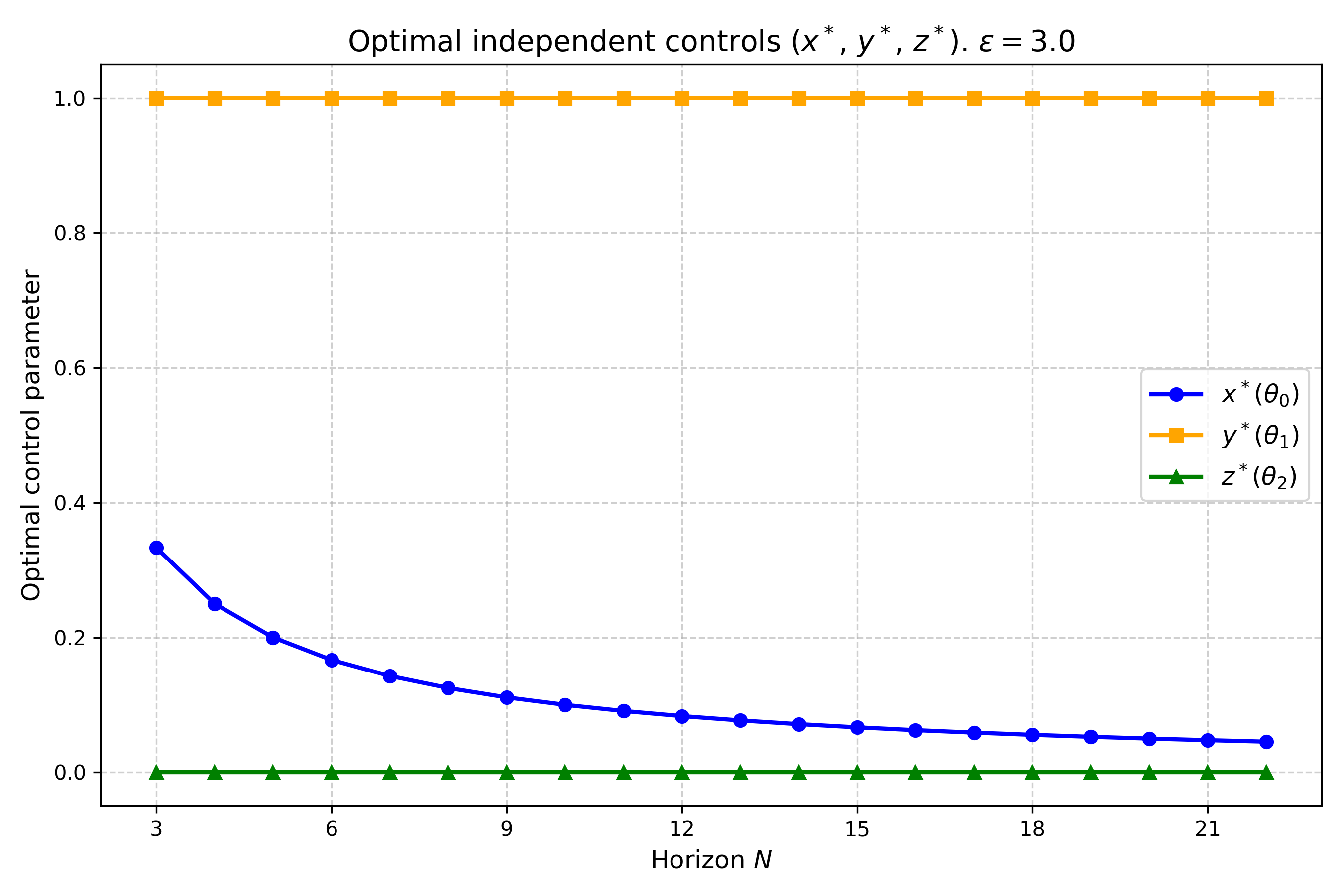}
 \caption{Optimal policy $(x^*,y^*,z^*)$ for the expected return
 \eqref{eq:Jpi3threethetas2} as a function of the horizon $N$, for
 $\varepsilon=0.75$ (left) and $\varepsilon=3.0$ (right).
 The saturations $y^*\to 1$ and $z^*\to 0$ are interpreted in the text.}
 \label{fig:Jpi-3level3VAR}
\end{figure}

Figure~\ref{fig:Jpi-3level3VAR} shows the results of the numerical
maximisation of $J_\varepsilon(x,y,z)$ across different trajectory lengths,
and reveals a clear migration of the optimal controls as $\varepsilon$ crosses
the threshold value $1$.

In the regime $\varepsilon\in(\frac{1}{2},1)$ (left panel), the reward
bracket in \eqref{eq:Jpi3threethetas2} penalises $y$ via
$-(1-\varepsilon)\textsf{n}_1 y$, so it is the probability weight
$(1-y)^{\textsf{n}_1-1-\mathsf{c}}y^{1+\mathsf{c}}$ that drives
$y^*\to 1$: as $y\to 1$ the weight collapses onto trajectories with
$\textsf{n}_1=1+\mathsf{c}$, i.e.\ minimal dwell at $|1\rangle$,
which corresponds to $\theta_1\to\frac{\pi}{2}$ and
$\pi_1|1\rangle\approx|2\rangle$. The cost asymmetry $\varepsilon>1-\varepsilon$
then keeps $x^*$ small: aggressive driving at $|0\rangle$ is penalised more
heavily than at $|1\rangle$.

In the regime $\varepsilon>1$ (right panel), the sign flip of
$-(1-\varepsilon)$ turns the $|1\rangle$ dwell into a reward source,
so reward and probability now both drive $y^*\to 1$. For $z^*\to 0$
the situation is reversed: the reward bracket favours $z>0$ via
$+(\textsf{n}_2-1)z$, but the factor $z^{\mathsf{c}}$ suppresses all
$\mathsf{c}\geqslant 1$ feedback trajectories as $z\to 0$, i.e.\
$\theta_2\to 0$ and $P_{20}=\sin^2\theta_2\to 0$. Here probability
dominates over reward: the gain from eliminating the feedback arc
outweighs the foregone $(\textsf{n}_2-1)z$ term, since each feedback
cycle would re-expose the system to the now-dominant $\varepsilon$-penalty
at $|0\rangle$.

%

 \section{4-level / two-qubit system}\label{sec:4levels}

 Last, we consider the two-qubit setting, a quantum system whose states are described by the Hilbert space
 \begin{equation}
  \cH\,=\,\mathbb{C}^2\otimes\mathbb{C}^2\,\cong\,\mathbb{C}^4\,.
 \end{equation}
 With the usual computational basis $\{|+\rangle,|-\rangle\}$ for the single qubit space $\mathbb{C}^2$, that is,
  \begin{equation}
 |1\rangle\,\equiv\,|+\rangle\,\equiv\,\begin{pmatrix} 1 \\ 0 \end{pmatrix}\,,\qquad |0\rangle\,\equiv\,|-\rangle \,\equiv\,\begin{pmatrix} 0 \\1 \end{pmatrix}\,,
 \end{equation}
 and with the usual representation
 \begin{equation}
  \begin{pmatrix} a_1 \\ a_2 \end{pmatrix} \otimes \begin{pmatrix} b_1 \\ b_2 \end{pmatrix}\,\equiv\,\begin{pmatrix} a_1 b_1 \\ a_1 b_2 \\ a_2 b_1 \\ a_2 b_2 \end{pmatrix}
 \end{equation}
 for $\mathbb{C}^2\otimes\mathbb{C}^2\cong\mathbb{C}^4$, the canonical computational basis for two qubits is $\{ |11\rangle,|10\rangle,|01\rangle,|00\rangle\}$ with
\begin{equation}\label{eq:compBasis2}
\begin{split}
|11\rangle \,\equiv\,  \begin{pmatrix} 1 \\ 0 \\ 0 \\ 0 \end{pmatrix}\!, \qquad |10\rangle \,\equiv\, \begin{pmatrix} 0 \\ 1 \\ 0 \\ 0 \end{pmatrix}\!, \qquad |01\rangle \,\equiv\, \begin{pmatrix} 0 \\ 0 \\ 1 \\ 0 \end{pmatrix}\!, \qquad |00\rangle \,\equiv\,~\begin{pmatrix} 0 \\ 0 \\ 0 \\ 1 \end{pmatrix}\!.
\end{split}
\end{equation}

 If the whole basis above is taken as family $\mathcal{S}$ of reference states, then from the perspective of the QRL models of the previous Sections one is effectively considering a 4-level system.

 In fact, further insight can be indeed gained by enriching the qutrits toy model by an additional intermediate-energy state. With respect to the notation therein, the setting now becomes
 \begin{equation}
  d=\mathrm{dim}\cH=4\,,\qquad\cH=\mathbb{C}^4\,,
 \end{equation}
 \begin{equation}\label{eq:states4}
\mathcal{S}\,=\,\big\{|0\rangle,|1\rangle,|1'\rangle,|2\rangle\big\}\;\; \textrm{ with }\;\; |0\rangle=\begin{pmatrix} 1 \\ 0 \\ 0 \\ 0 \end{pmatrix}\,,\;\;
   |1\rangle=\begin{pmatrix} 0 \\ 1 \\ 0 \\ 0 \end{pmatrix}\,,\;\;
   |1'\rangle=\begin{pmatrix} 0 \\ 0 \\ 1 \\ 0 \end{pmatrix}\,,\;\;
   |2\rangle=\begin{pmatrix} 0 \\ 0 \\ 0 \\ 1 \end{pmatrix}\,,
 \end{equation}
 \begin{equation}\label{eq:Hamilt4}
  \begin{split}
   & H\,=\,E_0|0\rangle\langle 0|+E_1|1\rangle\langle 1|+E_1'|1'\rangle\langle 1'|+E_2|2\rangle\langle 2| \\
   & \textrm{with }\quad E_0=0\,,\quad E_1=\varepsilon\in(0,1)\,,\quad E_1'=\varepsilon'\in(0,1)\,,\quad E_2=1\,.
  \end{split}
 \end{equation}

\medskip

 \textbf{Policies.} Each policy $\pi\equiv(\pi_0,\pi_1,\pi_{1'},\pi_2)\equiv(\theta,\theta')$, with $\theta,\theta'\in[0,2\pi)$ and unitaries
 \begin{equation}\label{eq:4policiespi0}
\begin{split}
  \pi_0 &:= \begin{pmatrix}
\cos^2\theta & \sin\theta & \sin\theta\cos\theta & 0 \\
-\sin\theta\cos\theta & \cos\theta & -\sin^2\theta & 0 \\
-\sin\theta & 0 & \cos\theta & 0 \\
0 & 0 & 0 & 1
\end{pmatrix}, \\
\pi_1 &:= \begin{pmatrix}
 1 & 0 & 0 & 0 \\
 0 & \cos\theta & 0 & -\sin\theta \\
 0 & 0 & 1 & 0 \\
 0 & \sin\theta & 0 & \cos\theta
 \end{pmatrix}, \\
\pi_{1'} &:= \begin{pmatrix}
1 & 0 & 0 & 0 \\
0 & 1 & 0 & 0 \\
0 & 0 & \cos\theta' & -\sin\theta' \\
0 & 0 & \sin\theta' & \cos\theta'
\end{pmatrix}, \\
\pi_2 &:= \begin{pmatrix}
 \cos\theta & 0 & 0 & -\sin\theta \\
 0 & 1 & 0 & 0 \\
 0 & 0 & 1 & 0 \\
 \sin\theta & 0 & 0 & \cos\theta
 \end{pmatrix}.
\end{split}
\end{equation}
 rotates the states as
 \begin{equation}
  |j\rangle\,\mapsto\,\pi_j|j\rangle\,,\qquad j\in\{0,1,1',2\}\,,
 \end{equation}
 hence
 \begin{eqnarray}
  & & |0\rangle\,\mapsto\,  \cos^2\theta|0\rangle - \sin\theta\cos\theta|1\rangle - \sin\theta|1'\rangle\,, \label{eq:unitaries4-1} \\
  & & |1\rangle\,\mapsto\, \cos\theta|1\rangle + \sin\theta|2\rangle\,, \label{eq:unitaries4-2} \\
  & & |1'\rangle\,\mapsto\,\cos\theta'|1'\rangle + \sin\theta'|2\rangle\,, \label{eq:unitaries4-3} \\
  & & |2\rangle\,\mapsto\, -\sin\theta|0\rangle + \cos\theta|2\rangle\,. \label{eq:unitaries4-4}
\end{eqnarray}

\medskip

\textbf{Trajectories.} Admitted trajectories are those realising the (`lowest-to-highest-energy-level') overall transition $|0\rangle\to |2\rangle$ through $N+1$ states ($N$ steps). A trajectory $\tau$ of this type has $\textsf{n}_j$ vectors of type $|j\rangle$, with $\textsf{n}_0+\textsf{n}_1+\textsf{n}_{1'}+\textsf{n}_2=N+1$, where
\begin{equation}\label{eq:njconstraints4}
 \begin{split}
   \textsf{n}_0\,&\in\,\{1,\dots,N-1\}\,, \\
   \textsf{n}_2\,&\in\,\{1,\dots,N-\textsf{n}_0\}\,, \\
   \textsf{n}_1\,&\in\,\{0,\dots, N+1-\textsf{n}_0-\textsf{n}_2\}\,, \\
   \textsf{n}_{1'}\,&=\,N+1-\textsf{n}_0-\textsf{n}_1-\textsf{n}_2\,.
 \end{split}
\end{equation}
(Observe that the constraint $\textsf{n}_0+\textsf{n}_2\leqslant N$ ensures that $\textsf{n}_1+\textsf{n}_{1'}\geqslant 1$, guaranteeing at least one visit to the intermediate layer.)

In terms of the parameters
\begin{equation}
 \mathsf{a}_{ij}\,:=\,\#[ij] \qquad \big(\textrm{number of transitions $|i\rangle\xrightarrow{\pi_i}\pi_i|i\rangle\xrightarrow{(\textrm{collapse})}|j\rangle$}\big)\,
\end{equation}
the transition counts must satisfy the following constraints. First, the total number of transitions equals the number of steps:
\begin{equation}\label{eq:4sumaij}
\sum_{i,j} \mathsf{a}_{ij} = N\,.
\end{equation}
Second, each state appearance (except the final state) produces exactly one outgoing transition:
\begin{equation}\label{eq:4outgoingaij}
\begin{split}
 \text{(from $|0\rangle$)} \quad &\mathsf{a}_{00} + \mathsf{a}_{01} + \mathsf{a}_{01'} + \mathsf{a}_{02} = \textsf{n}_0\,, \\
 \text{(from $|1\rangle$)} \quad &\mathsf{a}_{10} + \mathsf{a}_{11} + \mathsf{a}_{11'} + \mathsf{a}_{12} = \textsf{n}_1\,, \\
 \text{(from $|1'\rangle$)} \quad &\mathsf{a}_{1'0} + \mathsf{a}_{1'1} + \mathsf{a}_{1'1'} + \mathsf{a}_{1'2} = \textsf{n}_{1'}\,, \\
 \text{(from $|2\rangle$)} \quad &\mathsf{a}_{20} + \mathsf{a}_{21} + \mathsf{a}_{21'} + \mathsf{a}_{22} = \textsf{n}_2 - 1\,.
\end{split}
\end{equation}
Third, each state appearance (except the initial state) requires exactly one incoming transition:
\begin{equation}\label{eq:4incomingconstr}
\begin{split}
 \text{(to $|0\rangle$)} \quad &\mathsf{a}_{00} + \mathsf{a}_{10} + \mathsf{a}_{1'0} + \mathsf{a}_{20} = \textsf{n}_0 - 1\,, \\
 \text{(to $|1\rangle$)} \quad &\mathsf{a}_{01} + \mathsf{a}_{11} + \mathsf{a}_{1'1} + \mathsf{a}_{21} = \textsf{n}_1\,, \\
 \text{(to $|1'\rangle$)} \quad &\mathsf{a}_{01'} + \mathsf{a}_{11'} + \mathsf{a}_{1'1'} + \mathsf{a}_{21'} = \textsf{n}_{1'}\,, \\
 \text{(to $|2\rangle$)} \quad &\mathsf{a}_{02} + \mathsf{a}_{12} + \mathsf{a}_{1'2} + \mathsf{a}_{22} = \textsf{n}_2\,.
\end{split}
\end{equation}
These constraints automatically ensure the correct endpoint behaviour: starting at $|0\rangle$ (one extra outgoing transition from $|0\rangle$ compared to incoming) and ending at $|2\rangle$ (one extra incoming transition to $|2\rangle$ compared to outgoing).

In particular, comparing the outgoing and incoming constraints for states $|1\rangle$ and $|1'\rangle$ (which are neither initial nor final), one obtains the flow conservation relations
\begin{equation}\label{eq:4flowconserv}
\begin{split}
\text{(through $|1\rangle$)} \quad &\mathsf{a}_{10} + \mathsf{a}_{11'} + \mathsf{a}_{12} \,=\,\mathsf{a}_{01} + \mathsf{a}_{1'1} + \mathsf{a}_{21}\,, \\
\text{(through $|1'\rangle$)} \quad &\mathsf{a}_{1'0} + \mathsf{a}_{1'1} + \mathsf{a}_{1'2} \,=\, \mathsf{a}_{01'} + \mathsf{a}_{11'} + \mathsf{a}_{21'}\,.
\end{split}
\end{equation}
These express that for intermediate states, the total outgoing flow (excluding self-loops, which contribute equally to both sides) equals the total incoming flow. Note that the diagonal terms $\mathsf{a}_{11}$ and $\mathsf{a}_{1'1'}$ cancel when subtracting incoming from outgoing constraints, reflecting that self-loop transitions contribute equally to both outgoing and incoming counts for a given state.

The conservation relations for the end points, analogous to \eqref{eq:4flowconserv}, read
\begin{equation}\label{eq:4flowconserv_endpoints}
\begin{split}
\text{(through $|0\rangle$)} \quad
  &\mathsf{a}_{01} + \mathsf{a}_{01'} + \mathsf{a}_{02}
   - \mathsf{a}_{10} - \mathsf{a}_{1'0} - \mathsf{a}_{20} \,=\, +1\,, \\
\text{(through $|2\rangle$)} \quad
  &\mathsf{a}_{20} + \mathsf{a}_{21} + \mathsf{a}_{21'}
   - \mathsf{a}_{02} - \mathsf{a}_{12} - \mathsf{a}_{1'2} \,=\, -1\,.
\end{split}
\end{equation}
These are on the same footing as \eqref{eq:netflow3} from the previous three-level model. Together, \eqref{eq:4flowconserv} and \eqref{eq:4flowconserv_endpoints} account for all
four nodes of the directed graph $\mathcal{G}_4$ (Figure \ref{fig:graph-fourlevel}) and express the global Kirchhoff balance: state $|0\rangle$ has net outflow $+1$ (it is the trajectory's source), state $|2\rangle$ has net inflow $-1$ (it is the trajectory's sink), and the two intermediate states $|1\rangle$ and
$|1'\rangle$ have zero net flow (they are pure transit nodes).
The algebraic sum of the four net flows vanishes, $(+1) + 0 + 0 + (-1) = 0$,
consistently with the global conservation of probability flow through the system.

\medskip

\textbf{Allowed transitions.} The factors
\begin{equation}
 |\langle j | \pi_i | i \rangle |^2 \,=:\,P_{ij}\,,\qquad i,j\in\{0,1,1',2\}
\end{equation}
 constituting the general formula \eqref{eq:probtraj} for each transition $|i\rangle\xrightarrow{\pi_i}\pi_i|i\rangle\xrightarrow{(\textrm{collapse})}|j\rangle$ are computed by means of \eqref{eq:states4} and \eqref{eq:unitaries4-1}-\eqref{eq:unitaries4-4}, yielding the transition probability matrix
 \begin{equation}\label{eq:4transitionprobmatr}
\mathbf{P} \,:=\,(P_{ij})\,=\,\begin{pmatrix}
\cos^4\theta & \sin^2\theta\cos^2\theta & \sin^2\theta & 0 \\
0 & \cos^2\theta & 0 & \sin^2\theta \\
0 & 0 & \cos^2\theta' & \sin^2\theta' \\
\sin^2\theta & 0 & 0 & \cos^2\theta
\end{pmatrix}\,.
\end{equation}
 Transitions with null probability (zero entries in $\mathbf{P}$) are factored out of the product \eqref{eq:probtraj} expressing the trajectory's probability, by restricting their frequencies to zero. This results in the additional constraints
 \begin{equation}\label{eq:4additionalConstr}
  \textsf{a}_{02}\,=\,\textsf{a}_{10}\,=\,\textsf{a}_{11'}\,=\,\textsf{a}_{1'0}\,=\,\textsf{a}_{1'1}\,=\,\textsf{a}_{21}\,=\,\textsf{a}_{21'}\,=\,0\,.
 \end{equation}

 \medskip

 \textbf{Independent degrees of freedom.} By the same Kirchhoff argument as in Section~\ref{sec:Qtrits}, the net-flow equations
analogous to \eqref{eq:netflow3} for the directed graph $\mathcal{G}_4$ with node set
$\{|0\rangle,|1\rangle,|1'\rangle,|2\rangle\}$ and full off-diagonal edge set
(Figure~\ref{fig:graph-fourlevel}) yield exactly $4-1=3$ independent conditions
\cite[Theorem~2.3]{Bapat2014}. Together with the $4$ mutually independent outgoing
constraints \eqref{eq:4outgoingaij} and the $7$ forbidden-transition conditions
\eqref{eq:4additionalConstr}, one obtains $4+3+7=14$ independent conditions on the
$16$ transition counts $\mathsf{a}_{ij}$, leaving exactly $16-14=2$ degrees of
freedom. A convenient choice of the two free parameters is
\begin{equation}\label{eq:4parametersc}
 \begin{split}
  \mathsf{c}_{01} &:= \mathsf{a}_{01}\,, \quad \mathsf{c}_{01'} := \mathsf{a}_{01'}\,.
 \end{split}
\end{equation}
The remaining variables are uniquely determined by the state counts and these parameters:
\begin{equation}\label{eq:4Param4SIMPL}
\begin{aligned}
\mathsf{a}_{00} \,&=\, \textsf{n}_0 - \mathsf{c}_{01} - \mathsf{c}_{01'} \,, & \mathsf{a}_{01} \,&=\, \mathsf{c}_{01} \,, \\
\mathsf{a}_{01'} \,&=\, \mathsf{c}_{01'} \,, & \mathsf{a}_{02} \,&=\, 0 \,, \\
\mathsf{a}_{10} \,&=\, 0 \,, & \mathsf{a}_{11} \,&=\, \textsf{n}_1 - \mathsf{c}_{01} \,, \\
\mathsf{a}_{11'} \,&=\, 0 \,, & \mathsf{a}_{12} \,&=\, \mathsf{c}_{01} \,, \\
\mathsf{a}_{1'0} \,&=\, 0 \,, & \mathsf{a}_{1'1} \,&=\, 0 \,, \\
\mathsf{a}_{1'1'} \,&=\, \textsf{n}_{1'} - \mathsf{c}_{01'} \,, & \mathsf{a}_{1'2} \,&=\, \mathsf{c}_{01'} \,, \\
\mathsf{a}_{20} \,&=\, \mathsf{c}_{01} + \mathsf{c}_{01'} - 1 \,, & \mathsf{a}_{21} \,&=\, 0 \,, \\
\mathsf{a}_{21'} \,&=\, 0 \,, & \mathsf{a}_{22} \,&=\, \textsf{n}_2 - \mathsf{c}_{01} - \mathsf{c}_{01'} \,.
\end{aligned}
\end{equation}
For given $N$, the admissible ranges for $\textsf{n}_0$, $\textsf{n}_1$, $\textsf{n}_{1'}$, $\textsf{n}_2$, $\mathsf{c}_{01}$, and $\mathsf{c}_{01'}$ are determined by the constraint that all integers $\mathsf{a}_{ij}$ in \eqref{eq:4Param4SIMPL} are non-negative, compatibly with \eqref{eq:njconstraints4}.
This yields the conditions
\begin{equation}\label{eq:njconstraints4UPD}
\begin{split}
\textsf{n}_0\,&\in\, \{1, \ldots, N-1\}\,, \\
\textsf{n}_2\,&\in\, \{1, \ldots, N-\textsf{n}_0\}\,, \\
\textsf{n}_1\,&\in\, \{0, \ldots, N+1-\textsf{n}_0-\textsf{n}_2\}\,, \\
\textsf{n}_{1'}\,&=\, N+1-\textsf{n}_0-\textsf{n}_1-\textsf{n}_2\,, \\
\mathsf{c}_{01}\,&\in\, \{0, \ldots, \min\{\textsf{n}_1, \textsf{n}_0, \textsf{n}_2\}\}\,, \\
\mathsf{c}_{01'}\,&\in\, \{\max\{0, 1-\mathsf{c}_{01}\}, \ldots, \min\{\textsf{n}_{1'}, \textsf{n}_0-\mathsf{c}_{01}, \textsf{n}_2-\mathsf{c}_{01}\}\}\,.
\end{split}
\end{equation}

\begin{figure}[t!]
\centering
\begin{tikzpicture}[>=stealth, semithick]
  \node[draw, circle, inner sep=4pt] (n0)  at ( 0.0,  2.4) {$|0\rangle$};
  \node[draw, circle, inner sep=4pt] (n1)  at (-2.4,  0.0) {$|1\rangle$};
  \node[draw, circle, inner sep=4pt] (n1p) at ( 2.4,  0.0) {$|1'\rangle$};
  \node[draw, circle, inner sep=4pt] (n2)  at ( 0.0, -2.4) {$|2\rangle$};
  \draw[->] (n0)  to[bend right=15] (n1);
  \draw[->] (n1)  to[bend right=15] (n0);
  \draw[->] (n0)  to[bend right=15] (n1p);
  \draw[->] (n1p) to[bend right=15] (n0);
  \draw[->] (n0)  to[bend right=15] (n2);
  \draw[->] (n2)  to[bend right=15] (n0);
  \draw[->] (n1)  to[bend right=15] (n1p);
  \draw[->] (n1p) to[bend right=15] (n1);
  \draw[->] (n1)  to[bend right=15] (n2);
  \draw[->] (n2)  to[bend right=15] (n1);
  \draw[->] (n1p) to[bend right=15] (n2);
  \draw[->] (n2)  to[bend right=15] (n1p);
\end{tikzpicture}
\qquad\qquad
\begin{tikzpicture}[>=stealth, semithick]
  \node[draw, circle, minimum size=0.85cm] (n0)  at ( 0.0,  2.4) {$|0\rangle$};
  \node[draw, circle, minimum size=0.85cm] (n1)  at (-2.4,  0.0) {$|1\rangle$};
  \node[draw, circle, minimum size=0.85cm] (n1p) at ( 2.4,  0.0) {$|1^{\prime}\rangle$};
  \node[draw, circle, minimum size=0.85cm] (n2)  at ( 0.0, -2.4) {$|2\rangle$};
  \draw[->] (n0) to[bend right=15]
    node[left,  midway] {$\mathsf{c}_{01}\;$}  (n1);
  \draw[->] (n0) to[bend right=15]
    node[right, midway] {$\;\mathsf{c}_{01'}$} (n1p);
  \draw[->] (n1) to[bend right=15]
    node[left,  midway] {$\mathsf{c}_{01}\;$}  (n2);
  \draw[->] (n1p) to[bend right=15]
    node[right, midway] {$\;\mathsf{c}_{01'}$} (n2);
  \draw[->] (n2) to
    node[right=1mm, midway, align=center]
      {$\mathsf{c}_{01}{+}\mathsf{c}_{01'}$\\$-1$} (n0);
\end{tikzpicture}
\caption{Left: the directed graph $\mathcal{G}_4$ on node set
$\{|0\rangle,|1\rangle,|1'\rangle,|2\rangle\}$ with all twelve off-diagonal edges
($i\neq j$), whose incidence matrix has rank $4-1=3$ \cite[Theorem~2.3]{Bapat2014}.
Right: the reduced transition graph $\widetilde{\mathcal{G}}_4$, obtained from
$\mathcal{G}_4$ by retaining only the five edges $|i\rangle\to|j\rangle$ with
nonzero probability after imposing the forbidden-transition constraints
\eqref{eq:4additionalConstr}; edge labels indicate the transition multiplicities
$\mathsf{a}_{ij}$ from \eqref{eq:4Param4SIMPL}, for
$\mathsf{c}_{01}+\mathsf{c}_{01'}\geqslant 1$.}
\label{fig:graph-fourlevel}
\end{figure}

 \medskip

 \textbf{Trajectory probability.} Specialising \eqref{eq:probtraj} to the transitions allowed for the present model, hence taking the product of each non-zero entry $P_{ij}$ of \eqref{eq:4transitionprobmatr}, exponentiated to the corresponding transition count $\mathsf{a}_{ij}$ given by \eqref{eq:4Param4SIMPL}, yields the probability $P_\pi(\tau)$ of a lowest-to-highest-energy-level, length-$N$ trajectory $\tau$ with states counts $\textsf{n}_0,\textsf{n}_1,\textsf{n}_{1'},\textsf{n}_2$ and transition counts $\mathsf{a}_{ij}$, determined by a choice of $N$, $\textsf{n}_0$, $\textsf{n}_1$, $\textsf{n}_{1'}$, $\textsf{n}_2$, $\mathsf{c}_{01}$, $\mathsf{c}_{01'}$ compatible with \eqref{eq:njconstraints4UPD}:
\begin{equation}\label{eq:ProbTraj4}
\begin{split}
P_\pi(\tau) \,&=\, (\cos^4\theta)^{\mathsf{a}_{00}}\;(\sin^2\theta\cos^2\theta)^{\mathsf{a}_{01}}\; (\sin^2\theta)^{\mathsf{a}_{01'}}\;(\cos^2\theta)^{\mathsf{a}_{11}}\;(\sin^2\theta)^{\mathsf{a}_{12}} \,\times \\
&\qquad\times (\cos^2\theta')^{\mathsf{a}_{1'1'}}\;(\sin^2\theta')^{\mathsf{a}_{1'2}}\;(\sin^2\theta)^{\mathsf{a}_{20}}\;(\cos^2\theta)^{\mathsf{a}_{22}}\\
&=\, (\cos^4\theta)^{\textsf{n}_0 - \mathsf{c}_{01} - \mathsf{c}_{01'}}\;(\sin^2\theta\cos^2\theta)^{\mathsf{c}_{01}}\; (\sin^2\theta)^{\mathsf{c}_{01'}}\;(\cos^2\theta)^{\textsf{n}_1 - \mathsf{c}_{01}}\;(\sin^2\theta)^{\mathsf{c}_{01}} \,\times \\
&\qquad\times (\cos^2\theta')^{\textsf{n}_{1'} - \mathsf{c}_{01'}}\;(\sin^2\theta')^{\mathsf{c}_{01'}}\;(\sin^2\theta)^{\mathsf{c}_{01} + \mathsf{c}_{01'} - 1}\;(\cos^2\theta)^{\textsf{n}_2 - \mathsf{c}_{01} - \mathsf{c}_{01'}}\\
&=\, (\cos^2\theta)^{2\textsf{n}_0 + \textsf{n}_1 + \textsf{n}_2 - 3\mathsf{c}_{01} - 3\mathsf{c}_{01'}} \, (\sin^2\theta)^{3\mathsf{c}_{01} + 2\mathsf{c}_{01'} - 1} \, (\cos^2\theta')^{\textsf{n}_{1'} - \mathsf{c}_{01'}} \, (\sin^2\theta')^{\mathsf{c}_{01'}} \\
&=:\,P(N,\textsf{n}_0,\textsf{n}_1,\textsf{n}_{1'},\textsf{n}_2,\mathsf{c}_{01},\mathsf{c}_{01'},\theta,\theta')\,.
\end{split}
\end{equation}

\medskip

\textbf{Trajectory multiplicity.} Let $\mathcal{N}(N,\textsf{n}_0,\textsf{n}_1,\textsf{n}_{1'},\textsf{n}_2,\mathsf{c}_{01},\mathsf{c}_{01'})$ be the number of distinct lowest-to-highest-energy-level, length-$N$ trajectories with state counts $\textsf{n}_0, \textsf{n}_1, \textsf{n}_{1'}, \textsf{n}_2$ and transition counts $\mathsf{a}_{ij}$ given by \eqref{eq:4Param4SIMPL}.
The only allowed non-self-loop transitions are $|0\rangle\to|1\rangle$, $|0\rangle\to|1'\rangle$, $|1\rangle\to|2\rangle$, $|1'\rangle\to|2\rangle$, and $|2\rangle\to|0\rangle$. Setting
\begin{equation*}
 \mathsf{K}\,:=\,\mathsf{c}_{01}+\mathsf{c}_{01'}\,,
\end{equation*}
each admissible trajectory consists of $\mathsf{K}$ distinct cycles of the form $|0\rangle\to|\bullet\rangle\to|2\rangle\to|0\rangle$, with the last cycle terminating at $|2\rangle$ instead of returning to $|0\rangle$ (corresponding to $\mathsf{a}_{20}=\mathsf{K}-1$); in each cycle the intermediate state $|\bullet\rangle$ is either $|1\rangle$ or $|1'\rangle$, with exactly $\mathsf{c}_{01}$ cycles passing through $|1\rangle$ and $\mathsf{c}_{01'}$ cycles passing through $|1'\rangle$. The forbidden transitions $\mathsf{a}_{11'}=\mathsf{a}_{1'1}=0$ ensure that no cycle visits both intermediate states. Counting the admissible trajectories therefore amounts to the following four independent combinatorial choices:
\begin{itemize}
 \item distributing the $\textsf{n}_0$ visits to $|0\rangle$ into the $\mathsf{K}$ runs of consecutive $|0\rangle$-states (one run per cycle), each run non-empty: $\binom{\textsf{n}_0-1}{\mathsf{K}-1}$ ways;
 \item distributing the $\textsf{n}_2$ visits to $|2\rangle$ into the $\mathsf{K}$ runs of consecutive $|2\rangle$-states, each run non-empty: $\binom{\textsf{n}_2-1}{\mathsf{K}-1}$ ways;
 \item choosing which $\mathsf{c}_{01}$ of the $\mathsf{K}$ cycles pass through $|1\rangle$ (the remaining $\mathsf{c}_{01'}$ pass through $|1'\rangle$): $\binom{\mathsf{K}}{\mathsf{c}_{01}}$ ways;
 \item distributing independently the $\textsf{n}_1$ visits to $|1\rangle$ into the $\mathsf{c}_{01}$ runs that pass through $|1\rangle$, and the $\textsf{n}_{1'}$ visits to $|1'\rangle$ into the $\mathsf{c}_{01'}$ runs that pass through $|1'\rangle$: $\binom{\textsf{n}_1-1}{\mathsf{c}_{01}-1}\binom{\textsf{n}_{1'}-1}{\mathsf{c}_{01'}-1}$ ways.
\end{itemize}
The four choices are mutually independent, hence their multiplication yields the exact closed-form expression
\begin{equation}\label{eq:multiplicity4-closedform}
 \mathcal{N}(N,\textsf{n}_0,\textsf{n}_1,\textsf{n}_{1'},\textsf{n}_2,\mathsf{c}_{01},\mathsf{c}_{01'}) \,=\,
 \binom{\textsf{n}_0-1}{\mathsf{K}-1}
 \binom{\textsf{n}_2-1}{\mathsf{K}-1}
 \binom{\mathsf{K}}{\mathsf{c}_{01}}
 \binom{\textsf{n}_1-1}{\mathsf{c}_{01}-1}
 \binom{\textsf{n}_{1'}-1}{\mathsf{c}_{01'}-1}
\end{equation}
with the convention $\binom{-1}{-1}:=1$ (empty selection from the empty set), which correctly handles the degenerate cases $\mathsf{c}_{01}=0\wedge\textsf{n}_1=0$ (no cycle uses the $|1\rangle$ branch) and the symmetric one $\mathsf{c}_{01'}=0\wedge\textsf{n}_{1'}=0$. Notably, \eqref{eq:multiplicity4-closedform} is a product of \emph{five} binomial factors, rather than the four one might naively expect from the four state-occupation counts: the extra factor $\binom{\mathsf{K}}{\mathsf{c}_{01}}$ encodes the ordering of the intermediate-branch choices across the $\mathsf{K}$ cycles, a degree of freedom absent in the three-level model of Section~\ref{sec:Qtrits}. As a consistency check, the closed-form expression \eqref{eq:multiplicity4-closedform} reduces to the qutrit multiplicity \eqref{eq:3trajmult} of Section~\ref{sec:Qtrits} when one of the two intermediate branches is suppressed. Setting $\mathsf{c}_{01'}=0$ and $\textsf{n}_{1'}=0$ gives $\mathsf{K}=\mathsf{c}_{01}$, $\binom{\mathsf{K}}{\mathsf{c}_{01}}=1$, $\binom{-1}{-1}=1$, and
\begin{equation*}
 \mathcal{N}(N,\textsf{n}_0,\textsf{n}_1,0,\textsf{n}_2,\mathsf{c}_{01},0) \,=\, \binom{\textsf{n}_0-1}{\mathsf{c}_{01}-1}\binom{\textsf{n}_1-1}{\mathsf{c}_{01}-1}\binom{\textsf{n}_2-1}{\mathsf{c}_{01}-1}\,,
\end{equation*}
in agreement with \eqref{eq:3trajmult} upon the identification $\mathsf{c}=\mathsf{c}_{01}-1$. The symmetric reduction $\mathsf{c}_{01}=0$, $\textsf{n}_1=0$ yields the same expression with $\textsf{n}_{1'}$ in place of $\textsf{n}_1$.

Observe that the combinatorial argument leading to \eqref{eq:multiplicity4-closedform} requires that whenever the intermediate state $|1\rangle$ is visited ($\textsf{n}_1\geqslant 1$) at least one cycle must pass through it, i.e. $\mathsf{c}_{01}\geqslant 1$; symmetrically, $\textsf{n}_{1'}\geqslant 1\Rightarrow\mathsf{c}_{01'}\geqslant 1$. The lower bounds in \eqref{eq:njconstraints4UPD} are accordingly refined to
\begin{equation}\label{eq:njconstraints4UPDcorr}
\begin{split}
\mathsf{c}_{01}\,&\in\,\big\{\mathbf{1}_{\{\textsf{n}_1\geqslant 1\}},\,\ldots,\,\min\{\textsf{n}_1,\textsf{n}_0,\textsf{n}_2\}\big\}\,,\\
\mathsf{c}_{01'}\,&\in\,\big\{\max\!\big\{\mathbf{1}_{\{\textsf{n}_{1'}\geqslant 1\}},\,1-\mathsf{c}_{01}\big\},\,\ldots,\,\min\{\textsf{n}_{1'},\textsf{n}_0-\mathsf{c}_{01},\textsf{n}_2-\mathsf{c}_{01}\}\big\}\,,
\end{split}
\end{equation}
with $\mathbf{1}_{\{\,\cdot\,\}}$ the indicator function (value $1$ if the bracketed condition holds, $0$ otherwise).

\medskip

Table~\ref{tab:multiplicities} reports the case $N=5$, computed directly from \eqref{eq:multiplicity4-closedform}.

\begin{table}[h]
\centering
\begin{tabular}{cccccc|c}
$\textsf{n}_0$ & $\textsf{n}_1$ & $\textsf{n}_{1'}$ & $\textsf{n}_2$ & $\mathsf{c}_{01}$ & $\mathsf{c}_{01'}$ & $\mathcal{N}(N,\textsf{n}_0,\textsf{n}_1,\textsf{n}_{1'},\textsf{n}_2,\mathsf{c}_{01},\mathsf{c}_{01'})$ \\
\hline
1 & 0 & 4 & 1 & 0 & 1 & 1 \\
1 & 4 & 0 & 1 & 1 & 0 & 1 \\
1 & 0 & 3 & 2 & 0 & 1 & 1 \\
1 & 3 & 0 & 2 & 1 & 0 & 1 \\
1 & 0 & 2 & 3 & 0 & 1 & 1 \\
1 & 2 & 0 & 3 & 1 & 0 & 1 \\
1 & 0 & 1 & 4 & 0 & 1 & 1 \\
1 & 1 & 0 & 4 & 1 & 0 & 1 \\
2 & 0 & 3 & 1 & 0 & 1 & 1 \\
2 & 3 & 0 & 1 & 1 & 0 & 1 \\
2 & 0 & 2 & 2 & 0 & 1 & 1 \\
2 & 0 & 2 & 2 & 0 & 2 & 1 \\
2 & 1 & 1 & 2 & 1 & 1 & 2 \\
2 & 2 & 0 & 2 & 1 & 0 & 1 \\
2 & 2 & 0 & 2 & 2 & 0 & 1 \\
2 & 0 & 1 & 3 & 0 & 1 & 1 \\
2 & 1 & 0 & 3 & 1 & 0 & 1 \\
3 & 0 & 2 & 1 & 0 & 1 & 1 \\
3 & 2 & 0 & 1 & 1 & 0 & 1 \\
3 & 0 & 1 & 2 & 0 & 1 & 1 \\
3 & 1 & 0 & 2 & 1 & 0 & 1 \\
4 & 0 & 1 & 1 & 0 & 1 & 1 \\
4 & 1 & 0 & 1 & 1 & 0 & 1 \\
\hline
\multicolumn{6}{c|}{Total configurations:} & 23 \\
\multicolumn{6}{c|}{Total trajectories:} & 24 \\
\end{tabular}
\medskip
\caption{Trajectory multiplicities for $N=5$, computed from the closed-form expression \eqref{eq:multiplicity4-closedform}. Among the 23 admissible configurations, only one, namely $(\textsf{n}_0,\textsf{n}_1,\textsf{n}_{1'},\textsf{n}_2,\mathsf{c}_{01},\mathsf{c}_{01'})=(2,1,1,2,1,1)$, has multiplicity greater than one: it corresponds to a `mixed' configuration in which the two cycles use different intermediate branches, one passing through $|1\rangle$ and the other through $|1'\rangle$, with the factor $\binom{\mathsf{K}}{\mathsf{c}_{01}}=\binom{2}{1}=2$ in \eqref{eq:multiplicity4-closedform} encoding the choice of which cycle uses which branch.}
\label{tab:multiplicities}
\end{table}

\medskip

 \textbf{Reward along a trajectory.} Each transition $|i\rangle\xrightarrow{\pi_i}\pi_i|i\rangle\xrightarrow{(\textrm{collapse})}|j\rangle$ along the trajectory is rewarded by the energy difference
 \begin{equation}
   E_j - \langle i|\pi_i^*H\pi_i| i\rangle\,=:\, R_{ij}\,,\qquad i,j\in\{0,1,1',2\}\,.
\end{equation}
 Their computation, using \eqref{eq:states4}-\eqref{eq:Hamilt4} and \eqref{eq:unitaries4-1}-\eqref{eq:unitaries4-4} yields the transition reward matrix
 \begin{equation}\label{eq:4rewmatr}
  \begin{split}
     \mathbf{R} \,:=&\,(R_{ij}) \\
     =&\,\tiny\begin{pmatrix}
-\varepsilon\sin^2\theta\cos^2\theta - \varepsilon'\sin^2\theta & \varepsilon - \varepsilon\sin^2\theta\cos^2\theta - \varepsilon'\sin^2\theta & \varepsilon' - \varepsilon\sin^2\theta\cos^2\theta - \varepsilon'\sin^2\theta & 1 - \varepsilon\sin^2\theta\cos^2\theta - \varepsilon'\sin^2\theta \\
-\varepsilon\cos^2\theta - \sin^2\theta & -(1-\varepsilon)\sin^2\theta & \varepsilon' - \varepsilon\cos^2\theta - \sin^2\theta & (1-\varepsilon)\cos^2\theta \\
-\varepsilon'\cos^2\theta' - \sin^2\theta' & \varepsilon - \varepsilon'\cos^2\theta' - \sin^2\theta' & -(1-\varepsilon')\sin^2\theta' & (1-\varepsilon')\cos^2\theta' \\
-\cos^2\theta & \varepsilon - \cos^2\theta & \varepsilon' - \cos^2\theta & \sin^2\theta
\end{pmatrix}\,.
  \end{split}
 \end{equation}

\normalsize

Plugging the entries of \eqref{eq:4rewmatr} into \eqref{eq:rewardRnodisc} and replacing the $\mathsf{a}_{ij}$'s from \eqref{eq:4Param4SIMPL} yields the reward $R_\pi(\tau)$ of the considered trajectory:
\begin{equation}\label{eq:reward4}
\begin{split}
R_\pi(\tau) \,&=\, (\textsf{n}_0 - \mathsf{c}_{01} - \mathsf{c}_{01'})(-\varepsilon\sin^2\theta\cos^2\theta - \varepsilon'\sin^2\theta) \\
&\quad + \mathsf{c}_{01}(\varepsilon - \varepsilon\sin^2\theta\cos^2\theta - \varepsilon'\sin^2\theta) \\
&\quad + \mathsf{c}_{01'}(\varepsilon' - \varepsilon\sin^2\theta\cos^2\theta - \varepsilon'\sin^2\theta) \\
&\quad + (\textsf{n}_1 - \mathsf{c}_{01})(-(1-\varepsilon)\sin^2\theta) \\
&\quad + \mathsf{c}_{01}((1-\varepsilon)\cos^2\theta) \\
&\quad + (\textsf{n}_{1'} - \mathsf{c}_{01'})(-(1-\varepsilon')\sin^2\theta') \\
&\quad + \mathsf{c}_{01'}((1-\varepsilon')\cos^2\theta') \\
&\quad + (\mathsf{c}_{01} + \mathsf{c}_{01'} - 1)(-\cos^2\theta) \\
&\quad + (\textsf{n}_2 - \mathsf{c}_{01} - \mathsf{c}_{01'})(\sin^2\theta) \\
&=\, \textsf{n}_0 (-\varepsilon\sin^2\theta\cos^2\theta - \varepsilon'\sin^2\theta) + \textsf{n}_1 (-(1-\varepsilon)\sin^2\theta) \\
&\quad + \textsf{n}_{1'} (-(1-\varepsilon')\sin^2\theta') + \textsf{n}_2 (\sin^2\theta) \,+\, \cos^2\theta \\
&=:\,R_{\varepsilon,\varepsilon'}(\textsf{n}_0,\textsf{n}_1,\textsf{n}_{1'},\textsf{n}_2,\theta,\theta')\,.
\end{split}
\end{equation}
Remarkably, the linear dependence on the branching fluxes $\mathsf{c}_{01}$ and $\mathsf{c}_{01'}$ cancels out exactly (e.g., the energy cost of entering, traversing, and exiting the upper branch is balanced by the feedback loop cost).
The final result depends solely on the occupation numbers.

\medskip

  \begin{figure}[t!]
    \centering
   \includegraphics[height=6.4cm]{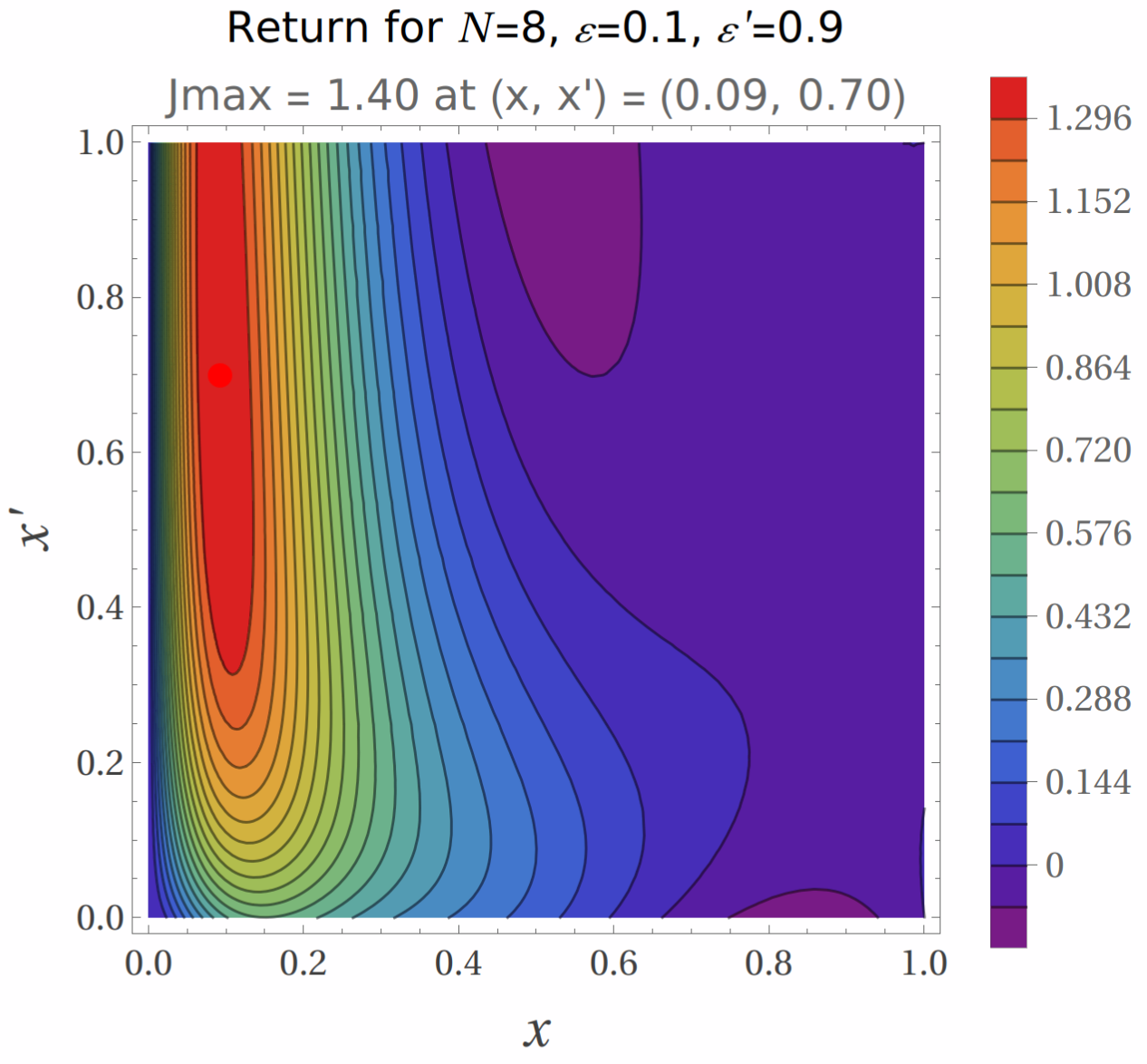}\quad
        \includegraphics[height=6.4cm]{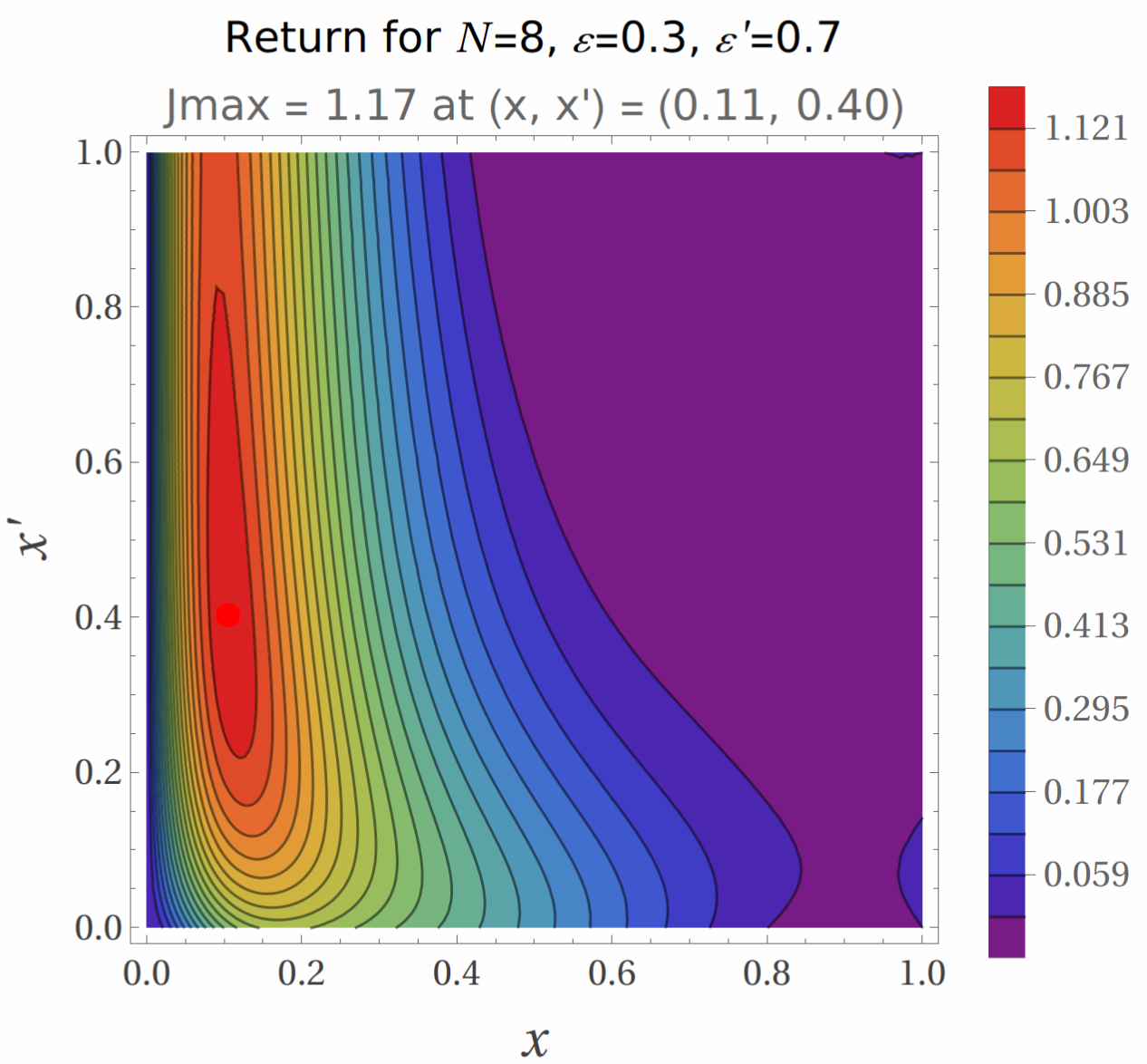}

 \medskip

   \includegraphics[height=6.4cm]{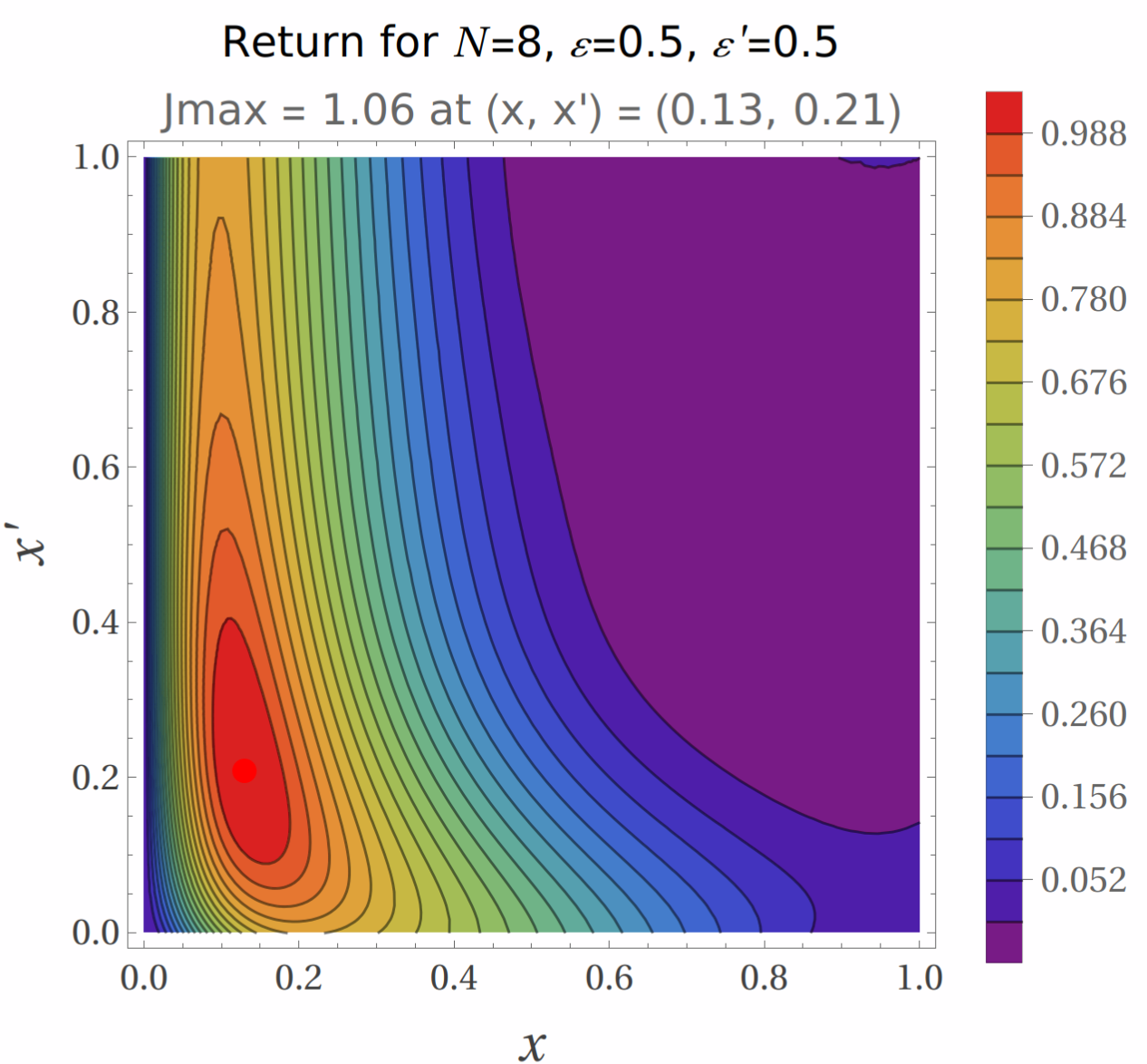}\quad
        \includegraphics[height=6.4cm]{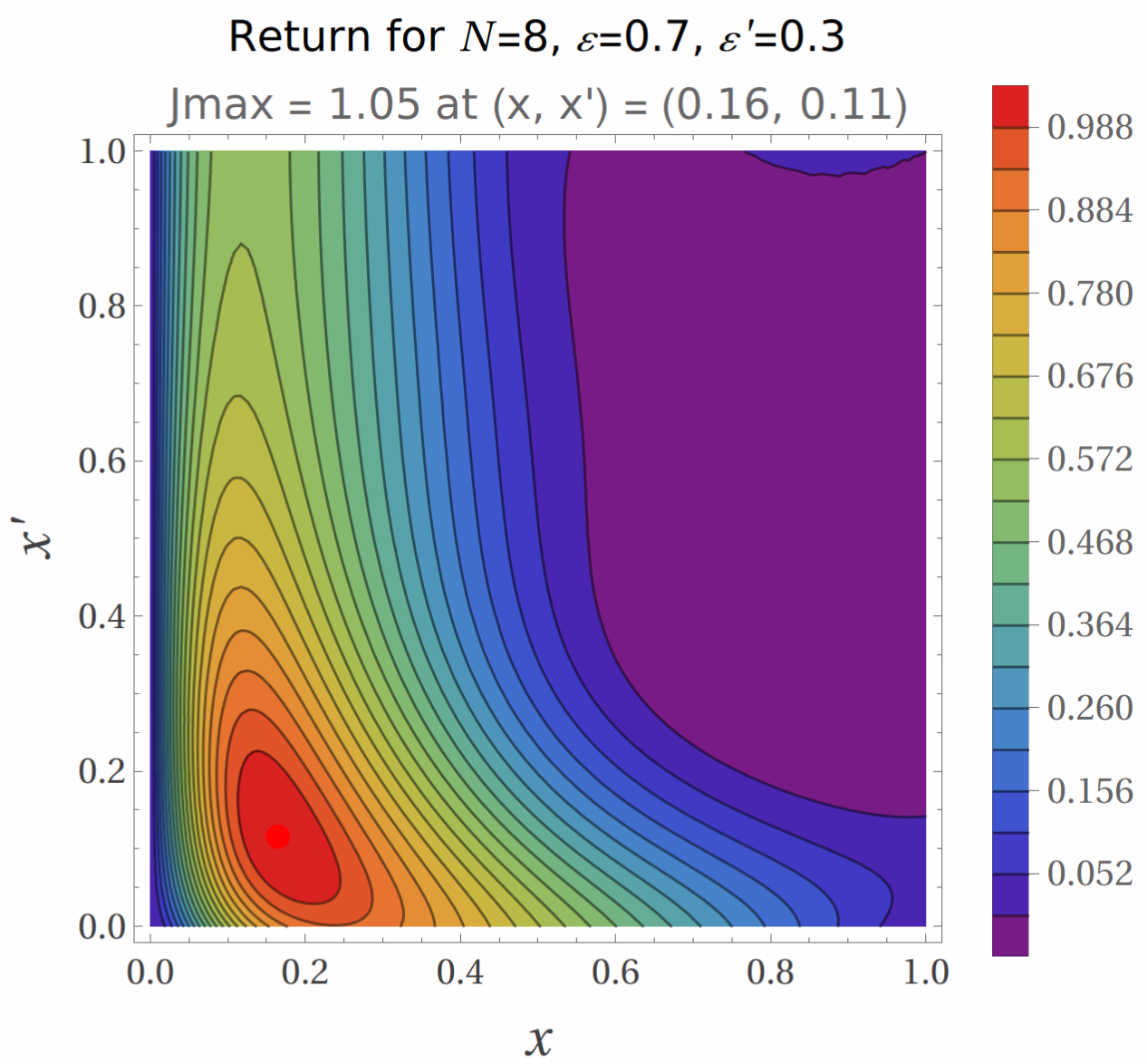
        }
    \caption{Numerical optimisation of the expected return \eqref{eq:Jpidi4long} on a horizon of $N=8$ long trajectories. The plots show the migration of the optimal policy inside the control space $[0,1]\times[0,1]$ as the energy penalties $(\varepsilon, \varepsilon')$ are varied.}
    \label{fig:Jpi-4level}
\end{figure}

  \begin{figure}[t!]
    \centering
   \includegraphics[height=6.8cm]{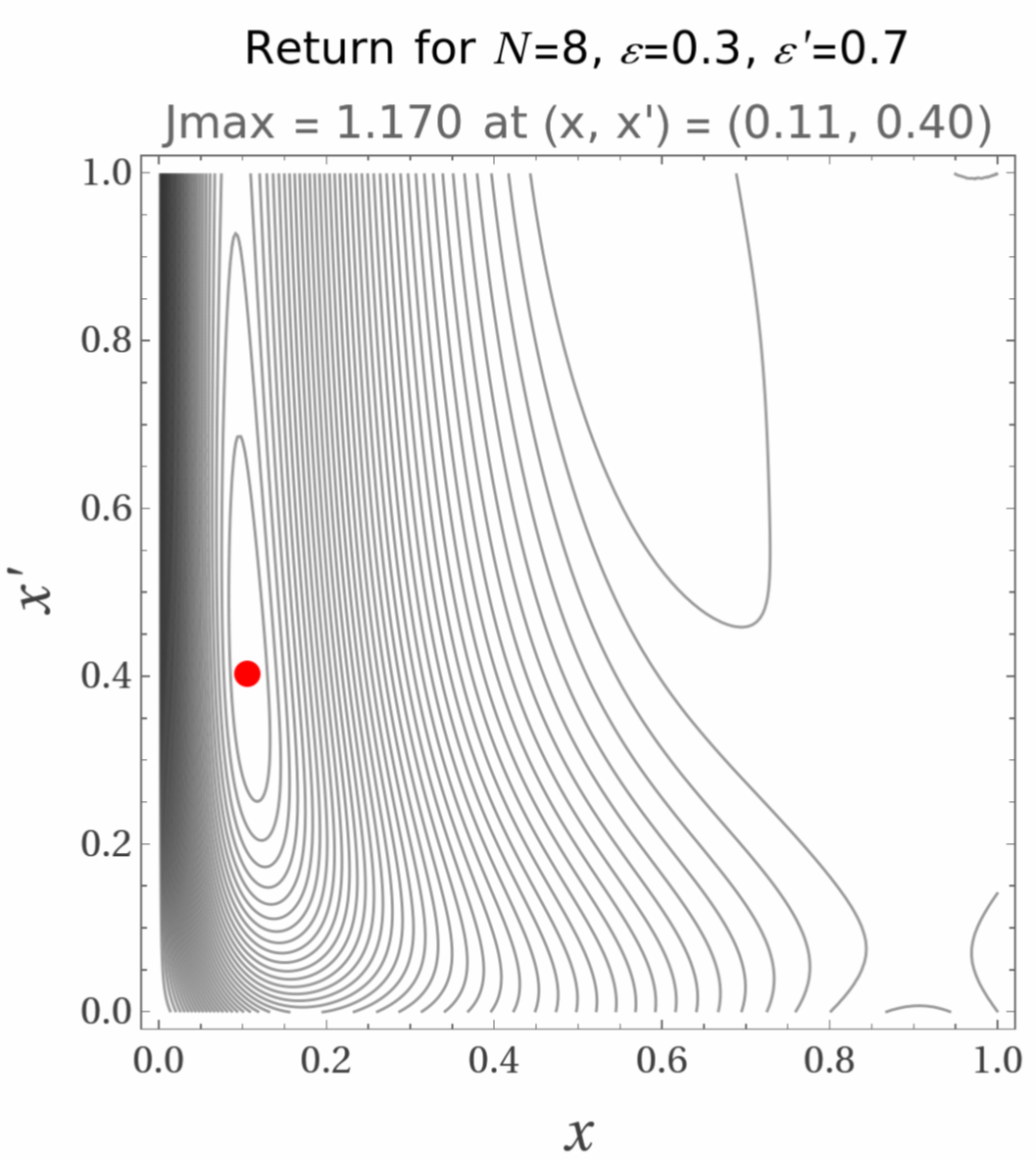}\quad
        \includegraphics[height=6.8cm]{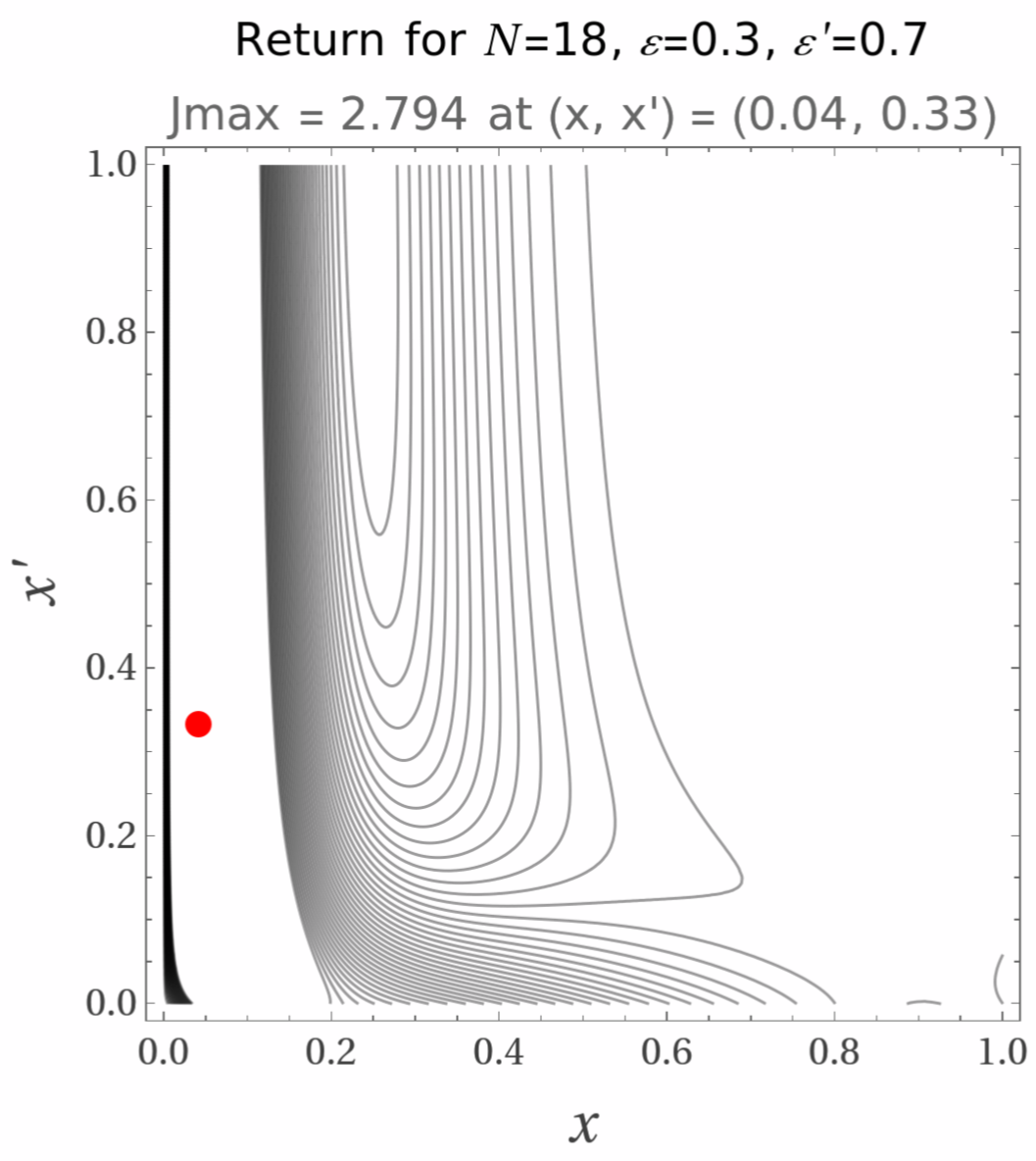}
\caption{Numerical evidence of the emergence of a plateau-like region in the profile of the return function \eqref{eq:Jpidi4long} around the point of optimal policy, as $N$ increases.}
    \label{fig:4levelPlateau}
\end{figure}

  \begin{figure}[t!]
    \centering
   \includegraphics[height=6.4cm]{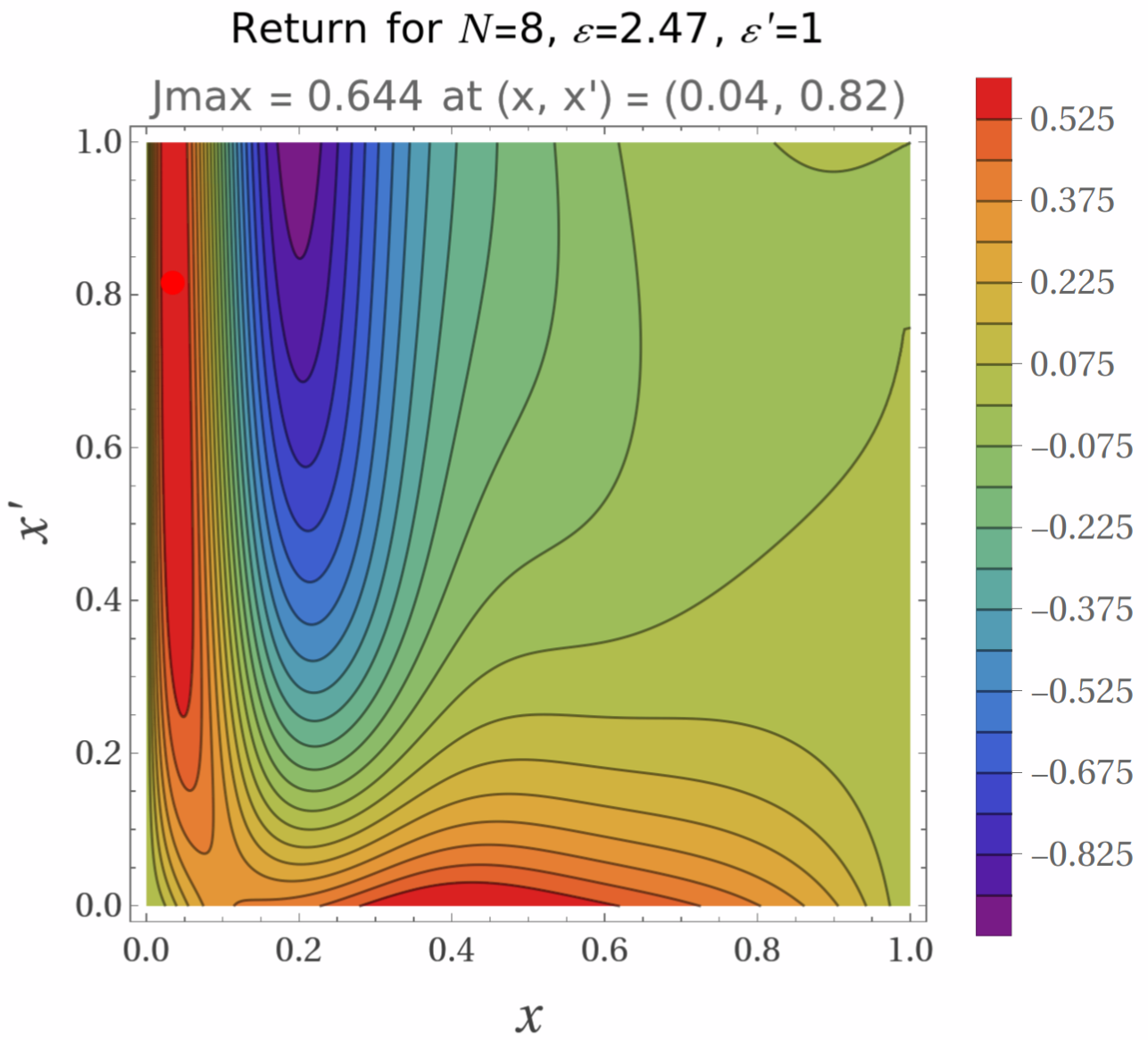}\quad
        \includegraphics[height=6.4cm]{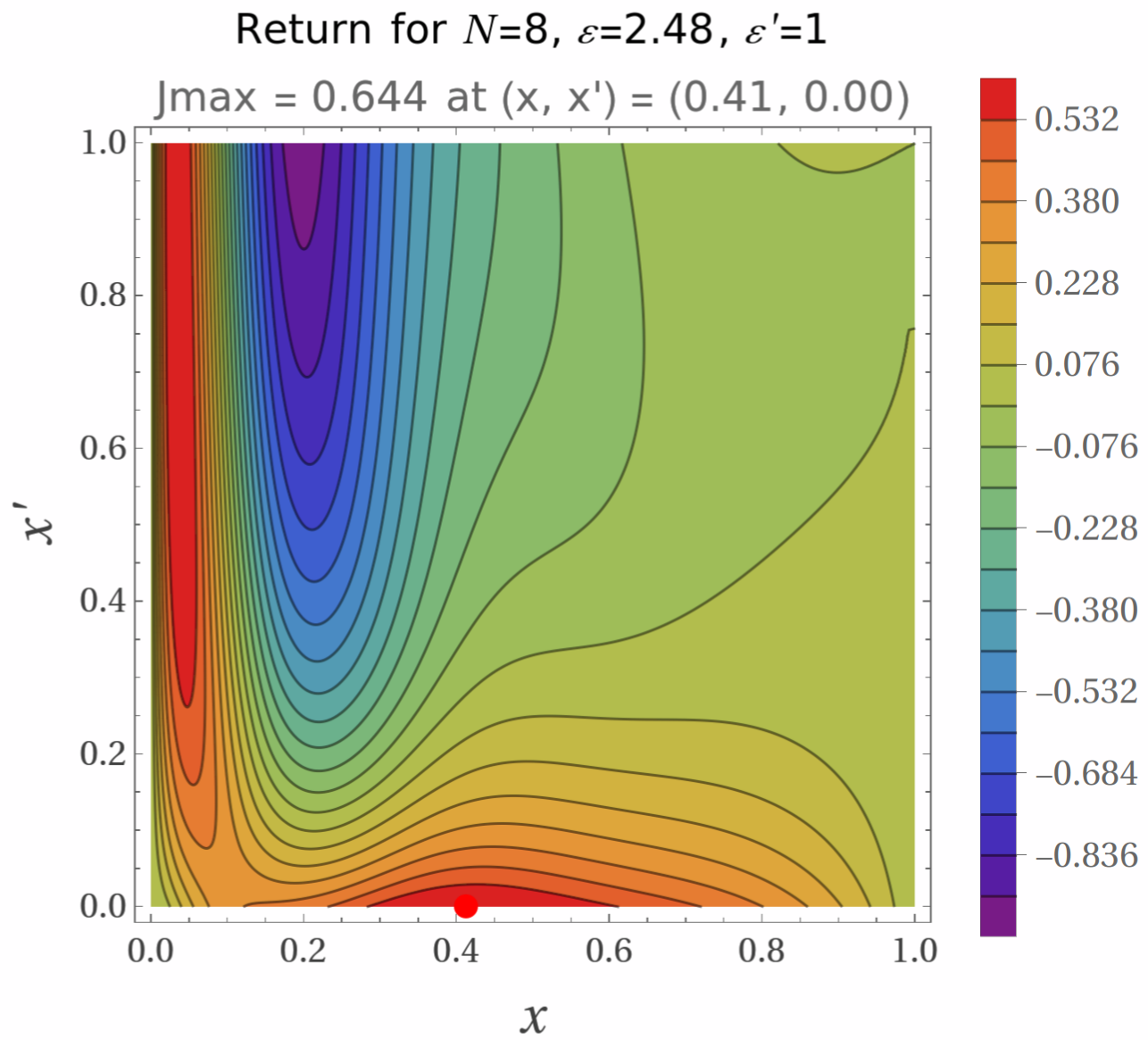}
\caption{Emergence of degenerate optimal policies for the expected return \eqref{eq:Jpidi4long} (horizon $N=8$). A crossover of the global maximum $(x_{\textrm{max}},x'_{\textrm{max}})$ between two distinct regions of the domain is observed when increasing $\varepsilon$ from $2.47$ to $2.48$. Since the return $J_{\varepsilon,1}$ is continuous in $\varepsilon$, this implies the existence of a critical value $\varepsilon^* \in (2.47, 2.48)$ where two distinct optimal policies coexist with $J_{\varepsilon^*,1}^{\textrm{max}}\approx 0.644$.}
    \label{fig:4levelDEG}
\end{figure}

\textbf{Expected return formula.} The expected return formula \eqref{eq:Jpidi} takes for this model the form
\begin{equation}\label{eq:Jpidi4}
\begin{split}
  J(\pi) \:= &\sum_{ \substack{ \mathsf{n}_0, \mathsf{n}_1, \mathsf{n}_{1'}, \mathsf{n}_2, \mathsf{c}_{01}, \mathsf{c}_{01'} \\ \textrm{with } \eqref{eq:njconstraints4UPDcorr}} } \: \mathcal{N}(N,\textsf{n}_0,\textsf{n}_1,\textsf{n}_{1'},\textsf{n}_2,\mathsf{c}_{01},\mathsf{c}_{01'})\: \times \\
  &\qquad\qquad\times R_{\varepsilon,\varepsilon'}(\textsf{n}_0,\textsf{n}_1,\textsf{n}_{1'},\textsf{n}_2,\theta,\theta')\:\times\: P(N,\textsf{n}_0,\textsf{n}_1,\textsf{n}_{1'},\textsf{n}_2,\mathsf{c}_{01},\mathsf{c}_{01'},\theta,\theta')\,.
\end{split}
\end{equation}
Substituting the closed-form multiplicity \eqref{eq:multiplicity4-closedform}, the trajectory probability \eqref{eq:ProbTraj4}, and the reward \eqref{eq:reward4} into \eqref{eq:Jpidi4} yields a fully analytic expression. Re-writing $J(\pi)\equiv J_{\varepsilon,\varepsilon'}(x,x')$ in terms of the new variables $x,x'\in[0,1]$ defined by
\begin{equation}
x\,:=\,\sin^2\theta\,,\qquad x'\,:=\,\sin^2\theta'\,,
\end{equation}
one obtains
\begin{equation}\label{eq:Jpidi4long}
\begin{split}
J_{\varepsilon,\varepsilon'}(x,x') \:=\: &\sum_{\textsf{n}_0=1}^{N-1} \sum_{\textsf{n}_2=1}^{N-\textsf{n}_0} \sum_{\textsf{n}_1=0}^{N+1-\textsf{n}_0-\textsf{n}_2} \sum_{\mathsf{c}_{01}=\mathbf{1}_{\{\textsf{n}_1\geqslant 1\}}}^{\min\{\textsf{n}_1,\textsf{n}_0,\textsf{n}_2\}} \sum_{\mathsf{c}_{01'}=\max\{\mathbf{1}_{\{\textsf{n}_{1'}\geqslant 1\}},\,1-\mathsf{c}_{01}\}}^{\min\{\textsf{n}_{1'},\textsf{n}_0-\mathsf{c}_{01},\textsf{n}_2-\mathsf{c}_{01}\}} \\
&\qquad\binom{\textsf{n}_0-1}{\mathsf{K}-1}\binom{\textsf{n}_2-1}{\mathsf{K}-1}\binom{\mathsf{K}}{\mathsf{c}_{01}}\binom{\textsf{n}_1-1}{\mathsf{c}_{01}-1}\binom{\textsf{n}_{1'}-1}{\mathsf{c}_{01'}-1} \\
&\qquad\times \Big[(1-x)+ \textsf{n}_2 x - \textsf{n}_1(1-\varepsilon)x  -\textsf{n}_{1'}(1-\varepsilon')x'\\
&\qquad\qquad -\textsf{n}_0(\varepsilon x(1-x) + \varepsilon' x) \Big] \\
&\qquad\times (1-x)^{2\textsf{n}_0 + \textsf{n}_1 + \textsf{n}_2 - 3\mathsf{c}_{01} - 3\mathsf{c}_{01'}} \, x^{3\mathsf{c}_{01} + 2\mathsf{c}_{01'} - 1} \\
&\qquad\times (1-x')^{\textsf{n}_{1'} - \mathsf{c}_{01'}} \, (x')^{\mathsf{c}_{01'}}\,, \\
& \textrm{with }\;\textsf{n}_{1'}=N+1-\textsf{n}_0-\textsf{n}_1-\textsf{n}_2\;\textrm{ and }\;\mathsf{K}=\mathsf{c}_{01}+\mathsf{c}_{01'}\,.
\end{split}
\end{equation}
The combinatorial multiplicity is embedded directly in the summand through the explicit product of binomial coefficients, with the convention $\binom{-1}{-1}=1$ inherited from \eqref{eq:multiplicity4-closedform}. The following is worth being remarked:
\begin{enumerate}
 \item[(1)] As with \eqref{eq:Jpidi2bis}, \eqref{eq:Jpidi2pm}, and \eqref{eq:Jpidi3} of the previous models, \eqref{eq:Jpidi4long} displays a substantial reduction of computational complexity with $N$: the sum scales as $O(N^5)$ (summing over $\textsf{n}_0,\textsf{n}_1,\textsf{n}_2,\mathsf{c}_{01},\mathsf{c}_{01'}$ with constraints), as compared to the general $O(4^N)$-scaling of the expression \eqref{eq:Jpidi}. The closed-form multiplicity \eqref{eq:multiplicity4-closedform} eliminates the need for any combinatorial enumeration of the trajectories, making the computation of $J_{\varepsilon,\varepsilon'}(x,x')$ fully analytical, in line with the qutrit case (Section~\ref{sec:Qtrits}, Remark~(2) after \eqref{eq:Jpidi3}).
 \item[(2)] Formula \eqref{eq:Jpidi4long} depends smoothly on the two control variables $(x,x')\in[0,1]^2$ and on the two energy parameters $(\varepsilon,\varepsilon')$, making the policy optimisation tractable numerically over a two-dimensional control space.
 \item[(3)] The reward kernel inside the square brackets in \eqref{eq:Jpidi4long} highlights the distinct roles of the three controls: $x$ couples to both upward transitions $|0\rangle\to|1\rangle$ and the feedback $|2\rangle\to|0\rangle$ (visible through the powers $x^{3\mathsf{c}_{01}+2\mathsf{c}_{01'}-1}$), while $x'$ couples exclusively to the upper branch $|1'\rangle$ (powers $(x')^{\mathsf{c}_{01'}}$). The reward structure similarly separates the energy costs: penalties proportional to occupancies $\textsf{n}_0$, $\textsf{n}_1$, $\textsf{n}_{1'}$ weighted by partial steps $\varepsilon$ and $\varepsilon'$, versus a gain proportional to $\textsf{n}_2$.
\end{enumerate}

\textbf{Numerical optimisation.} Figure \ref{fig:Jpi-4level} illustrates the migration of the optimal policy $(x_{\textrm{max}},x'_{\textrm{max}})$ inside $[0,1]\times[0,1]$ for the return $J_{\varepsilon,\varepsilon'}(x,x')$, as $\varepsilon$ increases from $0$ to $1$ and $\varepsilon'$ decreases from $1$ to $0$, at fixed $N$.

Numerics also show that, as the horizon length increases, the profile of $J_{\varepsilon,\varepsilon'}$ tends to develop local regions of plateau type around the optimal policy (Figure \ref{fig:4levelPlateau}). This indicates the emergence of an almost infinite local quasi-degeneracy of the optimal policy at large horizons.

 Most noticeably, beyond the standard regime of energies $0\leqslant\varepsilon,\varepsilon'\leqslant 1$, exploring the broader parameter space reveals critical phenomena akin to first-order phase transitions. Figure \ref{fig:4levelDEG} captures a crossover event where the global maximum discontinuously jumps between two distinct regions of the control space as $\varepsilon$ is tuned and, as a consequence, the existence of two distinct, widely separated optimal policies coexist with identical expected returns. This represents a genuine \emph{optimisation degeneracy}, where the agent faces a choice between two qualitatively different but quantitatively equal strategies.

\section{Complexity reduction}\label{sec:complexitycomparison}

The models explored so far, worked out both analytically and numerically in Sections \ref{sec:closequbit} through \ref{sec:4levels}, provide an instructive testing ground for the discussion of the first of the two structural features that this work set out to investigate (Section \ref{sec:intro}): the scaling with $N$ of the complexity of the expected return.

As will become clear in a moment, such a discussion is only made possible by the availability of the exact, closed-form computations carried out for each model -- which was indeed one of the very purposes of performing them. The mechanisms of the complexity reduction would have remained hidden in a black-box numerical evaluation.

For concreteness, the three-level model of Section \ref{sec:Qtrits} shall serve throughout as the main reference case. The subject is in fact considerably larger, and we aim at a systematic study in a subsequent analysis: the scope here is to highlight the evidence emerging from the instructive `toy-model' realisations of this work.

In all such models, the agent must optimise the return function $J(\pi)$ from \eqref{eq:Jpidi} over the set of allowed policies $\pi$.

The policies consist of certain choices of unitary transformations \eqref{eq:choosepol}, one for each state of the orthonormal basis \eqref{eq:statesS} of the $d$-dimensional system. If all unitaries are allowed, the optimisation is made over a $d^2$-real-parameter set (scenarios of Sections \ref{sec:closequbit} and \ref{sec:qubitantiper}). This number is reduced if a constraint on the type of policies is imposed (scenarios of Sections \ref{sec:Qtrits} and \ref{sec:4levels}). In either case, crucially, this count does not depend on $N$.

What is critical, then, when one wants to set up the model with longer and longer trajectories, is how the computation of each $J(\pi)$ scales with $N$.

The expressions \eqref{eq:probtraj}, \eqref{eq:rewardRnodisc}, and \eqref{eq:Jpidi} indicate that for a fixed policy $\pi$, the number of operations required to determine $J(\pi)$ nominally scales exponentially with $N$, since the set $\Gamma_{\pi,N}$ of length-$N$ trajectories $\tau$ in a $d$-level model (in the sense of \eqref{eq:statesS} above) has cardinality $d^{N-1}$ (out of the $N+1$ nodes, the two endpoints are fixed).

This count, however, neglects the specific structure of the quantum reinforcement models under consideration. Our analytic derivation of the return function, instead, exploits such crucial features.

The ensuing drop from the nominal exponential complexity to an actual power-law complexity is the combined effect of distinct structural mechanisms:
\begin{enumerate}
 \item[I.] a \emph{trajectory equivalence} (Section \ref{sec:reductionI}), operating for every admissible policy $\pi$;
 \item[II.] a \emph{sparsity of the transition graph} (Section \ref{sec:reductionII}), equally policy-independent;
 \item[III.] a \emph{spectral concentration at the optimum} (Section \ref{sec:reductionIII}), a mechanism of different nature, localised at the optimal policy.
\end{enumerate}

\subsection{Reduction I -- trajectory equivalence}\label{sec:reductionI}

The first key point is that, still at full level of generality, the elements of $\Gamma_{\pi,N}$ can be grouped into sub-classes, each of which is identified by an $N$-independent number of parameters:
\begin{itemize}
 \item the counts $\mathsf{n}^{(1)},\dots,\mathsf{n}^{(d)}$ of trajectory states for each of the considered energy level vectors $\{\Psi^{(1)},\dots,\Psi^{(d)}\}$ -- thus, $\mathsf{n}_0,\mathsf{n}_1,\mathsf{n}_2$ for the three-level system studied in Section \ref{sec:Qtrits}, or $\mathsf{n}_0,\mathsf{n}_1,\mathsf{n}_{1'},\mathsf{n}_2$ for the four-level system studied in Section \ref{sec:4levels};
 \item the inter-state transition counts $\mathsf{a}_{ij}$, $i,j\in\{1,\dots,d\}$.
\end{itemize}
Thus, $d+d^2$ parameters naturally label each trajectory sub-class in $\Gamma_{\pi,N}$.

The above completely general re-grouping of the trajectories of the collection $\Gamma_{\pi,N}$ is not a source of complexity reduction per se. Indeed, a priori, inside each sub-class labelled as above, further choices of trajectories could be possible, each contributing in a different way to the return function $J(\pi)$, so that exploring the whole $\Gamma_{\pi,N}$ still requires an exponentially-long-in-$N$ summation in \eqref{eq:Jpidi}.

Here the second crucial point enters: \emph{all the trajectories in each of the sub-classes par\-am\-e\-trised as above give instead exactly the same contribution to the computation of $J(\pi)$}. This fact is due to the structure of the trajectory probability \eqref{eq:probtraj} and the trajectory reward \eqref{eq:rewardRnodisc}: such quantities only depend on the labelling parameters identified above. In the three-level model this is read out explicitly from formulas \eqref{eq:ProbTraj3} and \eqref{eq:reward3}, where probability and reward of a trajectory only depend on the occupation numbers $\textsf{n}_0,\textsf{n}_1,\textsf{n}_2$ and on the transition counts $\mathsf{a}_{ij}$, the cardinality of each sub-class being given by the closed-form multiplicity \eqref{eq:3trajmult}.

As a consequence, the nominally exponential complexity in the computation of $J(\pi)$ drops down to an actual power-law complexity. For each of the parameters labelling each sub-class of trajectories there is a choice of at most $O(N)$ options. This brings the number of distinct sub-classes to $O(N^{d+d^2})$, a power law in $N$. All that remains is to compute the contribution `reward $\times$ probability' of a generic trajectory for each sub-class, and to multiply it by the cardinality $\mathcal{N}$ of that sub-class, be the magnitude of such cardinality polynomial or exponential in $N$. With $O(N^{d+d^2})$ of such operations, $J(\pi)$ is finally computed with power-law complexity.

The exponent $d+d^2$ established by this first mechanism is, however, an over-estimate: a second, independent mechanism lowers it substantially.

\subsection{Reduction II -- sparsity of the transition graph}\label{sec:reductionII}

In practice, owing to the existence of forbidden (zero-probability) transitions, which in turn is a consequence of the choice of the allowed unitaries and of the structure of the probability functional, a lower number $\mathcal{I}_{\textrm{par}}<d+d^2$ of independent parameters suffices to identify each trajectory class.

In the three-level model, the zero entries of the transition probability matrix \eqref{eq:3transP} act as super-selection rules enforcing the three forbidden transitions \eqref{eq:3additionalConstr}; together with the three outgoing constraints \eqref{eq:constr3-2} and the two independent net-flow (Kirchhoff) equations \eqref{eq:netflow3} -- the incidence matrix of the connected transition graph $\mathcal{G}_3$ having rank $3-1=2$ \cite[Theorem~2.3]{Bapat2014} (Figure \ref{fig:graph-qutrits}) -- this makes $3+2+3=8$ independent conditions on the $9$ transition counts $\mathsf{a}_{ij}$, leaving exactly one residual degree of freedom \eqref{eq:c20def}. Explicitly:
\begin{itemize}
 \item the parameters $\mathsf{n}_0,\mathsf{n}_{2},\mathsf{c}$ ($\mathcal{I}_{\textrm{par}}=3$) for the three-level model (Section \ref{sec:Qtrits});
 \item the parameters $\mathsf{n}_0,\mathsf{n}_1,\mathsf{n}_{2},\mathsf{c}_{01},\mathsf{c}_{01'}$ ($\mathcal{I}_{\textrm{par}}=5$) for the four-level model (Section \ref{sec:4levels}).
\end{itemize}
Correspondingly, the number of distinct sub-classes drops from $O(N^{d+d^2})$ to $O(N^{\mathcal{I}_{\textrm{par}}})$ -- see \eqref{eq:3constraintsn} or \eqref{eq:njconstraints4UPD} above -- and so does the complexity of the computation of $J(\pi)$: thus, $O(N^3)$ for the three-level model and $O(N^5)$ for the four-level model, as compared to the nominal $O(3^N)$- and $O(4^N)$-scaling of the expression \eqref{eq:Jpidi}.

Without altering the above conclusion, two technical points that emerged in the discussion are worth being highlighted in this regard. First, the conservation constraint $\sum_{j=1}^d\mathsf{n}^{(j)}=N+1$, and the constraints on forbidden transitions $\mathsf{a}_{ij}$ driven by the special choices of the unitary actions and of the expression of the trajectory probability, all contribute to reduce further the actual range of admissible values for the $\mathcal{I}_{\textrm{par}}$ independent parameters labelling each sub-class. This range can be considerably smaller than the nominal order $N$. Thus, $N^{\mathcal{I}_{\textrm{par}}}$ is an actual over-estimate of the final power-law complexity. This occurrence was discussed in Sections \ref{sec:modelQtrits} and \ref{sec:4levels}.

Furthermore, whereas all $\mathcal{I}_{\textrm{par}}$ independent parameters are needed in general to identify a sub-class in the whole trajectory landscape, for certain policies a strictly smaller number of parameters suffices, thereby resulting in a further reduction of the power-law complexity. This was observed in the last example of Section \ref{sec:2levelclosed}.

\subsection{A comparison between brute-force numerics and analytics}\label{sec:numbruteanal}

The three-level model provides a robust and effective testing ground to quantify the combined effect of the two reductions (Sections \ref{sec:reductionI} and \ref{sec:reductionII}), from exponential to power-law, on the evaluation of the return function $J_\varepsilon(x)$, for concreteness in the scenario of policies of common rotations (Section \ref{sec:optim3commonrot}).

\begin{figure}[t]
\centering
\includegraphics[width=0.8\textwidth]{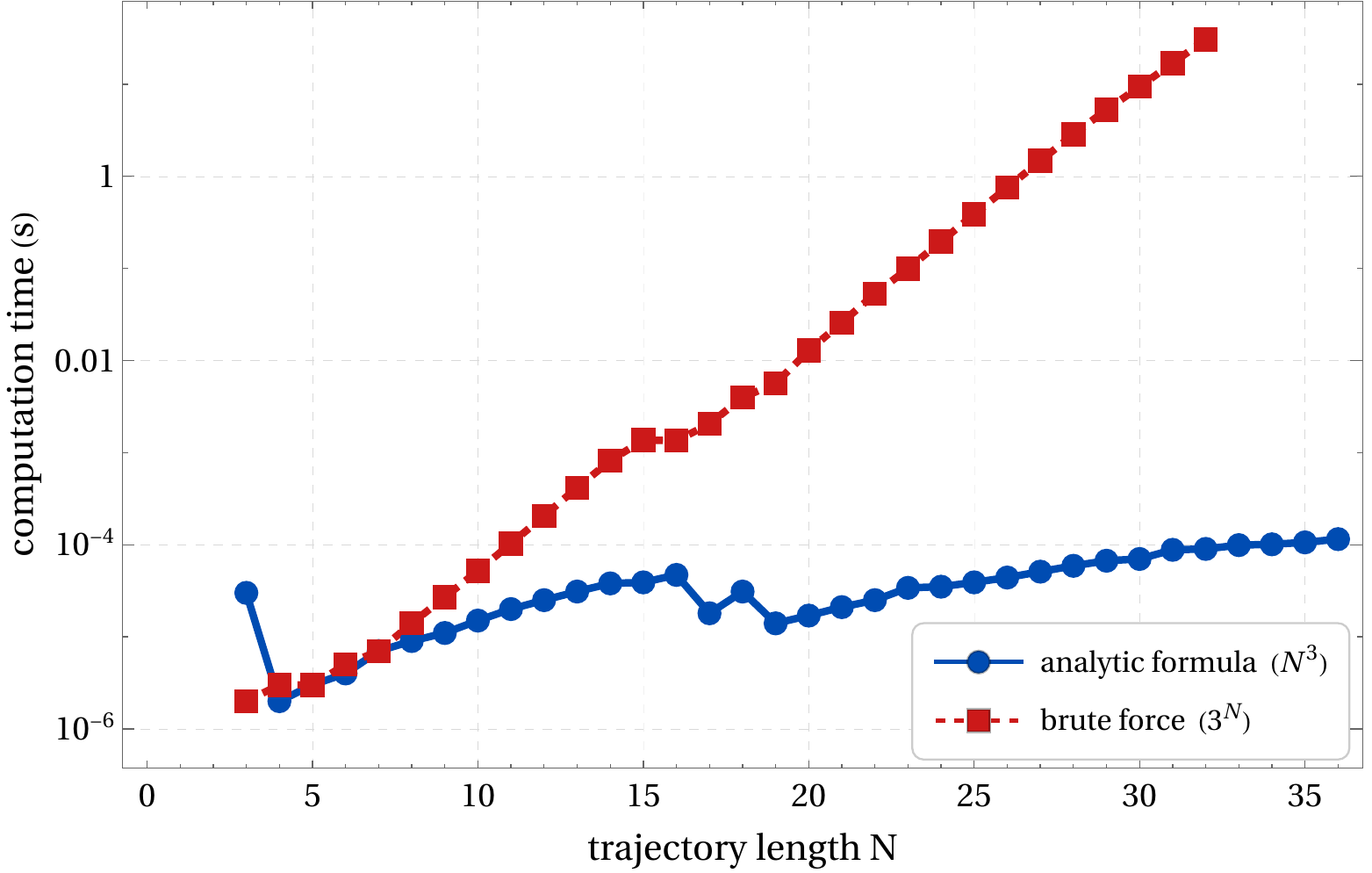}
\caption{Comparison of computational time required to evaluate the expected return $J_\varepsilon(x)$ as a function of trajectory length $N$. The `Brute Force' approach (red dashed line) scales exponentially as $O(3^N)$, while the analytic formula (blue solid line) scales polynomially as $O(N^3)$.}
\label{fig:complexity_comparison}
\end{figure}

The `brute-force' evaluation involves summing over all possible intermediate state trajectories. Based on the general definitions \eqref{eq:probtraj}, \eqref{eq:rewardRnodisc}, and \eqref{eq:Jpidi} of Section \ref{sec:generalQRLmodel},
the expected return of the model under consideration is given by
 \begin{equation}\label{eq:JpidiBF}
 \begin{split}
J_{\varepsilon}(x)\,&=\sum_{ \substack{ |\psi_0\rangle = |0\rangle,\; |\psi_N\rangle=|2\rangle, \\ |\psi_1\rangle\in\{|0\rangle,|1\rangle,|2\rangle\}, \\ .... \\ .... \\ |\psi_{N-1}\rangle\in\{|0\rangle,|1\rangle,|2\rangle\}   }  } \bigg(\sum_{i=0}^{N-1}\langle\psi_{i+1}|H\psi_{i+1}\rangle-\big\langle\pi(|\psi_i\rangle)\psi_i|H\pi(|\psi_i\rangle)\psi_i\big\rangle\bigg) \times \\
& \qquad \qquad\qquad\qquad \qquad \times \bigg( \prod_{k=0}^{N-1} \big|\big\langle\psi_{k+1}\big|\pi(|\psi_k\rangle)\psi_k\big\rangle\big|^2\bigg),
 \end{split}
   \end{equation}
   where
 \begin{equation}\label{eq:Hamilt3BF}
   \begin{split}
      H\,=\begin{pmatrix}
         0 & 0 & 0 \\
         0 & \varepsilon & 0 \\
         0 & 0 & 1
        \end{pmatrix},\quad |0\rangle=\begin{pmatrix} 1 \\ 0 \\ 0 \end{pmatrix}\,,\quad
   |1\rangle=\begin{pmatrix} 0 \\ 1 \\ 0 \end{pmatrix}\,,\quad
   |2\rangle=\begin{pmatrix} 0 \\ 0 \\ 1 \end{pmatrix}\,,
  \end{split}
 \end{equation}
 and with matrices $\pi_0\equiv\pi(|0\rangle)$, $\pi_1\equiv\pi(|1\rangle)$, $\pi_2\equiv\pi(|2\rangle)$ given by
 \begin{equation}
  \hspace{-1cm}\pi_0(\theta)\,=
  \begin{pmatrix}
   \cos\theta & -\sin\theta & 0 \\
   \sin\theta & \cos\theta & 0 \\
   0 & 0 & 1
  \end{pmatrix},\;
  \pi_1(\theta)\,=
  \begin{pmatrix}
   1 & 0 & 0 \\
   0 & \cos\theta & -\sin\theta \\
   0 & \sin\theta & \cos\theta
  \end{pmatrix},\;
 \pi_2(\theta)\,=
  \begin{pmatrix}
   \cos\theta & 0 & -\sin\theta \\
   0 & 1 & 0 \\
   \sin\theta & 0 & \cos\theta
  \end{pmatrix}
 \end{equation}
 (see \eqref{eq:states3}-\eqref{eq:unitaries3-1} above).
 Such matrices are $\theta$-dependent, with
 \begin{equation}
  \theta\,=\,\arcsin\sqrt{x}\,.
 \end{equation}

 The sum \eqref{eq:JpidiBF} runs over the $3^{N-1}$ possible configurations of the intermediate states, and is parametrised by the intermediate energy level $\varepsilon$.

This computational task is to be compared with the closed-form analytic expression derived in \eqref{eq:Jpi3simplified-3}, which reduces the summation complexity from $O(3^N)$ to $O(N^3)$.

The computational cost of evaluating the expected return $J_\varepsilon(x)$ using both methods is illustrated in Figure \ref{fig:complexity_comparison}. The numerical evidence strikingly confirms the theoretical prediction.

It is also worth noting a structural difference regarding the intermediate energy level $\varepsilon$. As proven in Section \ref{sec:optim3commonrot}, the exact expected return $J(x)$ is independent of $\varepsilon$ due to the perfect cancellation of symmetric trajectories. The analytic formula \eqref{eq:Jpi3simplified-3} respects this symmetry manifestly, as $\varepsilon$ does not appear in the expression.

Conversely, the brute-force summation \eqref{eq:JpidiBF} blindly incorporates $\varepsilon$ into the energy cost of every individual step. These contributions theoretically sum to zero, but the numerical algorithm performs $O(3^N)$ floating-point operations involving $\varepsilon$. This not only obscures the physical independence of the result from the intermediate energy level but also introduces potential for accumulated rounding errors (catastrophic cancellation), which could lead a purely numerical optimizer to detect false dependencies on $\varepsilon$ that are artifacts of finite precision rather than physical reality.

\subsection{Reduction III -- spectral concentration at the optimum}\label{sec:reductionIII}

The two mechanisms discussed so far lower the cost of \emph{evaluating} $J(\pi)$, and operate identically for every admissible policy $\pi$. The combined numerical and analytical evidence of the preceding Sections points to a third phenomenon, of different nature: it concerns the \emph{optimisation} itself, and indicates that at the optimal policy the value of $J(\pi)$ is asymptotically supported on a structurally identifiable minority of the $O(N^{\mathcal{I}_{\textrm{par}}})$ trajectory sub-classes -- those of \emph{polynomial} cardinality -- whereas the exponentially populated sub-classes, which account for almost the totality of $\Gamma_{\pi,N}$, contribute evanescently.

\medskip

\textbf{Numerical evidence and cardinality regimes.} For concreteness, consider the three-level model with policies of multiple independent rotations (Section \ref{sec:optim3indeprot}) in the regime $\varepsilon\in(\frac{1}{2},1)$. The numerical optimisation of \eqref{eq:Jpi3threethetas2} (Figure \ref{fig:Jpi-3level3VAR}, left panel) locates the optimal controls at $x^*,z^*\to 0$ as $N$ grows, with a decay compatible with a $O(N^{-1})$ scaling, and at $y^*\to 1$, the probability weight collapsing onto the sub-classes with $\textsf{n}_1=1+\mathsf{c}$ (minimal dwell at $|1\rangle$, Section \ref{sec:optim3indeprot}). On the other hand, the cardinality of the sub-classes contributing to \eqref{eq:Jpi3threethetas2} ranges from polynomial (\eqref{eq:cardA}-\eqref{eq:cardB}) to exponential \eqref{eq:cardexp} in $N$. The natural question is then which sub-classes the optimisation is actually sensitive to.

\medskip

\textbf{Suppression of the exponentially populated sub-classes.} The following estimate is rigorous and uniform over the policy region
\begin{equation}\label{eq:softregion}
 x\,\leqslant\,\kappa N^{-1}\,,\qquad z\,\leqslant\,\kappa N^{-1}
\end{equation}
for fixed $\kappa>0$ (a region which, according to the numerical evidence above, contains the optimal policy for all sufficiently large $N$). Fix $\gamma\in(0,1)$ and consider the aggregate contribution to \eqref{eq:Jpi3threethetas2} of all the sub-classes with $\mathsf{c}\geqslant\gamma N$. Each trajectory probability contains the factor $x^{1+\mathsf{c}}z^{\mathsf{c}}\leqslant(\kappa/N)^{2\mathsf{c}+1}$, all remaining probability factors being bounded by one; the reward bracket in \eqref{eq:Jpi3threethetas2} is bounded in absolute value by $N+1+2\kappa$; each of the three binomials in \eqref{eq:3trajmult} is bounded by $2^{\textsf{n}_j-1}$, whence $\mathcal{N}\leqslant 2^{N-2}$; and there are at most $N^3$ sub-classes altogether. Therefore,
\begin{equation}\label{eq:expsuppression}
\begin{split}
  \left|\sum_{\substack{(\textsf{n}_0,\textsf{n}_2,\mathsf{c}) \textrm{ with } \eqref{eq:3constraintsn} \\ \mathsf{c}\,\geqslant\,\gamma N}}\mathcal{N}\,R\,P\,\right|\;&\leqslant\; N^3\,(N+1+2\kappa)\;2^{N-2}\,\Big(\frac{\kappa}{N}\Big)^{2\gamma N+1} \\
  &=\;e^{\,N\,(\ln 2\,-\,2\gamma\ln N)\,(1+o(1))}\,,
\end{split}
\end{equation}
which vanishes super-exponentially fast as $N\to\infty$. Thus, on the whole region \eqref{eq:softregion} the entropic growth $e^{Nr(\alpha,\beta,\gamma)}$ of the cardinality \eqref{eq:cardexp} is over-compensated by the decay $N^{-(2\mathsf{c}+1)}$ of the per-trajectory probability: the exponentially populated sub-classes are effectively silent.

\medskip

\textbf{Survival of the polynomially populated sub-classes.} Conversely, consider, as suggested by the numerics, the ansatz
\begin{equation}\label{eq:softansatz}
 x\,=\,\frac{\xi}{N}\,,\qquad y\,=\,1\,,\qquad z\,=\,\frac{\zeta}{N}\qquad (\xi,\zeta>0 \textrm{ fixed})\,,
\end{equation}
under which the weight $(1-y)^{\textsf{n}_1-1-\mathsf{c}}\,y^{1+\mathsf{c}}$ selects the sub-classes with $\textsf{n}_1=1+\mathsf{c}$, and take the family of sub-classes with $N$-independent $\mathsf{c}$ and
\begin{equation}\label{eq:survfamily}
 \textsf{n}_0\,=\,\alpha N\,,\qquad \textsf{n}_2\,=\,\beta N\,,\qquad \alpha,\beta>0\,,\quad \alpha+\beta=1+O(N^{-1})\,.
\end{equation}
The minimal-dwell selection trivialises the middle binomial of \eqref{eq:3trajmult} (boundary case of \eqref{eq:polysubclassB}), lowering the cardinality from $O(N^{3\mathsf{c}})$ to
\begin{equation}\label{eq:survcard}
 \mathcal{N}(N,\textsf{n}_0,\textsf{n}_2,\mathsf{c})\,=\,\binom{\textsf{n}_0-1}{\mathsf{c}}\binom{\mathsf{c}}{\mathsf{c}}\binom{\textsf{n}_2-1}{\mathsf{c}}\,=\,\frac{(\alpha\beta)^{\mathsf{c}}}{(\mathsf{c}!)^2}\,N^{2\mathsf{c}}\,\big(1+O(N^{-1})\big)\,.
\end{equation}
Under \eqref{eq:softansatz} the trajectory probability of \eqref{eq:Jpi3threethetas2} evaluates to
\begin{equation}\label{eq:survprob}
 P\,=\,e^{-\alpha\xi-\beta\zeta}\;\xi^{1+\mathsf{c}}\zeta^{\mathsf{c}}\;N^{-(2\mathsf{c}+1)}\,\big(1+o(1)\big)\,,
\end{equation}
whence
\begin{equation}\label{eq:survNP}
 \mathcal{N}\cdot P\,=\,\frac{(\alpha\beta\,\xi\zeta)^{\mathsf{c}}\,\xi}{(\mathsf{c}!)^2}\;e^{-\alpha\xi-\beta\zeta}\;\frac{1}{N}\,\big(1+o(1)\big)\,,
\end{equation}
while the reward bracket in \eqref{eq:Jpi3threethetas2} amounts to $1-\varepsilon\alpha\xi-(1-\varepsilon)(1+\mathsf{c})+\beta\zeta+O(N^{-1})$, an $O(1)$ quantity. Each such sub-class therefore contributes $O(N^{-1})$ to the expected return, and the $O(N)$ admissible choices of $\textsf{n}_0$ within the family produce an aggregate $O(1)$ contribution, surviving in the limit $N\to\infty$ -- consistently with the bounded values of $J$ observed numerically. The factorial damping $(\mathsf{c}!)^{-2}$ in \eqref{eq:survNP} further indicates that, within the surviving families, those with small $\mathsf{c}$ dominate.

Two caveats are in order. First, the above is a \emph{self-consistency} argument, not a dominance theorem: the location \eqref{eq:softansatz} of the optimum is assumed from the numerical evidence, and the comparison between the two classes of contributions is performed at that point, not over the whole control domain $[0,1]^3$. Second, maximising the probability factor family by family is informative only for the families that actually dominate the sum, and produces spurious predictions otherwise.

\smallskip

The above evidence is summarised in the following

\smallskip

\noindent\textbf{Conjecture [Spectral concentration].}
\emph{In the limit of large horizon $N$, the expected return at the optimal policy is asymptotically exhausted by the contribution of the trajectory sub-classes of polynomial-in-$N$ cardinality: the exponentially populated sub-classes, despite comprising almost the totality of $\Gamma_{\pi,N}$, contribute evanescently. Equivalently, only a polynomially small subset of the $O(N^{\mathcal{I}_{\textrm{par}}})$ sub-classes effectively supports the optimisation.}

\medskip

\textbf{Beyond the three-level model.} The phenomenon is not specific to the three-level model. In the anti-periodic qubit chain of Section \ref{sec:qubitantiper} the optimisation saturates $x_+^{\max}=1$ exactly, which forces the trajectory probability to vanish unless $\textsf{p}=\mathsf{c}$; the multiplicity $\binom{\textsf{p}-1}{\mathsf{c}-1}\binom{N-\textsf{p}}{\mathsf{c}-1}$ of the surviving sub-classes collapses accordingly to $\binom{N-\textsf{p}}{\textsf{p}-1}$. In terms of the residual dwell $\textsf{d}=N-2\textsf{p}+1$ at $|-\rangle$, the latter equals $\binom{\textsf{p}-1+\textsf{d}}{\textsf{d}}$: a polynomial of degree $\textsf{d}$ in $N$ when $\textsf{d}$ is kept $N$-independent (whence $\textsf{p}\approx\frac{N}{2}$), and exponentially large for extensive $\textsf{d}$. Repeating the previous two-step evaluation at the numerically observed optimum ($x_+=1$ and $1-x_-=O(N^{-1})$, Figure \ref{fig:Jpi-2levelANTIPER}), the probability factor $x_-^{\textsf{p}-1}(1-x_-)^{\textsf{d}}$ suppresses the exponentially populated sub-classes super-exponentially, while assigning each fixed-$\textsf{d}$ sub-class an $O(1)$ weight. Notably, the surviving sub-classes have here extensively many sign changes ($\mathsf{c}=\textsf{p}\approx\frac{N}{2}$): the relevant invariant of the phenomenon is the \emph{cardinality} of the sub-class, not the scaling of its transition counts.

If confirmed in wider generality, this mechanism would constitute a further, `spectral' reduction of complexity, complementary to Reductions I and II. A proof of dominance, requiring the comparison of the competing contributions over the whole control domain rather than at the numerically identified optimum, and the systematic analysis of the remaining models of this work, are deferred to a subsequent investigation.

\section{Degeneracy of optimal policies}\label{sec:deg}

The models explored in Sections \ref{sec:closequbit} through \ref{sec:4levels} provide equally instructive evidence on the second of the two structural features that this work set out to investigate (Section \ref{sec:intro}): the possible degeneracy of the optimal policies. Also in this respect the subject is considerably vaster than the present scope, and we intend to proceed with a systematic analysis in a separate work: the goal here is to organise the evidence emerging from the models under consideration.

In reinforcement learning it is of both theoretical and practical relevance to be able to distinguish between the existence of a single optimal policy $\pi^*$ as in \eqref{eq:optpolpistar} and the possibility that a multiplicity of distinct policies produce the same absolute maximum of $J(\pi)$.
On the practical side, in the presence of a degeneracy of optimal policies, additional criteria might be introduced to further select optimisations that are more convenient. Moreover, it is crucial whether the numerical methods can efficiently detect multiple distinct absolute maxima of the return function.

 Naturally, what is referred to here as (non-)degeneracy is the (non-)uniqueness of optimal policies \emph{up to symmetries} or trivial redundancies in the RL problem's formulation. For instance, it was evident throughout the concrete scenarios of Sections \ref{sec:closequbit}-\ref{sec:4levels} that an obvious periodic optimality degeneracy in the policies' rotation angles $\theta$'s is removed by transforming to variables of the type $x=\sin^2\theta$.

 In this respect, various instructive non-trivial behaviours were detected.
\begin{itemize}
  \item Modulo angles symmetries, the optimal policy turns out to be unique in the considered qubits scenarios (Sections \ref{sec:closequbit} and \ref{sec:qubitantiper}), as well as in the qutrits scenario with common rotations (Section \ref{sec:optim3commonrot}), exhibiting a precise asymptotic behaviour with  $N$ (Figures \ref{fig:Jpi-2levelPER}, \ref{fig:Jpi-2levelANTIPER}, and \ref{fig:Jpi-3levelCOMMON}).
  \item In the four-level system (Section \ref{sec:4levels}), as the horizon length $N$ increases, the profile of the return function tends to develop local regions of plateau type around the optimal policy (Figure \ref{fig:4levelPlateau}). This indicates the emergence of a continuum of infinite local quasi-degeneracy of the optimal policy at large horizons, posing a challenge for gradient-based optimisers which may stall in these flat regions.
    \item In the four-level system, in suitable regimes of the energy levels $\varepsilon,\varepsilon'$ a genuine discrete degeneracy of the optimal policy emerges (Figure \ref{fig:4levelDEG}), where the agent faces a choice between two qualitatively different but quantitatively equal strategies.
\end{itemize}

These types of actually occurring degeneracies are a characteristic of QRL that numerical methods for more realistic and complex models must account for to avoid misidentifying the global optimum or failing to converge due to landscape flatness.

\section{Concluding remarks}\label{sec:conclusions}

Let us finally revisit the class of QRL models presented in Section \ref{sec:generalQRLmodel}, in view of the insights from the simplified scenarios explored in Sections \ref{sec:closequbit} through \ref{sec:4levels}.

It is worth stressing that the rigorous analytical treatment
of the models under consideration is not intended to advocate generally for analytics
over numerics; clearly, closed-form derivations become more demanding as the complexity
of the underlying system increases, and the present work is deliberately focused on a
class of models chosen for their tractability and structural transparency.

The perspective advanced here is, rather, that structural insight -- of the kind
developed in this paper -- is a necessary complement to numerical approaches. Without
it, quantum reinforcement learning protocols proceed by black-box evaluation and
maximisation of the expected return, even when sophisticated ad hoc numerical methods
are employed.

The analytical treatment of the present models has concretely demonstrated three
specific vulnerabilities of purely numerical approaches: they may obscure structural
features (such as the independence of $J(\pi)$ from intermediate energy levels
$\varepsilon$, established in Section~\ref{sec:Qtrits}); they introduce unnecessarily
high computational complexity by failing to exploit the trajectory equivalence classes
identified in Section~\ref{sec:complexitycomparison}; and they may fail to detect or
correctly characterise the degeneracy of optimal strategies discussed in
Section~\ref{sec:deg}. These findings suggest concrete diagnostics that should inform
the design of numerical solvers for more realistic and complex QRL models.


\def\cprime{$'$}

\end{document}